\documentclass[10pt,twoside]{article}
\usepackage{graphicx}
\usepackage{amsmath}
\usepackage{amssymb}
\usepackage{amsthm}
\usepackage{hyperref}
\usepackage{subcaption}
\usepackage{cite}
\usepackage[capitalize]{cleveref}

\usepackage{tikz}
\usetikzlibrary{matrix,arrows,decorations.pathmorphing,decorations.markings}
\usepackage{tikz-cd}

\def\YEAR{\year}\newcount\VOL\VOL=\YEAR\advance\VOL by-1995
\def\firstpage{1}\def\lastpage{1000}
\def\received{}\def\revised{}
\def\communicated{}

\makeatletter
\def\magnification{\afterassignment\m@g\count@}
\def\m@g{\mag=\count@\hsize6.5truein\vsize8.9truein\dimen\footins8truein}
\makeatother

\font\eightrm=cmr8
\font\caps=cmcsc10                    
\font\Caps=cmcsc10 scaled \magstep1   

\makeatletter
\@specialpagefalse\headheight=8.5pt
\def\DocMath{}
\renewcommand{\@evenhead}{%
    \ifnum\thepage>\lastpage\rlap{\thepage}\hfill%
    \else\rlap{\thepage}\slshape\leftmark\hfill{\caps\SAuthor}\hfill\fi}%
\renewcommand{\@oddhead}{%
    \ifnum\thepage=\firstpage{\DocMath\hfill\llap{\thepage}}%
    \else{\slshape\rightmark}\hfill{\caps\STitle}\hfill\llap{\thepage}\fi}%
\makeatother

\def\TSkip{\bigskip}
\newbox\TheTitle{\obeylines\gdef\GetTitle #1
\ShortTitle  #2
\SubTitle    #3
\Author      #4
\ShortAuthor #5
\EndTitle
{\setbox\TheTitle=\vbox{\baselineskip=20pt\let\par=\cr\obeylines%
\halign{\centerline{\Caps##}\cr\noalign{\medskip}\cr#1\cr}}%
	\copy\TheTitle\TSkip\TSkip%
\def\next{#2}\ifx\next\empty\gdef\STitle{#1}\else\gdef\STitle{#2}\fi%
\def\next{#3}\ifx\next\empty%
    \else\setbox\TheTitle=\vbox{\baselineskip=20pt\let\par=\cr\obeylines%
    \halign{\centerline{\caps##} #3\cr}}\copy\TheTitle\TSkip\TSkip\fi%
\centerline{\caps #4}\TSkip\TSkip%
\def\next{#5}\ifx\next\empty\gdef\SAuthor{#4}\else\gdef\SAuthor{#5}\fi%
\ifx\received\empty\relax
    \else\centerline{\eightrm Received: \received}\fi%
\ifx\revised\empty\TSkip%
    \else\centerline{\eightrm Revised: \revised}\TSkip\fi%
\ifx\communicated\empty\relax
    \else\centerline{\eightrm Communicated by \communicated}\fi\TSkip\TSkip%
\catcode'015=5}}\def\Title{\obeylines\GetTitle}
\def\Abstract{\begingroup\narrower
    \parskip=\medskipamount\parindent=0pt{\caps Abstract. }}
\def\EndAbstract{\par\endgroup\TSkip}

\long\def\MSC#1\EndMSC{\def\arg{#1}\ifx\arg\empty\relax\else
     {\par\narrower\noindent%
     2010 Mathematics Subject Classification: #1\par}\fi}

\long\def\KEY#1\EndKEY{\def\arg{#1}\ifx\arg\empty\relax\else
	{\par\narrower\noindent Keywords and Phrases: #1\par}\fi\TSkip}

\newbox\TheAdd\def\Addresses{\vfill\copy\TheAdd\vfill
    \ifodd\number\lastpage\vfill\eject\phantom{.}\vfill\eject\fi}
{\obeylines\gdef\GetAddress #1
\Address #2 
\Address #3
\Address #4
\EndAddress
{\def\xs{4.3truecm}\parindent=0pt
\setbox0=\vtop{{\obeylines\hsize=\xs#1\par}}\def\next{#2}
\ifx\next\empty 
     \setbox\TheAdd=\hbox to\hsize{\hfill\copy0\hfill}
\else\setbox1=\vtop{{\obeylines\hsize=\xs#2\par}}\def\next{#3}
\ifx\next\empty 
     \setbox\TheAdd=\hbox to\hsize{\hfill\copy0\hfill\copy1\hfill}
\else\setbox2=\vtop{{\obeylines\hsize=\xs#3\par}}\def\next{#4}
\ifx\next\empty\ 
     \setbox\TheAdd=\vtop{\hbox to\hsize{\hfill\copy0\hfill\copy1\hfill}
                \vskip20pt\hbox to\hsize{\hfill\copy2\hfill}}
\else\setbox3=\vtop{{\obeylines\hsize=\xs#4\par}}
     \setbox\TheAdd=\vtop{\hbox to\hsize{\hfill\copy0\hfill\copy1\hfill}
	        \vskip20pt\hbox to\hsize{\hfill\copy2\hfill\copy3\hfill}}
\fi\fi\fi\catcode'015=5}}\gdef\Address{\obeylines\GetAddress}

\makeatletter
\newenvironment{subtheorem}[1]{%
  \def\subtheoremcounter{#1}%
  \refstepcounter{#1}%
  \protected@edef\theparentnumber{\csname the#1\endcsname}%
  \setcounter{parentnumber}{\value{#1}}%
  \setcounter{#1}{0}%
  \expandafter\def\csname the#1\endcsname{\theparentnumber.\Alph{#1}}%
  \ignorespaces
}{%
  \setcounter{\subtheoremcounter}{\value{parentnumber}}%
  \ignorespacesafterend
}
\makeatother
\newcounter{parentnumber}

\theoremstyle{plain}
  \newtheorem{thm}{Theorem}
  \newtheorem{defn}{Definition}
  \newtheorem{prop}{Proposition}
  \newtheorem{cor}{Corollary}
  \newtheorem{lem}{Lemma}

\theoremstyle{definition}
  \newtheorem{example}{Example}
  \newtheorem{rem}{Remark}

\newcommand{\mbb}{\mathbb}
\newcommand{\mbf}{\mathbf}
\newcommand{\mf}{\mathfrak}
\newcommand{\mc}{\mathcal}
\newcommand{\on}{\operatorname}
\newcommand{\ann}{\on{ann}}
\newcommand{\g}{\mathfrak{g}}
\newcommand{\h}{\mathfrak{h}}

\newcommand{\Hom}{\operatorname{Hom}}

\newcommand{\ad}{\mathbf{ad}}
\newcommand{\Ad}{\mathbf{Ad}}
\newcommand{\gr}{\on{gr}}
\newcommand{\hol}{\on{hol}}
\newcommand{\half}{\frac{1}{2}}
\renewcommand{\d}{{\sf{d}}}

\newcommand{\mmat}[2][3em]{\matrix (#2) [matrix of math nodes, row sep=#1,
  column sep=#1, text height=1.5ex, text depth=0.25ex]}
\tikzset{node distance=2cm, auto}

\crefname{pluralequation}{Eqs.}{Eqs.}
\Crefname{pluralequation}{Eqs.}{Eqs.}

\begin{document}
\Title
Moduli spaces for quilted surfaces and Poisson structures
\ShortTitle 
Moduli spaces for quilted surfaces
\SubTitle   
\Author 
David Li-Bland\footnote{\rm D.L-B. was supported by the National Science Foundation under Award No. DMS-1204779.} and Pavol \v{S}evera\footnote{\rm P.\v S. was partially supported by the Swiss National Science Foundation (grants 140985 and 141329).}
\ShortAuthor 
D. Li-Bland and P. \v{S}evera
\EndTitle
\Abstract 
Let $G$ be a Lie group endowed with a bi-invariant pseudo-Riemannian metric. Then the moduli space of flat connections on a principal $G$-bundle, $P\to \Sigma$, over a compact oriented surface with boundary, $\Sigma$, carries a Poisson structure. If we trivialize $P$ over a finite number of points on $\partial\Sigma$ then the moduli space carries a quasi-Poisson structure instead. Our first result is to describe this quasi-Poisson structure in terms of an intersection form on the fundamental groupoid of the surface, generalizing results of Massuyeau and Turaev \cite{Massuyeau:2012uw,Turaev:2007jh}. 

Our second result is to extend this framework to \emph{quilted surfaces}, i.e. surfaces where the structure group varies from region to region and a reduction (or relation) of structure occurs along the borders of the regions, extending results of the second author \cite{Severa:2011ug,Severa98,Severa:2005vla}.

We describe the Poisson structure on the moduli space for a quilted surface in terms of an operation on \emph{spin networks}, i.e. graphs immersed in the surface which are endowed with some additional data on their edges and vertices. This extends the results of various authors \cite{Goldman:1986eh,Goldman:1984hr,Roche:2000ws,Andersen:1996ur}.

\EndAbstract
\MSC 
Primary 53D30, 53D17
\EndMSC
\KEY 
\EndKEY


\tableofcontents

\section{Introduction}
Let $G$ denote a quadratic Lie group, i.e. a Lie group endowed with a bi-invariant pseudo-Riemannian metric. 
If $\Sigma$ is a closed oriented surface, the corresponding moduli space of flat $G$-bundles over $\Sigma$ (note: we always take our $G$-bundles to be principal)
$$\Hom(\pi_1(\Sigma),G)/G$$
carries a symplectic form \cite{Atiyah:1983dt}; more generally, if $\Sigma$ has a boundary, 
then the moduli space carries a Poisson structure.

If $\Sigma$ is connected, and one marks a point on one of the boundary components 
\begin{center}
\begingroup%
  \makeatletter%
  \providecommand\color[2][]{%
    \errmessage{(Inkscape) Color is used for the text in Inkscape, but the package 'color.sty' is not loaded}%
    \renewcommand\color[2][]{}%
  }%
  \providecommand\transparent[1]{%
    \errmessage{(Inkscape) Transparency is used (non-zero) for the text in Inkscape, but the package 'transparent.sty' is not loaded}%
    \renewcommand\transparent[1]{}%
  }%
  \providecommand\rotatebox[2]{#2}%
  \ifx\svgwidth\undefined%
    \setlength{\unitlength}{114.7432251bp}%
    \ifx\svgscale\undefined%
      \relax%
    \else%
      \setlength{\unitlength}{\unitlength * \real{\svgscale}}%
    \fi%
  \else%
    \setlength{\unitlength}{\svgwidth}%
  \fi%
  \global\let\svgwidth\undefined%
  \global\let\svgscale\undefined%
  \makeatother%
  \begin{picture}(1,0.54869191)%
    \put(0,0){\includegraphics[width=\unitlength]{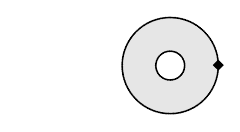}}%
    \put(0.45325387,0.25741511){\color[rgb]{0,0,0}\makebox(0,0)[rb]{\smash{$\Sigma=$}}}%
  \end{picture}%
\endgroup%

\end{center}
and trivializes the principal bundle over that point, the moduli space
$$\Hom(\pi_1(\Sigma),G)$$
becomes \emph{quasi-Poisson} \cite{Alekseev97,Alekseev00,Alekseev99}. In a recent paper \cite{Massuyeau:2012uw}, Massuyeau and Turaev described this quasi-Poisson structure in terms of an intersection form on the loop algebra $\mbb{Z}\pi_1(\Sigma)$, extending a result of Goldman \cite{Goldman:1984hr,Goldman:1986eh}.
The first result of our paper is to generalize their result to the case where $\Sigma$ has multiple marked points (possibly on the same boundary component):
\begin{center}
\begingroup%
  \makeatletter%
  \providecommand\color[2][]{%
    \errmessage{(Inkscape) Color is used for the text in Inkscape, but the package 'color.sty' is not loaded}%
    \renewcommand\color[2][]{}%
  }%
  \providecommand\transparent[1]{%
    \errmessage{(Inkscape) Transparency is used (non-zero) for the text in Inkscape, but the package 'transparent.sty' is not loaded}%
    \renewcommand\transparent[1]{}%
  }%
  \providecommand\rotatebox[2]{#2}%
  \ifx\svgwidth\undefined%
    \setlength{\unitlength}{114.7432251bp}%
    \ifx\svgscale\undefined%
      \relax%
    \else%
      \setlength{\unitlength}{\unitlength * \real{\svgscale}}%
    \fi%
  \else%
    \setlength{\unitlength}{\svgwidth}%
  \fi%
  \global\let\svgwidth\undefined%
  \global\let\svgscale\undefined%
  \makeatother%
  \begin{picture}(1,0.54869191)%
    \put(0,0){\includegraphics[width=\unitlength]{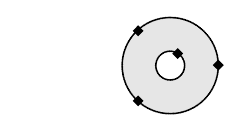}}%
    \put(0.45325387,0.25741511){\color[rgb]{0,0,0}\makebox(0,0)[rb]{\smash{$\Sigma=$}}}%
  \end{picture}%
\endgroup%

\end{center}
These surfaces allow for more economical description of the moduli spaces --- in particular, we show how to obtain them from a collection of discs with two marked points each via iterated fusion.

Blowing up at each of the marked points, we obtain a surface which we call a \emph{domain}:
\begin{center}
\begingroup%
  \makeatletter%
  \providecommand\color[2][]{%
    \errmessage{(Inkscape) Color is used for the text in Inkscape, but the package 'color.sty' is not loaded}%
    \renewcommand\color[2][]{}%
  }%
  \providecommand\transparent[1]{%
    \errmessage{(Inkscape) Transparency is used (non-zero) for the text in Inkscape, but the package 'transparent.sty' is not loaded}%
    \renewcommand\transparent[1]{}%
  }%
  \providecommand\rotatebox[2]{#2}%
  \ifx\svgwidth\undefined%
    \setlength{\unitlength}{114.7432251bp}%
    \ifx\svgscale\undefined%
      \relax%
    \else%
      \setlength{\unitlength}{\unitlength * \real{\svgscale}}%
    \fi%
  \else%
    \setlength{\unitlength}{\svgwidth}%
  \fi%
  \global\let\svgwidth\undefined%
  \global\let\svgscale\undefined%
  \makeatother%
  \begin{picture}(1,0.54869191)%
    \put(0,0){\includegraphics[width=\unitlength,page=1]{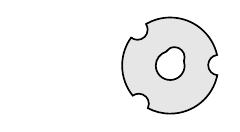}}%
    \put(0.45325387,0.25741511){\color[rgb]{0,0,0}\makebox(0,0)[rb]{\smash{$\hat\Sigma=$}}}%
    \put(0,0){\includegraphics[width=\unitlength,page=2]{Domain.pdf}}%
  \end{picture}%
\endgroup%

\end{center}
We refer to the preimage of any marked point as a \emph{domain wall} (these are the thickened segments of the boundary in the image above).
Our second result is the following: Suppose one chooses a reduced structure separately for each domain wall $\mbf{w}$, i.e. a subgroup $L_\mbf{w}\subset G$.
 If the Lie algebras $\mf{l}_{\mbf{w}}\subset\g$ corresponding to $L_\mbf{w}\subset G$ each satisfy $\mf{l}^\perp_\mbf{w}\subseteq\mf{l}_\mbf{w}$, then the moduli space of 
 \begin{itemize}
 \item flat $G$-bundles over $\hat\Sigma$ equipped with a flat\footnote{Let $P\rvert_{\mbf{w}}\to \mbf{w}$ be a flat principal $G$-bundle. By a  a \emph{flat} reduction of the structure from $G$ to $L_{\mbf{w}}$, we mean a choice of 
principal $L_{\mbf{w}}$-subbundle $Q_{\mbf{w}}\to \mbf{w}$ of $P\rvert_{\mbf{w}}$ which is flat with respect to the connection on $P\rvert_{\mbf{w}}$.} reduction of the structure group from $G$ to $L_{\mbf{w}}$ over each domain wall $\mbf{w}$
 \end{itemize}
is Poisson. 
We may think of this as `coloring' each domain wall with a reduced structure group $L_\mbf{w}\subseteq G$, as pictured below:
\begin{center}
\begingroup%
  \makeatletter%
  \providecommand\color[2][]{%
    \errmessage{(Inkscape) Color is used for the text in Inkscape, but the package 'color.sty' is not loaded}%
    \renewcommand\color[2][]{}%
  }%
  \providecommand\transparent[1]{%
    \errmessage{(Inkscape) Transparency is used (non-zero) for the text in Inkscape, but the package 'transparent.sty' is not loaded}%
    \renewcommand\transparent[1]{}%
  }%
  \providecommand\rotatebox[2]{#2}%
  \ifx\svgwidth\undefined%
    \setlength{\unitlength}{114.7432251bp}%
    \ifx\svgscale\undefined%
      \relax%
    \else%
      \setlength{\unitlength}{\unitlength * \real{\svgscale}}%
    \fi%
  \else%
    \setlength{\unitlength}{\svgwidth}%
  \fi%
  \global\let\svgwidth\undefined%
  \global\let\svgscale\undefined%
  \makeatother%
  \begin{picture}(1,0.54869191)%
    \put(0,0){\includegraphics[width=\unitlength]{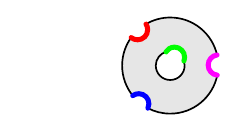}}%
    \put(0.45325387,0.25741511){\color[rgb]{0,0,0}\makebox(0,0)[rb]{\smash{$\Sigma=$}}}%
  \end{picture}%
\endgroup%

\end{center}
In this way we obtain, in particular, the Poisson structures inverting the symplectic forms carried by the moduli spaces of colored surfaces, introduced in \cite{Severa:2011ug} (see also \cite{Severa98,Severa:2005vla}).

Suppose now that $G'$ is a second Lie group whose Lie algebra $\g'$ carries an invariant metric. 
Once again, we choose a reduced structure 
over each domain wall $\mbf{w}'$ on $\Sigma'$, i.e. a subgroup $L_{\mbf{w}'}\subset G'$. If we simultaneously consider flat $G$-bundles over $\Sigma$ and $G'$-bundles over $\Sigma'$ which are 
compatible with the reduced structure on each domain wall, then (as before) the moduli space is Poisson. We picture this as follows:
\begin{center}
\begingroup%
  \makeatletter%
  \providecommand\color[2][]{%
    \errmessage{(Inkscape) Color is used for the text in Inkscape, but the package 'color.sty' is not loaded}%
    \renewcommand\color[2][]{}%
  }%
  \providecommand\transparent[1]{%
    \errmessage{(Inkscape) Transparency is used (non-zero) for the text in Inkscape, but the package 'transparent.sty' is not loaded}%
    \renewcommand\transparent[1]{}%
  }%
  \providecommand\rotatebox[2]{#2}%
  \ifx\svgwidth\undefined%
    \setlength{\unitlength}{116.16311035bp}%
    \ifx\svgscale\undefined%
      \relax%
    \else%
      \setlength{\unitlength}{\unitlength * \real{\svgscale}}%
    \fi%
  \else%
    \setlength{\unitlength}{\svgwidth}%
  \fi%
  \global\let\svgwidth\undefined%
  \global\let\svgscale\undefined%
  \makeatother%
  \begin{picture}(1,0.62256615)%
    \put(0,0){\includegraphics[width=\unitlength]{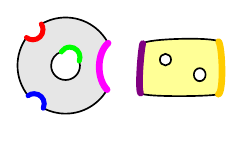}}%
    \put(0.75086479,0.08319391){\color[rgb]{0,0,0}\makebox(0,0)[b]{\smash{$\Sigma'$}}}%
    \put(0.2670103,0.08319391){\color[rgb]{0,0,0}\makebox(0,0)[b]{\smash{$\Sigma$}}}%
  \end{picture}%
\endgroup%

\end{center}
However, one might instead wish to choose a common reduction of structure for two domain walls, $\mbf{w}$ and $\mbf{w}'$ (on $\Sigma$ and $\Sigma'$, resp.). More precisely, to \emph{sew} the domain walls $\mbf{w}$ and $\mbf{w}'$ together is to choose an (orientation reversing) identification $\phi:\mbf{w}\to\mbf{w}'$; $\Sigma\cup_{\phi}\Sigma'$ is then the \emph{sewn surface}. We understand the common image of $\mbf{w}$ and $\mbf{w}'$ to define a domain wall in the sewn surface.
	A quilted surface $\Sigma_{quilt}$ is formed by sewing a collection of domains $\{\Sigma_i\}_{i=1}^n$ together along domain walls, and by 
choosing 
\begin{itemize}
\item a quadratic Lie group $G_i$ for each domain $\Sigma_i$,
\item a reduced structure $L_{\mbf{w}}\subset G_i$ for each unsewn domain wall $\mbf{w}\subset \Sigma_i$; such that $\mf{l}_{\mbf{w}}\subseteq\mf{l}_{\mbf{w}}^\perp$, and
\item a reduced structure $L_{\mbf{w}}\subset G_i\times G_j$ for each sewn domain wall $\mbf{w}\subset \Sigma_i\cap \Sigma_j$; such that $\mf{l}_{\mbf{w}}\subseteq\mf{l}_{\mbf{w}}^\perp$ as a Lie subalgebra of $ \g_i\oplus\bar\g_j$ (where $\bar\g_j$ denotes the Lie algebra $\g_j$ with the metric negated).
\end{itemize}
\begin{center}
\begingroup%
  \makeatletter%
  \providecommand\color[2][]{%
    \errmessage{(Inkscape) Color is used for the text in Inkscape, but the package 'color.sty' is not loaded}%
    \renewcommand\color[2][]{}%
  }%
  \providecommand\transparent[1]{%
    \errmessage{(Inkscape) Transparency is used (non-zero) for the text in Inkscape, but the package 'transparent.sty' is not loaded}%
    \renewcommand\transparent[1]{}%
  }%
  \providecommand\rotatebox[2]{#2}%
  \ifx\svgwidth\undefined%
    \setlength{\unitlength}{130.64716797bp}%
    \ifx\svgscale\undefined%
      \relax%
    \else%
      \setlength{\unitlength}{\unitlength * \real{\svgscale}}%
    \fi%
  \else%
    \setlength{\unitlength}{\svgwidth}%
  \fi%
  \global\let\svgwidth\undefined%
  \global\let\svgscale\undefined%
  \makeatother%
  \begin{picture}(1,0.48183211)%
    \put(0,0){\includegraphics[width=\unitlength]{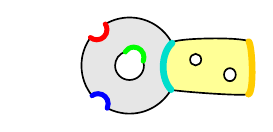}}%
    \put(0.27501314,0.22120111){\color[rgb]{0,0,0}\makebox(0,0)[rb]{\smash{$\Sigma_{quilt}=$}}}%
  \end{picture}%
\endgroup%

\end{center}
Such surfaces (or closely related ones) have played a role in recent developments in both Chern-Simons theory \cite{Kapustin:2010we,Kapustin:2011cn,Kapustin:2010up} and Floer theory \cite{Wehrheim:2010fa,Wehrheim:2010cg}.
	Our second main result is to show that the moduli space, $\mc{M}_{\Sigma_{quilt}}$, of 
\begin{itemize}
\item tuples $\{P_i\to \Sigma_i\}_{i=1}^n$ of flat $G_i$-bundles, equipped with 
\begin{enumerate}
\item a flat\footnote{i.e. 
 a principal $L_\mbf{w}$-subbundle $Q_\mbf{w}\to \mbf{w}$ of $P_i\rvert_\mbf{w}$ which is flat with respect to the connection on $P_i\rvert_\mbf{w}$.} reduction of the structure group of $P_i\rvert_\mbf{w}$ from $G_i$ to $L_\mbf{w}$ over each unsewn domain wall $\mbf{w}\subset \Sigma_i$, and
\item a flat reduction of the structure group of $P_i\times_{\mbf{w}}P_j$ from $G_i\times G_j$ to $L_\mbf{w}$ over each sewn domain wall $\mbf{w}\subset \Sigma_i\cap \Sigma_j$,
\end{enumerate}
\end{itemize}
is Poisson. We will often call $\mc{M}_{\Sigma_{quilt}}$ a \emph{moduli space of flat bundles over a quilted surface}.
	
We provide a description of this Poisson structure in terms of \emph{spin networks} \cite{Penrose:1971vl,Baez:1996hh}, as in \cite{Roche:2000ws,Andersen:1996ur,Goldman:1986eh,Goldman:1984hr}. More precisely, we identify functions $f\in C^\infty(\mc{M}_{\Sigma_{quilt}})$ on the moduli space of flat connections over a quilted surface with spin networks in the quilted surface. Such a spin network $[\Gamma,\ast]$ consists of an immersed graph $\Gamma\to \Sigma_{quilt}$,
\begin{center}
\begingroup%
  \makeatletter%
  \providecommand\color[2][]{%
    \errmessage{(Inkscape) Color is used for the text in Inkscape, but the package 'color.sty' is not loaded}%
    \renewcommand\color[2][]{}%
  }%
  \providecommand\transparent[1]{%
    \errmessage{(Inkscape) Transparency is used (non-zero) for the text in Inkscape, but the package 'transparent.sty' is not loaded}%
    \renewcommand\transparent[1]{}%
  }%
  \providecommand\rotatebox[2]{#2}%
  \ifx\svgwidth\undefined%
    \setlength{\unitlength}{99.94812012bp}%
    \ifx\svgscale\undefined%
      \relax%
    \else%
      \setlength{\unitlength}{\unitlength * \real{\svgscale}}%
    \fi%
  \else%
    \setlength{\unitlength}{\svgwidth}%
  \fi%
  \global\let\svgwidth\undefined%
  \global\let\svgscale\undefined%
  \makeatother%
  \begin{picture}(1,0.62982675)%
    \put(0,0){\includegraphics[width=\unitlength]{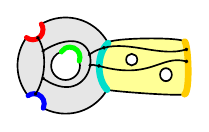}}%
  \end{picture}%
\endgroup%

\end{center}
together with some decoration\footnote{Note, when our structure groups are compact, the graph is decorated with a representation on each edge, and each vertex is decorated with an intertwinor of the representations on the surrounding edges, as in \cite{Roche:2000ws}.} of the edges and vertices of the graph, which (in the introduction) we will denote abstractly by $\ast$.
The Poisson bracket of two spin networks $[\Gamma,\ast]$ and $[\Gamma',\ast']$ is computed as a sum over their intersection points $p\in \Gamma\times_{\Sigma_{quilt}}\Gamma'$,
$$\big\{[\Gamma,\ast],[\Gamma',\ast']\big\}=\sum_{p\in \Gamma\times_{\Sigma_{quilt}}\Gamma'} \pm [\Gamma\cup_p\Gamma',\ast''],$$
where $\Gamma\cup_p\Gamma'$ denotes the union of the two graphs with a common vertex added at the intersection point $p$. This formula generalizes the one found in \cite{Roche:2000ws}.

The basic technical tool we use is a new type of reduction of quasi-Poisson $G$-manifolds by subgroups of $G$.

 In this paper we study the moduli spaces from the (quasi-)Poisson point of view. An alternative approach via (quasi-)symplectic 2-forms (or more generally, Dirac geometry) appears in \cite{LiBland:2013ue}; and the equivalence between these two approaches is explained in \cite{Severa:2014tx}. It should be mentioned that the approach taken in \cite{LiBland:2013ue} allows the construction of more general moduli spaces, while the approach taken in the current paper is more easily quantized (cf. \cite{LiBland:2014da}).
 \begin{rem}
Of course, one could also consider sewing $k$ domains together along a single domain wall, together with a compatible reduction of structure along that domain wall. The moduli space of flat bundles over the resulting branched surfaces is also Poisson \cite{LiBland:2013ue}, and there is an analogous description of the Poisson structure in terms of spin networks. To keep our presentation simple, however, we restrict to the case of quilted surfaces ($k=1,2$).
 \end{rem}

\subsection{Acknowledgements}
The authors would like to thank Eckhard Meinrenken, Alan Weinstein, Marco Gualtieri, Anton Alekseev, Alejandro Cabrera, Dror Bar-Natan, and Jiang-Hua Lu for helpful discussions, explanations, and advice. We'd also like to thank the helpful suggestions we received from various referees which helped to improve the overall readability of the paper.

\section{Quasi-Poisson manifolds}
In this section we recall the basic definitions from the theory of quasi-Poisson manifolds, as introduced by Alekseev, Kosmann-Schwarzbach, and Meinrenken \cite{Alekseev00,Alekseev99}.
 
Let $G$ be a Lie group with Lie algebra $\g$ and  with a chosen Ad-invariant symmetric quadratic tensor, $s\in S^2\g$. Let $\phi\in\bigwedge^3\g$ be the Ad-invariant element defined by
\begin{equation}\label{eq:phi}
\phi(\alpha,\beta,\gamma)=\frac{1}{4}\alpha\bigl([s^\sharp\beta,s^\sharp\gamma]\bigr)\quad (\alpha,\beta,\gamma\in\g^*),
\end{equation}
where $s^\sharp:\g^*\to\g$ is given by $\beta(s^\sharp\alpha)=s(\alpha,\beta)$.

Suppose $\rho:G\times M\to M$ is an action of $G$ on a manifold $M$. Abusing notation slightly, we denote the corresponding Lie algebra action $\rho:\g\to\Gamma(TM)$, by the same symbol. We extend $\rho$ to a Gerstenhaber algebra morphism $\rho:\bigwedge\g\to\Gamma(\bigwedge TM)$.

\begin{defn}
A \emph{quasi-Poisson $G$-manifold} is a triple $(M,\rho,\pi)$, where $M$ is a manifold, $\rho$ an action of $G$ on $M$, and $\pi\in\Gamma(\bigwedge^2 TM)$ a $G$-invariant bivector field, satisfying 
\begin{equation}\label{eq:qPoisCond}\frac{1}{2}[\pi,\pi]=\rho(\phi).\end{equation}
\end{defn}

This definition depends on the choice of $s$. If $G_1$, $G_2$ are Lie groups with chosen elements $s_i\in S^2\g_i$ ($i=1,2$), we set $s=s_1+s_2\in S^2(\g_1\oplus\g_2)$, so that we can speak about quasi-Poisson $G_1\times G_2$-manifolds.  In particular, if $(M_i,\rho_i,\pi_i)$ is a quasi-Poisson $G_i$-manifold ($i=1,2$) then 
$$(M_1,\rho_1,\pi_1)\times(M_2,\rho_2,\pi_2)=(M_1\times M_2, \rho_1\times\rho_2, \pi_1+\pi_2)$$
 is a quasi-Poisson $G_1\times G_2$-manifold. 

\begin{example}\label{ex:G-with-0}
$G$ is a quasi-Poisson $G\times G$-manifold, with the action $\rho(g_1,g_2)\cdot g = g_1 g g_2^{-1}$ and with $\pi=0$.
\end{example}

\begin{rem}\label{rem:CanNegates} Since $s$ appears twice in \cref{eq:phi}, it follows that any quasi-Poisson $(G,s)$-manifold is also a quasi-Poisson $(G,-s)$-manifold. Likewise, any quasi-Poisson $(G_1\times G_2,s_1\oplus s_2)$-manifold is also a quasi-Poisson $(G_1\times G_2,s_1\oplus-s_2)$-manifold.
\end{rem}

Let $\psi\in\bigwedge^2(\g\oplus\g)$ be given by
$$\psi=\frac{1}{2}\sum_{i,j}s^{ij}\,(e_i,0)\wedge(0,e_j)$$
where $s=\sum_{i,j}s^{ij}\,e_i\otimes e_j$ in some basis $e_i$ of $\g$.
\begin{defn}\label{def:Fusion}
If $(M,\rho,\pi)$ is a quasi-Poisson $G\times G\times H$-manifold then its \emph{fusion} is the quasi-Poisson $G\times H$-manifold $(M,\rho^*,\pi^*)$, where $\rho^*(g,h)=\rho(g,g,h)$ and 
$$\pi^*=\pi-\rho(\psi).$$
\end{defn}
Fusion is associative (but not commutative): if $M$ is a quasi-Poisson $G\times G\times G\times H$-manifold then the two $G\times H$-quasi-Poisson structures obtained by the iterated fusions coincide. If $M$ is a quasi-Poisson $G^n\times H$-manifold then its (iterated) fusion to a quasi-Poisson $G\times H$-manifold is given by
\begin{equation}\label{eq:multifusion}
\pi^*=\pi-\sum_{i<j}\rho(\psi_{i,j}),
\end{equation}
where $\psi_{i,j}\in\bigwedge^2(\g^{\oplus n})$ is the image of $\psi$ under the inclusion $\g\oplus\g\to \g^{\oplus n}$ sending the two $\g$'s to $i$'th and $j$'th place respectively.

\section{Reduction and moment maps}\label{sect:red}
A Lie subgroup $C\subseteq G$ will be called \emph{reducing} if its Lie algebra $\mf{c}\subseteq\g$ satisfies
$$\phi(\alpha,\beta,\gamma)=0\quad\forall\alpha,\beta,\gamma\in\ann(\mf{c})$$
where $\ann(\mf{c})\subseteq\g^*$ is the annihilator of $\mf{c}$.
Equivalently,
$$[s^\sharp\alpha,s^\sharp\beta]\in\mf{c}\quad\forall\alpha,\beta\in\ann(\mf{c}).$$
In particular, if $C\subseteq G$ is coisotropic, i.e.\ if $s^\sharp(\ann(\mf{c}))\subseteq\mf{c}$, then $C$ is reducing.

\begin{subtheorem}{thm}\label{thm:Reduction}
\begin{thm}\label{thm:Red}
Suppose that $(M,\rho,\pi)$ is a quasi-Poisson $G$-manifold and that $C\subseteq G$ is a reducing subgroup. Then 
\begin{equation}\label{eq:redBrk}
\{f,g\}:=\pi(\d f,\d g),\quad f,g\in C^\infty(M)^C
\end{equation}
is a Poisson bracket on the space of $C$-invariant functions. In particular, if
the action of $C$ on $M$ induces a regular equivalence relation on $M$\footnote{i.e.\ the orbit space, $M/C$, is a manifold, and the projection $M\to M/C$ is a surjective submersion. Equivalently, by Godement's criterion, the subset $M\times_{M/C} M\subseteq M\times M$ is a closed, embedded submanifold, and the projection onto either factor, $M\times_{M/C} M\to M$, is a submersion.}, then the bivector field $\pi$ descends to define a Poisson structure on $M/C$.
\end{thm}
\begin{proof}
The proof is essentially the same as that of \cite[Theorem 4.2.2]{Alekseev99}, but we include it for completeness.

First, we observe that $\{f,g\}\in C^\infty(M)^C$, since $f,g$ and $\pi$ are each $C$-invariant.

To see that the bracket \eqref{eq:redBrk} satisfies the Jacobi identity, notice that
$$\{f_1,\{f_2,f_3\}\}+c.p.=\half
[\pi,\pi](\d f_1,\d f_2,\d f_3)=\rho(\phi)(\d f_1,\d f_2,\d f_3),$$
for any $f_i\in C^\infty(M)^C$ ($i=1,2,3$).
Now $\rho^*\d f_i\in \ann(\mf{c})$ and $C$ is reducing, hence $\rho(\phi)(\d f_1,\d f_2,\d f_3)=0$.
\end{proof}

For $\xi\in\g$ let $\xi^L$ and $\xi^R$ denote the corresponding left and right invariant vector field on $G$.

\begin{defn}\label{def:twistedMomMap}
Let $(M,\rho,\pi)$ be a quasi-Poisson $G$-manifold and let $\tau:G\to G$ be an $s$-preserving automorphism. A map $\mu:M\to G$ is a \emph{($\tau$-twisted) moment map} if it is equivariant for the action 
\begin{equation}\label{eq:TwistAct}g\cdot \tilde g=\tau(g)\,\tilde g\,g^{-1}\end{equation}
 of $G$ on $G$, and if the image of $\pi$ under
$$\mu_*\otimes id:TM\otimes TM\to TG\otimes TM$$
is 
$$-\half \sum_{i,j}s^{ij}\bigl(e_i^L+\tau(e_i)^R\bigr)\otimes\rho(e_j).$$
\end{defn}

We shall use moment maps to get Poisson submanifolds of $M/C$, in analogy with Marsden-Weinstein reduction (under certain non-degeneracy conditions these submanifolds will be the symplectic leaves of $M/C$). First we need an analogue of coadjoint orbits.

\begin{lem}
If $C\subseteq G$ is a reducing subgroup then
$$\hat{\mf c}:=\{(\xi+s^\sharp\alpha,\xi-s^\sharp\alpha);\, \xi\in\mf{c}, \alpha\in\ann(\mf{c})\}$$
is a Lie subalgebra of $\g\oplus\g$. 
\end{lem}
\begin{proof}
Let $\xi,\eta\in\mf c$, $\alpha,\beta\in\ann(\mf c)$. Since
$$[\xi\pm s^\sharp\alpha,\eta\pm s^\sharp\beta]=\bigl([\xi,\eta]+[s^\sharp\alpha,s^\sharp\beta]\bigr)
\pm s^\sharp(\ad^*_\xi\beta-\ad^*_\eta\alpha)$$
and $[s^\sharp\alpha,s^\sharp\beta]\in\mf c$, the space $\hat{\mf c}$ is closed under the Lie bracket.
\end{proof}

Let $\hat C\subseteq G\times G$ be a Lie group with the Lie algebra $\hat{\mf c}$; and we suppose that the diagonal inclusion $\mf{c}\subseteq \hat{\mf{c}}$ lifts to an inclusion $C\subseteq \hat C\subseteq G\times G$. The group $G\times G$ acts on $G$ by 
\begin{equation}\label{eq:tautwist}(g_1,g_2)\cdot g=\tau(g_1)\,g\,g_2^{-1}.\end{equation} The orbits of $\hat C\subseteq G\times G$ on $G$ will serve as analogues of coadjoint orbits.

\begin{thm}\label{thm:mmred}
Let $(M,\rho,\pi)$ be a quasi-Poisson $G$-manifold with a moment map $\mu:M\to G$ and $C\subseteq G$ a reducing subgroup. Suppose that the $C$ action on $M$ induces a regular equivalence relation. Let $\mathcal{O}\subseteq G$ be a $\hat C$-orbit. If $\mathcal O$ is in a clean position relative to $\mu$,\footnote{That is, $\on{gr}(\mu)\cap(\mc{O}\times M)\subseteq G\times M$ is a submanifold, and $T(\on{gr}(\mu)\cap(\mc{O}\times M))=T\on{gr}(\mu)\cap (T\mc{O}\times TM)$, where $\on{gr}(\mu)=\{(\mu(x),x);\  x\in M\}$. It happens, in particular, if $\mu$ is transverse to $\mathcal O$.} then
$$\mu^{-1}(\mathcal{O})/C\subseteq M/C$$
 is a Poisson submanifold. More generally, if $S\subseteq G$ is a $\hat C$-stable submanifold and $S$ is in a clean position relative to $\mu$, then $\mu^{-1}(S)/C\subseteq M/C$ is a Poisson submanifold.
\end{thm}
\begin{proof} This follows from \cref{thm:GHred} when $H$ is trivial.
%
\end{proof}


\subsection*{Partial reduction}
One can generalize both \cref{thm:Red} and \cref{thm:mmred} in order to reduce quasi-Poisson $G\times H$-manifolds to quasi-Poisson $H$-manifolds:

\begin{thm}\label{thm:GHred}
Suppose that $(M,\rho,\pi)$ is a quasi-Poisson $G\times H$-manifold and  $C\subseteq G$ is a reducing subgroup.
\begin{enumerate}
\item
If the $C$ action on $M$ induces a regular equivalence relation, then
the bivector field $\pi$ descends to define a quasi-Poisson $H$-structure on $M/C$.
\item
Let $\tau_G$ and $\tau_H$ be automorphisms of $G$ and $H$, and let  $$(\mu_G,\mu_H):M\to G\times H$$ be a $(\tau_G,\tau_H)$-twisted moment map. If $S\subseteq G$ is a $\hat C$-stable submanifold, and the following topological conditions hold:
\begin{itemize}
\item $S$ is in a clean position relative to  $\mu_G$, and
\item the action of $C$ on $\mu_G^{-1}(S)$ induces a regular equivalence relation,
\end{itemize}
then $\mu_G^{-1}(S)/C$ is a quasi-Poisson $H$-manifold, and $\mu_H$ descends to define a $\tau_H$-twisted moment map,
$$\mu_H:\mu_G^{-1}(S)/C\to H.$$ 
\end{enumerate}
Moreover, if all the stated assumptions hold, then $\mu_G^{-1}(S)/C\subseteq M/C$ is a quasi-Poisson $H$-submanifold.
\end{thm}
\end{subtheorem}
\begin{proof}
\begin{enumerate}
\item 
The proof of this first statement is very similar to that of \cref{thm:Red}.

Let $q:M\to M/C$ denote the quotient map (a surjective submersion).
$\pi$ is $G\times H$-invariant, thus it descends to define a bivector field $\pi'$ on $M/C$. Likewise, the $H$-action on $M$ descends to define an action $\rho':H\times M/C\to M/C$.

Now the $\Ad$-invariant element $s\in S^2(\g\oplus\h)$ splits as the sum of $s_G\in S^2(\g)$ and $s_H\in S^2(\h)$; similarly $\phi\in \bigwedge^3(\g\oplus\h)$ splits as the sum of $\phi_G\in \bigwedge^3\g$ and $\phi_H\in \bigwedge^3\h$.
By definition, we have
$$\half
[\pi,\pi]=\rho(\phi)=\rho(\phi_G)+\rho(\phi_H).$$
However, for any $\alpha_i\in \Omega^1(M/C)$, ($i=1,2,3$) the pullbacks $\rho^*q^*\alpha_i\in \ann(\mf{c})$, and $C$ is reducing, hence $\rho(\phi_G)(q^*\alpha_1,q^*\alpha_2,q^*\alpha_3)=0$. Thus
$$\half
[\pi',\pi']=\rho'(\phi_H),$$
as desired.

\item First of all, since $S$ is in a clean position relative to $\mu_G$, $\mu_G^{-1}(S)\subseteq M$ is a submanifold. Moreover, since the action of $C$ induces a regular equivalence relation, $\mu_G^{-1}(S)/C$ is a manifold and the quotient map is a surjective submersion.
Since $H$ acts trivially on $S\subseteq G$, so $H$ restricts to an action on $\mu_G^{-1}(S)/C\subseteq M/C$.

For any $x\in \mu_G^{-1}(S)$ and $\alpha\in T^*_xM$ the moment map condition gives
\begin{equation}\label{eq:MMaux}
(\mu_G)_*\bigl(\pi(\cdot,\alpha)\bigr)=-\half\bigl((s_G^\sharp\rho_x^*\alpha)^L+\tau(s_G^\sharp\rho_x^*\alpha)^R\bigr).
\end{equation}
If $\alpha$ is annihilates $\rho_x(\mf{c})$, then $\rho_x^*\alpha\in\ann(\mf c)$. The vector \eqref{eq:MMaux} is thus the action of $\half(s_G^\sharp\rho_x^*\alpha,-s_G^\sharp\rho_x^*\alpha)\in\hat{\mf c}$ on $G$. In particular, since $S$ is in a clean position relative  to $\mu_G$, it follows that $\pi(\cdot,\alpha)$ is tangent to $\mu_G^{-1}(S)$. This implies that $\pi$ descends to  $\mu_G^{-1}(S)/C$ to define
a quasi-Poisson $H$-structure. 


 Finally, since $\mu_H:\mu_G^{-1}(S)\to H$ is $G\times H$-equivariant, it descends to a map on $\mu_G^{-1}(S)/C$.
 The image of $\pi$ under $(\mu_H)_*\otimes \on{id}: TM\otimes TM\to TH\otimes TM$ is 
 $$-\half \sum_{i,j} s_H^{ij}\bigl(e_i^L+\tau(e_i)^R\bigr)\otimes\rho(e_j),$$
 where $s_H\in S^2(\mf{h})^H$ denotes the chosen invariant symmetric tensor,
 and hence this also holds for the reduced bivector field on 
 $\mu_G^{-1}(S)/C$, proving that $\mu_H$ descends to define a moment map.
 
 If the action of $C$ on $M$ induces a regular equivalence relation, then $\mu_G^{-1}(S)/C\subseteq M/C$ is certainly an $H$-invariant submanifold, and the previous argument shows that the bivector field is tangent to $\mu_G^{-1}(S)/C$. That is, $\mu_G^{-1}(S)/C\subseteq M/C$ is a quasi-Poisson $H$-submanifold.
 \end{enumerate}
 
\end{proof}

\begin{example}\label{ex:FuseRedIsRed}
Suppose that $(M,\rho,\pi)$ is a quasi-Poisson $G\times G\times H$ manifold. Let 
$$G_\Delta:=\{(\xi,\xi)\in G\times G\}\subseteq G\times G$$ 
denote the diagonal subgroup.

One may first fuse the two $G$ factors of the quasi-Poisson structure on $M$ (in either order), yielding a quasi-Poisson $G\times H$ structure on $M$, and then reduce by $G$, to obtain a quasi-Poisson $H$ structure on $M/G_\Delta$. Alternatively, notice that $G_\Delta\subseteq G\times G$ is a reducing subgroup, so one may reduce $M$ directly by $G_\Delta$, as in Theorem~\ref{thm:GHred}, obtaining a quasi-Poisson $H$ structure on $M/G_\Delta$. Both these quasi-Poisson $H$ structures on $M/G_\Delta$ are identical.

\end{example}

\subsection{Induction and the Hat construction}\label{sec:HatConstr}

When $C\subseteq \hat C$ is closed, and $\mc{O}\subset G$ is a $\hat C$-orbit, the reduced space $\mu_G^{-1}(\mc{O})/C$ can be conveniently described in the following way.  Let $\hat M$ be the manifold obtained from $M$ by induction from $C$ to $\hat C$ (using the diagonal embedding $C\subset \hat C$). Concretely,
\begin{equation}\label{eq:Mhat}
\hat M=(\hat C\times M)/C
\end{equation}
where the $C$-action on $\hat C \times M$ is $c\cdot(\hat c, m)=(\hat c\, c^{-1},c\cdot m)$.

Now suppose $\mc{O}=\hat C\cdot g_0$ for some $g_0\in G$, then 
 \begin{equation}\label{eq:stab}\mu_G^{-1}(\mc{O})/C\cong \hat\mu_G^{-1}(g_0)/\on{Stab}(g_0),\end{equation}
  where $$\hat\mu_G:\hat M\xrightarrow{(\hat c,m)\to \hat c\cdot \mu_G(m)} G,$$ is the induced $\hat C$-equivariant moment map and 
  $$\on{Stab}(g_0)=\{\hat c\in\hat C;\,\hat c\cdot g_0=g_0\}\subseteq \hat C.$$

More generally, suppose that $S\subseteq G$ is a $\hat C$-invariant submanifold. Then for any submanifold $X\subseteq G$ such that $S=\hat C\cdot X$, we have natural bijections
\begin{equation}\label{eq:GenIndQuo}\mu_G^{-1}(S)/C\cong \hat\mu_G^{-1}(S)/\hat C\cong \hat\mu_G^{-1}(X)/\on{Stab}(X),\end{equation}
where the stabilizer groupoid, $\on{Stab}(X)\rightrightarrows X$, is defined as the restriction of the action groupoid $\hat C\ltimes G$ to $X\subseteq G$.
\begin{rem}[The stabilizer groupoid]\label{rem:Stab}
Recall that the stabilizer groupoid is defined as
$\on{Stab}(X):=\{(\hat c,x)\in \hat C\times X\mid \hat c\cdot x\in X\}$
 with source and target maps given by $\mbf{s}(\hat c,x)= x$, and $\mbf{t}(\hat c,x)= \hat c\cdot x$, respectively, and 
multiplication given by $(\hat c',\hat c\cdot x)\circ (\hat c, x):=(\hat c'\cdot \hat c,x).$ 

The action of $\on{Stab}(X)$ on $\hat\mu_G^{-1}(X)$ is the obvious one: the moment map is $\hat\mu_G$, and for any $\hat m\in\hat\mu_G^{-1}(X)$ and $\big(\hat c,\hat\mu_G(\hat m)\big)\in\on{Stab}(X)$, we have
$\big(\hat c,\hat\mu_G(\hat m)\big)\cdot \hat m= \hat c\cdot \hat m.$
In particular, two elements $\hat m,\hat m'\in \hat\mu_G^{-1}(X)$ are in the same $\on{Stab}(X)$-orbit if and only if there exists a $\hat c\in\hat C$ such that  $\hat m'= \hat c\cdot \hat m$.
\end{rem}

The following proposition gives sufficient conditions on $X\subseteq G$ for the quotient \eqref{eq:GenIndQuo} to be a smooth quasi-Poisson manifold.

\begin{prop}\label{prop:hatclean}
Suppose that $C\subseteq G$ is a reducing subgroup, and $C\subseteq \hat C$ is closed.
Suppose further that $(M,\rho,\pi)$ is a quasi-Poisson $G\times H$-manifold with a $(\tau_G,\tau_H)$-twisted moment map, $$(\mu_G,\mu_H):M\to G\times H.$$ If $X\subseteq G$, $\on{Stab}(X)\subseteq \hat C\ltimes G$, and $S:=\hat C\cdot X\subseteq G$ are embedded submanifolds, 
%
and the following topological conditions hold:
\begin{itemize}
\item $X$ is in a clean position relative to $\hat\mu_G$, and
\item the action of $\on{Stab}(X)$ on $\hat\mu_G^{-1}(X)$ induces a regular equivalence relation,
\end{itemize}
then $\hat\mu_G^{-1}(X)/\on{Stab}(X)$ is a quasi-Poisson $H$-manifold, and $\mu_H$ descends to define a $\tau_H$-twisted moment map,
$$\mu_H:\hat\mu_G^{-1}(X)/\on{Stab}(X)\to H.$$ 
\end{prop}
\begin{proof}
%
 Since $\mu_G^{-1}(S)/C\cong \hat\mu_G^{-1}(X)/\on{Stab}(X)$, 
 it suffices to show that the assumptions of Theorem~\ref{thm:GHred}.2 hold.  

We begin by showing that $S$ is in a clean position relative to $\mu_G$. First, consider the diagonal action map $$A:\hat C\ltimes (G\times \hat M)\xrightarrow{(\hat c;g,\hat m)\to (\hat c\cdot g,\hat c\cdot\hat m)} (G\times \hat M).$$

By assumption, $U:=\hat C\times (X\times \hat M)\subseteq \hat C\ltimes (G\times \hat M)$ and $V:=\hat C\times \gr(\hat\mu_G) \subseteq \hat C\ltimes (G\times \hat M)$ intersect cleanly, with $U\cap V\cong \hat C\times \hat\mu_G^{-1}(X)$. Since $V$ is a union of fibres of $A$, it follows that $A(U)=S\times \hat M$ and $A(V)=\gr(\hat\mu_G)$ intersect cleanly (cf. Lemma~\ref{lem:CleanInter}). That is, $S$ is in a clean position relative to $\hat \mu_G$. Moreover, the restriction $A:U\cap V\to A(U)\cap A(V)$ defines a surjective submersion
$$\hat C\times \hat\mu_G^{-1}(X)\to \hat\mu_G^{-1}(S).$$

Next, consider the $\hat C$-equivariant map $B:G\times \hat M\xrightarrow{(g,[\hat c,m])\to [\hat c]} \hat C/C$. The restrictions of $B$ to both $S\times \hat M$ and $\gr(\hat\mu_G)$ are surjective submersions (since both submanifolds are $\hat C$-invariant). Thus $S\times M=\left(B\rvert_{S\times \hat M}\right)^{-1}(C)$ and $\gr(\mu_G)=\left(B\rvert_{\gr(\hat\mu_G)}\right)^{-1}(C)$ intersect cleanly. That is, $S$ is in a clean position relative to $\mu_G$.

Since the action Lie groupoids $C\ltimes \mu_G^{-1}(S)$ and $\on{Stab}(X)\ltimes_X\hat\mu_G^{-1}(X)$ are Morita equivalent, and the equivalence relation defined by the action of $\on{Stab}(X)$ on $\hat\mu_G^{-1}(X)$ is regular, it follows that the equivalence relation corresponding to the action of $C$ on $\mu_G^{-1}(S)$ is also regular (cf. Lemma~\ref{lem:MoritaEqReg}).

\end{proof}

\section{Quasi-Poisson structures on moduli spaces}\label{sec:QPsMs}
Let $\Sigma$ be a compact oriented surface with boundary, and let $V\subset\partial\Sigma$ be a finite collection of ``\emph{marked points}'' such that every component of $\Sigma$ intersects $V$. Let $\Pi_1(\Sigma,V)$ denote the fundamental groupoid of $\Sigma$ with the base set $V$. The composition in $\Pi_1(\Sigma,V)$ is from right to left: $ab$ means path $b$ followed by path $a$. For $a\in\Pi_1(\Sigma,V)$ let $\mbf{s}(a)$ denote the source and $\mbf{t}(a)$ the target of $a$; $ab$ is defined if $\mbf{t}(b)=\mbf{s}(a)$.

 Let
$$M_{\Sigma,V}(G)=\Hom(\Pi_1(\Sigma,V),G).$$
$M_{\Sigma,V}(G)$ can be seen as the moduli space of pairs: 
\begin{equation}\label{eq:FramedBndls}(P\to \Sigma, \{\hat v\in P\rvert_{v}\}_{v\in V}),\end{equation} consisting of a flat (principal) $G$-bundle over $\Sigma$ together with a framing over $V$.\footnote{By a framing over $V$, we mean a choice of a lift $\hat v\in P_v$ of each marked point $v\in V$ to the fibre over $v$.}

For any arrow $a\in\Pi_1(\Sigma,V)$  let
$$\hol_a:M_{\Sigma,V}(G)\to G$$
denote evaluation at $a$ (in terms of flat connections it is the holonomy along $a$). There is a natural action $\rho=\rho_{\Sigma,V}$ of the group $G^V$ on $M_{\Sigma,V}(G)$ which is defined by 
\begin{equation}\label{eq:VertexAction}\hol_a(\rho(g)x)=g_{\mbf{t}(a)}\hol_a(x)g_{\mbf{s}(a)}^{-1}.\end{equation} Infinitesimally,
$$(\hol_a)_*(\rho(\xi))=-\xi^R_{\mbf{t}(a)}+\xi^L_{\mbf{s}(a)}$$
 for any $\xi\in\g^V$, where $\xi^{L/R}$ denotes the left/right invariant vector field on $G$ corresponding to $\xi\in\g$.

By a \emph{skeleton} of $(\Sigma,V)$ we mean a graph $\Gamma\subset\Sigma$ with the vertex set $V$, such that there is a deformation retraction of $\Sigma$ to $\Gamma$.\footnote{It is a simple exercise to show that there exists a skeleton for every marked surface. However, it should be emphasized that this skeleton is not unique! The edges of a skeleton determine generators of $\Pi_1(\Sigma,V)$, but the converse fails (generators for $\Pi_1(\Sigma,V)$ may be chosen which intersect in an essential way).} If we choose an orientation of every edge of $\Gamma$ then
 $M_{\Sigma,V}(G)$ gets identified (via $(\hol_a,a\in E_\Gamma)$) with $G^{E_\Gamma}$, where ${E_\Gamma}$ is the set of edges of $\Gamma$. In particular, if $\Sigma$ is a disc and $V$ has two elements then we get $M_{\Sigma,V}(G)= G$.

Our convention in this paper is to orient the $\partial\Sigma$ \emph{against} the induced boundary orientation.
Cutting the boundary of $\Sigma$ at the marked points, $V$, splits it  into oriented \emph{arcs} (the components of $\partial\Sigma$ that don't contain a marked point\footnote{that is, an element of $V$} are not considered to be arcs).  If we choose an ordered pair $(v_{\mbf{s}},v_{\mbf{t}})$ of marked points ($v_{\mbf{s}}\neq v_{\mbf{t}}\in V$) then the corresponding \emph{fused surface} $\Sigma^*$ is obtained by gluing a short piece of the arc ending at $v_{\mbf{s}}$ with a short piece of the arc starting at $v_{\mbf{t}}$ (so that $v_{\mbf{s}}$ and $v_{\mbf{t}}$ get identified).
The subset $V^*\subset\partial\Sigma^*$ is obtained from $V$ by identifying $v_{\mbf{s}}$ and $v_{\mbf{t}}$. Notice that the map 
$$M_{\Sigma^*,V^*}(G)\to M_{\Sigma,V}(G),$$
 coming from the map $(\Sigma,V)\to(\Sigma^*,V^*)$, is a diffeomorphism: if $\Sigma$ retracts to a skeletal graph $\Gamma$ then $\Sigma^*$ retracts to its image $\Gamma^*$, and the two graphs have the same number of edges. We can thus identify the manifolds $M_{\Sigma^*,V^*}(G)$ and $M_{\Sigma,V}(G)$.
\begin{figure}[h!]
\centering
\def\svgwidth{10cm}
\input{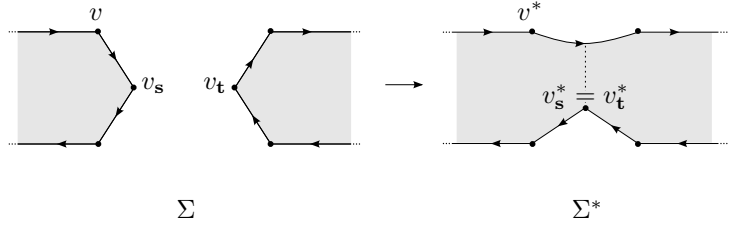}
\caption{Fusion. The marked points $v_{\mbf{s}},v_{\mbf{t}}\in V$ are identified after fusion, while every other marked point $v\in V$ ($v_{\mbf{s}}\neq v\neq v_{\mbf{t}}$) is unaffected. Notice that $v_{\mbf{s}}$ remains the source of the same edge after fusion, and $v_{\mbf{t}}$ remains the target of the same edge after fusion.}
\end{figure}

Every $(\Sigma,V)$ can be obtained by fusion from a collection of discs, each with two marked points: If $\Gamma\subset\Sigma$ is a skeleton then the subset of $\Sigma$ that retracts onto an edge $e\in E_\Gamma$ is a disc $D_e$, and $\Sigma$ is obtained from $D_e$'s by repeated fusion.

\begin{figure}[h!]
\centering
\includegraphics[width=6cm]{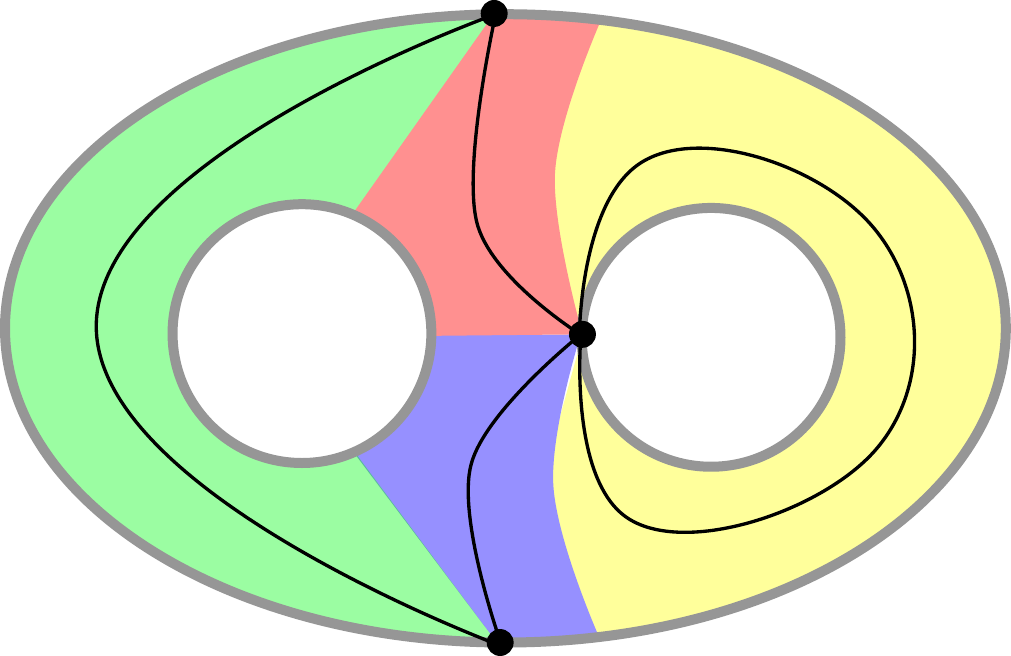}
\caption{A surface with a skeleton, fused from four discs. The discs have been assigned colors in the picture, but this figure should not be confused with a quilted surface.}
\end{figure}

\begin{thm}\label{thm:qPoiss-mod}
There is a natural bivector field $\pi_{\Sigma,V}$ on $M_{\Sigma,V}(G)$ such that $$(M_{\Sigma,V}(G),\rho_{\Sigma,V},\pi_{\Sigma,V})$$
 is a quasi-Poisson manifold, uniquely determined by the properties
\begin{enumerate}
\item if $\Sigma$ is a disc and $V$ has two elements then $\pi_{\Sigma,V}=0$
\item if $(\Sigma,V)=(\Sigma_1,V_1)\sqcup(\Sigma_2,V_2)$ then $\pi_{\Sigma,V}=\pi_{\Sigma_1,V_1}+\pi_{\Sigma_2,V_2}$
\item \label{thm:qPoiss-mod_aux}if $(\Sigma^*,V^*)$ is obtained from $(\Sigma,V)$ by fusion, then\\ $(M_{\Sigma^*,V^*}(G),\rho_{\Sigma^*,V^*},\pi_{\Sigma^*,V^*})$ is obtained from $(M_{\Sigma,V}(G),\rho_{\Sigma,V},\pi_{\Sigma,V})$ by the corresponding fusion.
\end{enumerate}
\end{thm}

If there is no danger of confusion, we shall denote $\pi_{\Sigma,V}$ simply by $\pi$.

\begin{rem}
Alejandro Cabrera has independently studied quasi-Hamiltonian $G^V$-structures for the marked surfaces described above.
\end{rem}

Once we choose a skeleton $\Gamma$ of $(\Sigma,V)$, \cref{thm:qPoiss-mod} gives us a formula for the quasi-Poisson structure on $M_{\Sigma,V}$, as $(\Sigma,V)$ is a fusion of a collection of discs with two marked points. Let us denote the resulting bivector field on $M_{\Sigma,V}(G)$ by $\pi_\Gamma$. Theorem  \ref{thm:qPoiss-mod} follows from the following Lemma:

\begin{lem}\label{lem:indep}
The bivector field $\pi_\Gamma$ on $M_{\Sigma,V}(G)$ is independent of the choice of $\Gamma$. 
\end{lem}

\begin{rem}
The lemma follows from the special case where $(\Sigma,V)$ is a disc with 3 marked points (see Example \ref{ex:triangle} below). However, we shall give a different proof in the next section.
\end{rem}

\begin{proof}[Proof of Theorem \ref{thm:qPoiss-mod}]
By the lemma we have a well-defined quasi-Poisson structure on $M_{\Sigma,V}(G)$. Properties (1)--(3) of the theorem are satisfied by the construction of $\pi_\Gamma$.
\end{proof}

Let us now describe the calculation of $\pi=\pi_\Gamma$ in more detail. Notice that for any vertex $v$ of $\Gamma$, the (half)edges adjacent to $v$ are linearly ordered: a cyclic order is given by the orientation of $\Sigma$. Since $v$ is on the boundary, the cyclic order is actually a linear order. $\Gamma$ is a \emph{ciliated graph} in the terminology of Fock and Rosly \cite{Fock:1999wz}. 

We choose an orientation of every edge of $\Gamma$ to get an identification $M_{\Sigma,V}(G)=G^{E_\Gamma}$. First we see it as a $G^{E_\Gamma}\times G^{E_\Gamma}$-quasi-Poisson space with zero bivector (i.e.\ as $M_{\Sigma',V'}(G)$, where $(\Sigma',V')$ is a disjoint union of discs with two marked points each). Then, fusing at each vertex using the \emph{reversed}\footnote{We reverse the linear order when fusing to account for the minus sign appearing in \eqref{eq:multifusion} (cf. \cite{Alekseev00}).} linear order, we obtain a $G^V$-quasi-Poisson space. 

\begin{example}
As the simplest example, suppose $(\Sigma,V)$ is an annulus with a single marked point (on one of the boundary circles). Then $(\Sigma,V)$ may be obtained by fusion from a disc $(\Sigma',V')$ with two marked points, as in \cref{fig:AnnulusFusion}. Now $M_{\Sigma',V'}=G$ with the the quasi-Poisson $G\times G$-structure described in \cref{ex:G-with-0}: the bivector field is trivial and $G\times G$ acts by $(g_1,g_2)\cdot g=g_1gg_2^{-1}$. Thus $M_{\Sigma,V}=G$, the $G$-action is by conjugation, and
$$\pi=\frac{1}{2}\sum_{i,j}s^{ij}\,e_i^R\wedge e_j^L.$$
\begin{figure}
\begin{center}
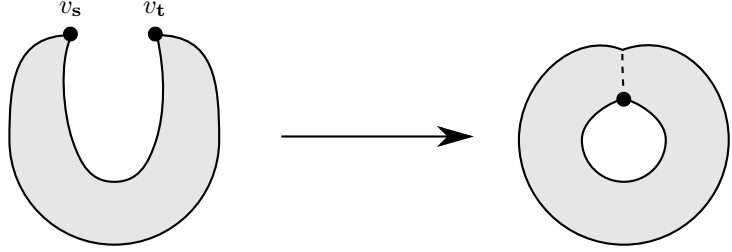
\end{center}
\caption{The annulus with one marked point is obtained by fusion from the disc with two marked points.\label{fig:AnnulusFusion}}
\end{figure}

\end{example}

\begin{example}\label{ex:triangle}
Let $\Sigma$ be a triangle and $V$ is the set of its vertices.
$$
\begin{tikzpicture}[thick,decoration={markings,mark=at position 0.5 with {\arrow{>}}}]
\fill[black!20!white] (0,0)--(2,0)--(1,1.73)--cycle; 
\draw[postaction={decorate}](0,0)--node[auto,swap]{$a$}(2,0);
\draw[postaction={decorate}](2,0)--node[auto,swap]{$b$}(1,1.73);
\draw[postaction={decorate}](1,1.73)--node[auto,swap]{$c$}(0,0);
\end{tikzpicture}
$$
We can identify $M_{\Sigma,V}$ with $G^2$ via $(\hol_{a^{-1}},\hol_{b})$, i.e.\ $\Gamma$ is the graph with the oriented edges $a^{-1}$, $b$. In this case
\begin{equation}\label{eq:piTriang}\pi=-\frac{1}{2}\sum_{i,j}s^{ij}\,e_i^L(1)\wedge e_j^L(2)\end{equation}
(where $e_i^L(k)$ denotes the left-invariant vector field which is tangent to the $k^{th}$ factor of $G^2$ ($k=1,2$)). Equivalently,
$$\pi(\hol_{a^{-1}}^*\theta^L, \hol_{b}^*\theta^L)=-s\in\g\otimes\g,$$
where $\theta^L\in\Omega^1(G,\g)$ is the left-invariant Maurer-Cartan form. An easy calculation shows
$$\pi(\hol_{b^{-1}}^*\theta^L, \hol_{c}^*\theta^L)=-s\in\g\otimes\g,$$
confirming that $\pi$ is independent of the choice of $\Gamma$.
\end{example}

For a general surface $(\Sigma,V)$ with a choice of a skeleton $\Gamma$, we get  an identification $M_{\Sigma,V}=G^{E_\Gamma}$. Applying \eqref{eq:multifusion} we obtain
\begin{equation}\label{eq:FR}
\pi_{\Sigma,V}=\frac{1}{2}\sum_{v\in V}\sum_{a<b}\sum_{i,j}s^{ij}\,e_i(a,v)\wedge e_j(b,v)
\end{equation}
where $a,b$ run over the (half)edges adjacent to $v$,
$$e_i(a,v)=
 \begin{cases}
 -e_i^R(a) & a\text{ goes into }v\\
 e_i^L(a) & a\text{ goes out of }v
 \end{cases}
$$
and for $a\in E_\Gamma$,  $e_i^{R,L}(a)$ denotes the right/left-invariant vector field on $G^{E_\Gamma}$ equal to $$(0,\dots,0,\overbrace{e_i}^{a},0,\dots,0)\in\g^{E_\Gamma}$$ at the identity element. 

\begin{rem}\label{rem:FR}
Essentially the same formula  was discovered by Fock and Rosly \cite{Fock:1999wz}, for Poisson structures on $M_{\Sigma,V}$ obtained by a choice of a classical $r$-matrix. More precisely, if one considers the bivector described in Equation (16) of \cite{Fock:1999wz}, but replacing the $r$-matrix $r^{ij}$ with it's symmetrization $s^{ij}$, then one arrives at the same formula as \eqref{eq:FR}.

 Meanwhile, Skovborg studied the corresponding formula in the absence of an $r$-matrix for invariant functions \cite{Skovborg:2009ve}.
\end{rem}

\section{The homotopy intersection form and quasi-Poisson structures}

Massuyeau and Turaev \cite{Massuyeau:2012uw} made a beautiful observation that, in the case of one marked point and $G=GL_n$, the quasi-Poisson structure on $M_{\Sigma,V}(G)$ can be expressed in terms of the homotopy intersection form on $\pi_1(\Sigma)$, introduced by Turaev in \cite{Turaev:2007jh}. Here we extend their result to the case of arbitrary $(G,s)$ and arbitrary $V$.

\begin{figure}[h]
\begin{center}
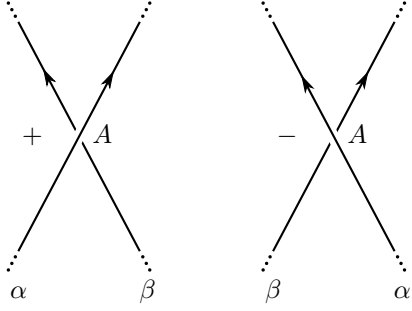
\end{center}
\caption{\label{fig:Signp} $\on{sign}(A)=\pm 1$ is determined by comparing the orientation of $\alpha$ and $\beta$ with that of $\Sigma$.}
\end{figure}

Let us first extend (a skew-symmetrized version of) Turaev's homotopy intersection form to fundamental groupoids.
If $a,b\in\Pi_1(\Sigma,V)$, let us represent them by transverse smooth paths $\alpha,\beta$. For any point $A$ in their intersection, let
\begin{align*}
\lambda(A)&=
\begin{cases}
1 &\text{if }A\in\partial\Sigma\\
2 &\text{otherwise}
\end{cases}\\
\on{sign}(A)&:=
\begin{cases}
1 &\text{if }(\dot\alpha|_A,\dot\beta|_A)\text{ is positively oriented}\\
-1 &\text{otherwise.}
\end{cases}
\end{align*}
as in \cref{fig:Signp}.

Let $\alpha_A$ denote the portion of $\alpha$ parametrized from the source of $\alpha$ up to the point $A$. Finally, let
\begin{equation}\label{eq:IntPairing}(a,b):=\sum_A \lambda(A)\on{sign}(A) [\alpha_A^{-1}\beta_A]\in\mathbb{Z}\Pi_1(\Sigma,V),\end{equation}
where $\mbb{Z}\Pi_1(\Sigma,V)$ denotes the groupoid ring\footnote{As a $\mbb{Z}$-module, $\mathbb{Z}\Pi_1(\Sigma,V)$ is freely generated by $\Pi_1(\Sigma,V)$. For generators $a,b\in \Pi_1(\Sigma,V)$, their product in $\mathbb{Z}\Pi_1(\Sigma,V)$ is $$ab:=\begin{cases}
a\circ b &\text{if } a \text{ and } b \text{ are composable, i.e. }\mbf{s}(a)=\mbf{t}(b)\\0&\text{otherwise.}
\end{cases}$$} generated by $\Pi_1(\Sigma,V)$ over $\mbb{Z}$.
As in \cite{Turaev:2007jh} one can check that $(a,b)$ is well defined, i.e.\ independent of the choice of $\alpha$ and $\beta$.

Let us list the properties of $(a,b)$.
 For $x\in\mathbb{Z}\Pi_1(\Sigma,V)$, $x=\sum n_i a_i$, let $\bar x=\sum n_i a_i^{-1}$.

\begin{prop}\label{prop:hpairing}
The pairing 
$$(\cdot,\cdot):\Pi_1(\Sigma,V)\times\Pi_1(\Sigma,V)\to\mathbb{Z}\Pi_1(\Sigma,V)$$
satisfies
\begin{enumerate}
\item $(a,b)$ is a linear combination of paths from the source of $b$ to the source of $a$
\item $(b,a)=-\overline{(a,b)}$
\item \label{prop:hpairing_auxx}$(a,bc)=(a,c)+(a,b)c$
\item \label{prop:hpairing_aux} if $(\Sigma^*,V^*)$ is obtained from $(\Sigma,V)$ by fusing $\Sigma$ at $v_{\mbf{s}},v_{\mbf{t}}\in V$ and $a^*,b^*$ denotes the image of $a,b$ in $\Pi_1(\Sigma^*,V^*)$, then
\begin{multline*}
(a^*,b^*)=(a,b)^*-(\delta_{\mbf{t}(a),v_{\mbf{s}}}\,a^{-1}-\delta_{\mbf{s}(a),v_{\mbf{s}}}1_{v_{\mbf{s}}^*})(\delta_{\mbf{t}(b), v_{\mbf{t}}}\,b-\delta_{\mbf{s}(b), v_{\mbf{t}}}1_{v_{\mbf{s}}^*})+\\
+(\delta_{\mbf{t}(a),v_{\mbf{t}}}\,a^{-1}-\delta_{\mbf{s}(a),v_{\mbf{t}}}1_{v_{\mbf{s}}^*})(\delta_{\mbf{t}(b), v_{\mbf{s}}}\,b-\delta_{\mbf{s}(b), v_{\mbf{s}}}1_{v_{\mbf{s}}^*}).
\end{multline*}
\end{enumerate}
It is the only pairing satisfying these properties.
\end{prop}

\begin{proof}
The fact that \eqref{eq:IntPairing} satisfies these properties is readily verified.

We now show that these properties uniquely determine a pairing \begin{equation}\label{eq:intPairAxiomatic}(\cdot,\cdot):\Pi_1(\Sigma,V)\times\Pi_1(\Sigma,V)\to\mathbb{Z}\Pi_1(\Sigma,V).\end{equation}

First suppose $(\Sigma,V)$ is a disjoint union of disks each of which is marked with two points, then any pairing which satisfies the first three properties must be trivial: From property 1, we see that $(a,b)=0$ unless $a$ and $b$ lie in the same connected component of $\Sigma$. This leaves the following cases to check:
\begin{subequations}
\begin{align}
\label{eq:aunit}a&=1_v\text{ is a unit, for some }v\in V\\
\label{eq:bunit}b&=1_v\text{ is a unit, for some }v\in V\\
\label{eq:aisb}b&=a\\
\label{eq:aisbinv}b&=a^{-1}
\end{align}
\end{subequations}
Now in case \eqref{eq:bunit} property 3 (with $c=1_v$) implies that $(a,1_v)=0$. 
%
%
Now in case \eqref{eq:aunit}, property 2 implies $(a,b)=(1_v,b)=-\overline{(b,1_v)}=0$. 
In case \eqref{eq:aisb}, $(a,a)=k1_{\mbf{s}(a)}$ (by property 1), for some $k\in\mbb{Z}$. Then property 2 implies that $k1_{\mbf{s}(a)}=-\overline{k1_{\mbf{s}(a)}}:=-k1_{\mbf{s}(a)}$, which forces $k=0$, and hence $(a,b)=(a,a)=0$. 
Finally, in case \eqref{eq:aisbinv}, property 3 implies $(a,1_{\mbf{t}(a)})=(a,a^{-1})+(a,a)a^{-1}$. Since $(a,1_{\mbf{t}(a)})=0=(a,a)$, we have $(a,b)=(a,a^{-1})=0$.

Next, notice that properties 3 and 4 determine the pairing for a fused surface $(\Sigma^*,V^*)$ in terms of the pairing for $(\Sigma, V)$.
Since any marked surface can be obtained by fusing a collection of disks (with two marked points each), a finite number of times, this shows that \eqref{eq:intPairAxiomatic} is uniquely determined by properties 1-4.
\end{proof}

For $a,b\in \Pi_1(\Sigma,V)$ let us consider the $\g\otimes\g$-valued function on $M_{\Sigma,V}$,
$$\pi(\hol_{a}^*\theta^L,\hol_{b}^*\theta^L).$$
These functions, in turn, specify $\pi$ completely.

We can now state our version of the result of Massuyeau and Turaev \cite{Massuyeau:2012uw}.

\begin{thm}\label{thm:pi_hpair}
For any $a,b\in\Pi_1(\Sigma,V)$ we have
\begin{equation}\label{eq:pi_hpair}
\pi(\hol_{a}^*\theta^L,\hol_{b}^*\theta^L)=\frac{1}{2}(\Ad_{\hol_{(a,b)}}\otimes 1)\,s.
\end{equation}
\end{thm}
\begin{rem}
Essentially the same formula was discovered independently by Xin Nie \cite{Nie:2013ux}.
\end{rem}

\begin{proof}[Proof of Theorem \ref{thm:pi_hpair} and of Lemma \ref{lem:indep}]
To prove both Theorem \ref{thm:pi_hpair} and Lemma \ref{lem:indep} (and thus finish the proof of Theorem \ref{thm:qPoiss-mod}) we need to check that
\begin{equation}\label{eq:pi_Gamma_hpair}
\pi_\Gamma(\hol_{a}^*\theta^L,\hol_{b}^*\theta^L)=\frac{1}{2}(\Ad_{\hol_{(a,b)}}\otimes 1)\,s
\end{equation}
for every $a,b\in\Pi_1(\Sigma,V)$ and any skeleton $\Gamma$ of $(\Sigma,V)$.

Notice (via $\hol_{bc}^*\theta^L=\hol_c^*\theta^L+\Ad_{(\hol_c)^{-1}}\hol_b^*\theta^L$) that
 \begin{multline*}\pi_\Gamma(\hol_{a}^*\theta^L,\hol_{bc}^*\theta^L)=\pi_\Gamma(\hol_{a}^*\theta^L,\hol_{c}^*\theta^L)+\\(1\otimes \Ad_{(\hol_c)^{-1}})\pi_\Gamma(\hol_{a}^*\theta^L,\hol_{b}^*\theta^L).\end{multline*}
As a result, if \eqref{eq:pi_Gamma_hpair} is true for all $a,b$ in a set of generators of $\Pi_1(\Sigma,V)$, it is then true (by Proposition \ref{prop:hpairing} Part \ref{prop:hpairing_auxx}) for all elements of $\Pi_1(\Sigma,V)$.

Equation \eqref{eq:pi_Gamma_hpair} is true if $(\Sigma,V)$ is a disc with two marked points, as both sides of the equation vanish. The same is true for the disjoint union of a collection of such discs. As any $(\Sigma,V,\Gamma)$ can be obtained from such a collection by a repeated fusion, it remains to check that \eqref{eq:pi_Gamma_hpair} is preserved under fusion.

Suppose that \eqref{eq:pi_Gamma_hpair} is satisfied for some $(\Sigma,V)$ and its skeleton $\Gamma$. Let $(\Sigma^*,V^*)$ be a fusion of $(\Sigma,V)$ and let $\Gamma^*$ be the image of $\Gamma$ in $\Sigma^*$. Then $\pi_{\Gamma^*}$ is obtained from $\pi_\Gamma$ by the corresponding quasi-Poisson fusion. By Proposition \ref{prop:hpairing} Part \ref{prop:hpairing_aux} we then get
$$\pi_{\Gamma^*}(\hol_{a}^*\theta^L,\hol_{b}^*\theta^L)=\frac{1}{2}(\Ad_{\hol_{(a^*,b^*)}}\otimes 1)\,s.$$
In other words, \eqref{eq:pi_Gamma_hpair} is satisfied also for $(\Sigma^*,V^*,\Gamma^*)$ for the elements of $\Pi_1(\Sigma^*,V^*)$ in the image of $\Pi_1(\Sigma,V)$. As the image generates $\Pi_1(\Sigma^*,V^*)$, we conclude that \eqref{eq:pi_Gamma_hpair} is satisfied for $(\Sigma^*,V^*,\Gamma^*)$.
\end{proof}

\begin{rem}
Theorem \ref{thm:pi_hpair} can be used as an alternative definition of $\pi$. Properties 1--3 of the homotopy intersection form (cf. \cref{prop:hpairing}) mean that there is a unique $G^V$-invariant bivector field $\pi$ satisfying \eqref{eq:pi_hpair}. Property 4 means that $\pi$ is compatible with fusion.
\end{rem}

If $f:(\Sigma',V')\to(\Sigma,V)$ is an embedding then clearly 
$$(f_*a,f_*b)=f_*(a,b)$$
for every $a,b\in\Pi_1(\Sigma',V')$.
From Theorem \ref{thm:pi_hpair} we thus get the following result.
\begin{cor}\label{cor:emb}
If $f:(\Sigma',V')\to(\Sigma,V)$ is an embedding then
$$f^*:M_{\Sigma,V}(G)\to M_{\Sigma',V'}(G)$$
is a quasi-Poisson map, i.e.\ $\pi_{\Sigma,V}$ and $\pi_{\Sigma',V'}$ are $f^*$-related.
\end{cor}
Let us now consider the special case $\Sigma'=\Sigma$, $V'\subset V$. Recall that if $(M,\rho,\pi)$ is a $G\times H$-quasi-Poisson manifold and if $M/G$ is a manifold (e.g.\ if the action of $G$ is free and proper) then $\pi$ descends to a bivector field $\pi'$ on $M/G$ such that $p_*\pi=\pi'$, where $p:M\to M/G$ is the projection, and that $M/G$ thus becomes a $H$-quasi-Poisson manifold, called the \emph{quasi-Poisson reduction of $M$ by $G$} (cf. \cite{Alekseev99,Alekseev00} or Theorem~\ref{thm:GHred}). Using this terminology, Corollary \ref{cor:emb} becomes
\begin{cor}
If $V'\subset V$ then $M_{\Sigma,V'}(G)$ is the quasi-Poisson reduction of $M_{\Sigma,V}(G)$ by $G^{V\backslash V'}$.
\end{cor}

\begin{figure}
\begin{center}
  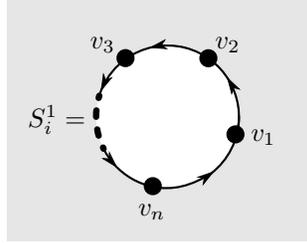
\end{center}
  \caption{\label{fig:perm} Pictured, is a connected component of the boundary, $S^1_i\subseteq \partial\Sigma$. $S^1_i$ is oriented \emph{against} the induced boundary orientation.
  The marked points $\{v_1,\dots,v_n\}= V\cap S^1_i$ inherit a cyclic order, and the permutation is defined as $\sigma(v_i)=v_{i+1}$ (where the index $i+1$ is computed modulo $n$).}
\end{figure}

Finally, again following Massuyeau and Turaev \cite{Massuyeau:2012uw}, we can define a moment map for the quasi-Poisson $G^V$-manifold $M_{\Sigma,V}(G)$. Let us orient $\partial\Sigma$ \emph{against} the orientation induced from $\Sigma$, and let $S^1_1\sqcup\dots\sqcup S^1_n=\partial\Sigma$ be the decomposition into connected components. The subset of marked points $V\cap S^1_i$ which all lie on a given component of the boundary inherit a cyclic order from the orientation of $\partial\Sigma$. Let 
$\sigma:V\to V$ be the corresponding permutation (as in Figure~\ref{fig:perm}).
%
 Let $\tau:G^V\to G^V$ be the automorphism defined for any $g\in G^V$ and $v\in V$ by
\begin{equation*}\tau(g)_v=g_{\sigma(v)},\end{equation*}
so that the $\tau$-twisted action of $G^V\times G^V$ on $G^V$ is
\begin{equation}\label{eq:sigTwist}
\big((g_1,g_2)\cdot g\big)_v=(g_1)_{\sigma(v)}\:g_v \:(g_2)_v^{-1}
\end{equation}
(cf. \cref{eq:tautwist}). 
 For every $v\in V$ let $\mu_v:M_{\Sigma,V}(G)\to G$ be
$$\mu_v=\hol_{a_v}$$
where $a_v$ is the boundary arc from $v$ to $\sigma(v)$. Let us combine the maps $\mu_v$ to a single map $\mu:M_{\Sigma,V}(G)\to G^V$.

\begin{thm}\label{thm:modmoment}
The map $\mu:M_{\Sigma,V}(G)\to G^V$ is a $\tau$-twisted moment map.
\end{thm}
\begin{proof}
The equivariance of $\mu$ is obvious. 
If $v\in V$ and $b\in\Pi_1(\Sigma,V)$ then by the definition of $(a_v,b)$ we have
$$(a_v,b)=-\delta_{\mbf{s}(b),v}\,1_v-\delta_{\mbf{s}(b),\sigma(v)}\,a_v^{-1}
+\delta_{\mbf{t}(b),v}\,b+\delta_{\mbf{t}(b),\sigma(V)}\,a_v^{-1}b.$$
Theorem \ref{thm:pi_hpair} then implies that $\mu$ is a moment map, as the 1-forms $\hol_b^*\theta^L$ span the cotangent bundle of $M_{\Sigma,V}(G)$. 
\end{proof}

\begin{rem}
An alternative way of proving Theorem \ref{thm:modmoment} is to verify it for the case of a disc with two marked points on the boundary, and then to use fusion.
\end{rem}

\section{Surfaces with boundary data}

For every point $v\in V$ let us choose a coisotropic subgroup $C_v\subseteq G$.\footnote{If $G$ is semisimple then any parabolic subgroup is coisotropic; if $G$ is simple then these are all the coisotropic subgroups}
Let $C=\prod_{v\in V} C_v\subseteq G^V$. Suppose the action of $C$ on $M_{\Sigma,V}(G)$ defines a regular equivalence relation, and consider the orbit space
$$M_{\Sigma,V}(G,(C_v)_{v\in V}):=M_{\Sigma,V}(G)/C.$$
Geometrically, $M_{\Sigma,V}(G,(C_v)_{v\in V})$ is the moduli space of 
flat (principal) $G$-bundles $P\to\Sigma$ equipped with a reduction to $C_v$ over $v$ for every $v\in V$.

Since $C\subseteq G^V$ is coisotropic, by Theorem \ref{thm:Red} we know that $M_{\Sigma,V}(G,(C_v)_{v\in V})$ is Poisson. More generally, if $V'\subseteq V$ then  \cref{thm:GHred} implies that
$$M_{\Sigma,V}(G)/\prod_{v\in V'}C_v$$
is a quasi-Poisson $G^{V\smallsetminus V'}$-manifold.

\begin{example}\label{ex:PoisLie1}
Let $\Sigma$ be a triangle and $V$ the set of its vertices. Suppose that $s\in S^2\g$ is non-degenerate, and that $(\g,\mf a,\mf b)$ is a Manin triple. Let $A,B\subset G$ be the corresponding subgroups, and let us suppose that $A\cap B=\{1\}$ and $BA=G$, i.e. $B\times A\xrightarrow{(b,a)\to ba}G$ is a diffeomorphism.

Let us choose the subgroup, $C_v$, at two of the vertices to be $B$ and the remaining vertex to be $A$, as in the picture below
\begin{center}
\begingroup%
  \makeatletter%
  \providecommand\color[2][]{%
    \errmessage{(Inkscape) Color is used for the text in Inkscape, but the package 'color.sty' is not loaded}%
    \renewcommand\color[2][]{}%
  }%
  \providecommand\transparent[1]{%
    \errmessage{(Inkscape) Transparency is used (non-zero) for the text in Inkscape, but the package 'transparent.sty' is not loaded}%
    \renewcommand\transparent[1]{}%
  }%
  \providecommand\rotatebox[2]{#2}%
  \ifx\svgwidth\undefined%
    \setlength{\unitlength}{124.62785645bp}%
    \ifx\svgscale\undefined%
      \relax%
    \else%
      \setlength{\unitlength}{\unitlength * \real{\svgscale}}%
    \fi%
  \else%
    \setlength{\unitlength}{\svgwidth}%
  \fi%
  \global\let\svgwidth\undefined%
  \global\let\svgscale\undefined%
  \makeatother%
  \begin{picture}(1,0.8355146)%
    \put(0,0){\includegraphics[width=\unitlength]{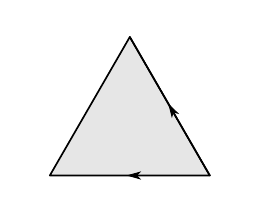}}%
    \put(0.50405007,0.07932995){\color[rgb]{0,0,0}\makebox(0,0)[b]{\smash{$g_1$}}}%
    \put(0.66792859,0.45215367){\color[rgb]{0,0,0}\makebox(0,0)[lb]{\smash{$g_2$}}}%
    \put(0.18448677,0.07362547){\color[rgb]{0,0,0}\makebox(0,0)[rb]{\smash{$B$}}}%
    \put(0.81541925,0.07362547){\color[rgb]{0,0,0}\makebox(0,0)[lb]{\smash{$A$}}}%
    \put(0.50405007,0.72255329){\color[rgb]{0,0,0}\makebox(0,0)[b]{\smash{$B$}}}%
  \end{picture}%
\endgroup%

\end{center}
Now the holonomies $g_1,g_2\in G$ along the edges pictured above identify $M_{\Sigma,V}(G)$ with $G\times G$, where the bivector field was described in \cref{eq:piTriang}. The diffeomorphism $$\big((b,a)\to ba\big):B\times A\to G$$ identifies $M_{\Sigma,V}$ with $B\times A\times B\times A$, and the action of $C=B\times A\times B$ on $M_{\Sigma,V}(G)$ becomes
$$(b,a,b')\cdot(b_1,a_1,b_2,a_2)=(bb_1,a_1a^{-1},b'b_2,a_2a^{-1}).$$
Thus the map $(b_1,a_1,b_2,a_2)\to a_2a_1^{-1}$ identifies $M_{\Sigma,V}(G)/C$ with the Lie group $A$. Using \eqref{eq:piTriang} we see that the resulting bivector field on $A\cong B\backslash G$ is the pushforward of 
$$\pi =\frac{1}{2}\sum_i{\xi_i}^L\wedge{\eta^i}^L,$$ where $\{\eta^i\}\subset\mf{b}$ and $\{\xi_i\}\subset\mf{a}$ are bases in duality.

A comparison shows that the Poisson structure on $M_{\Sigma,V}(G)/C\cong A$ is the Poisson-Lie structure on $A$ described in \cite{thesis-3}.

The assumption that $G=BA$ is rather strong. In Example~\ref{ex:PoissLie2}, we will identify the Poisson Lie group with a moduli space in the absence of this assumption.
\end{example}

\begin{example}
If $\Sigma,V$ is a disc with two marked points and $C\subseteq G$ is a coisotropic subgroup which we embed as a subgroup of the second factor of $G\times G$, then $M_{\Sigma,V}(G)/C=G/C$ is a quasi-Poisson $G$-manifold, with $\pi=0$.
\end{example}

Since, according to Theorem \ref{thm:modmoment}, the holonomies along the boundary arcs define a moment map $\mu:M_{\Sigma,V}(G)\to G^V$, we can apply the moment map reduction (\cref{thm:mmred}) to get Poisson submanifolds of $M_{\Sigma,V}(G,(C_v)_{v\in V})=M_{\Sigma,V}(G)/C$.\footnote{If $s\in S^2\g$ is nondegenerate and if every component of $\partial\Sigma$ contains a marked point then these submanifolds are in fact the symplectic leaves. See \cite{LiBland:2013ue} for more details.} To give geometric descriptions of these Poisson submanifolds it will be convenient to use the ``hat''-construction from Section~\ref{sec:HatConstr}. We begin be describing $\hat M$, $\hat\mu$ and $\hat C$.

Recall that $\sigma:V\to V$ is the permutation obtained by walking along $\partial\Sigma$ against the orientation induced from $\Sigma$, and that $a_v$ is the boundary arc from $v$ to $\sigma(v)$.

We now describe the group $\hat C$. Let $\mf{c}_v^\perp=s^\sharp\ann(\mf c_v)$; it is an ideal in ${\mf c}_v$. Let $C_v^\perp\subseteq C_v$ denote the corresponding connected Lie group.\footnote{If $G$ is simple, so that $C_v$ is parabolic, then $C_v^\perp\subseteq C_v$ is the nilpotent radical.} Then $C_v$ normalizes $C_v^\perp$. Hence
$$\hat C_v=\{(g_1,g_2)\in C_v\times C_v\mid \ g_1g_2^{-1}\in C_v^\perp\}\subset G\times G$$
is a subgroup. 
We take
$$\hat C=\prod\hat C_v\subset G^{V}\times G^V.$$

\begin{figure}[h!]
\centering
\def\svgwidth{9cm}
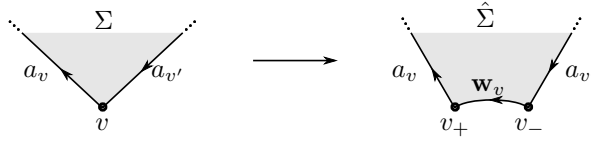
\caption{$\hat\Sigma$ is obtained from $\Sigma$ by blowing up at each point $v\in V$. We denote the exceptional divisor (a segment on the boundary of $\Sigma$) by $\mbf{w}_v$, its initial point by $v_-$ and its endpoint by $v_+$. Note that $v=\sigma(v')$, and the orientation of the arcs is opposite to the induced boundary orientation.}\label{fig:Sigmahat}
\end{figure}

Let us now give a geometric description of the induced $\hat C$-manifold $\hat M$ (see Eq.~\eqref{eq:Mhat}) and of the induced $\hat C$-equivariant moment map $\hat\mu:\hat M\to G^V$ in the case of $M=M_{\Sigma,V}(G)$.
First, let $\hat\Sigma$ be the surface obtained from $\Sigma$ by blowing up at each point $v\in V$, as in \cref{fig:Sigmahat}. We let $\mbf{w}_v$ denote the exceptional divisor obtained by blowing up at $v$. We label the initial and end points of the segment $\mbf{w}_v$ by $v_-$ and $v_+$, respectively. We let $\mbf{Wall}$ denote the set of $\mbf{w}_v$'s, and we let $V_-$ and $ V_+$ denote the set of initial and end points of the $\mbf{w}_v$'s. Thus $(\hat\Sigma, V_-\cup V_+)$ is a marked surface.

Then
\begin{equation}\label{eq:hatM}\hat M=\{f\in M_{\hat\Sigma,V_-\cup V_+}(G);\ f(\mbf{w}_v)\in C^\perp_v\ (\forall v\in V)\}.\end{equation}
The group $G^V\times G^V\cong G^{V_-\cup V_+}$ acts naturally on $M_{\hat\Sigma, V_-\cup V_+}(G)$ (cf. \cref{eq:VertexAction}) and the subgroup $\hat C$ preserves $\hat M\subset M_{\hat\Sigma, V_-\cup V_+}(G)$. The elements of $\hat C$ are  $(g_{ v_-},g_{v_+}\in C_v)_{v\in V}$ such that $g_{ v_+} g_{ v_-}^{-1}\in C_v^\perp$. 
The map $\hat\mu:\hat M\to G^V$ is given by
$$\hat\mu(f)=(f(a_v))_{v\in V}.$$

Suppose we now choose an element $h=(h_v)_{v\in V}\in G^V$. Recall from \cref{eq:sigTwist} that the  action of $\hat C\subset G^{V_-\cup V_+}$ on $G^V$ is by
$$(g\cdot h)_v=g_{\sigma(v)_-}h_v g_{v_+}^{-1}.$$
 Let $\mathcal{O}\subseteq G^V$ denote the $\hat C$-orbit containing $h\in G^V$. Using Eq.\ \eqref{eq:stab} we get
$$\mu^{-1}(\mathcal O)/C\cong\hat\mu^{-1}(\mathcal O)/\hat C\cong \{f\in\hat M;\ f(a_v)=h_v\ (\forall v\in V)\}/\on{Stab}(h).$$

\begin{thm}\label{thm:bdrdata}
Suppose that the graph of $\hat\mu$ composes cleanly with $h\in G^V$, 
and the $\on{Stab}(h)$-orbits of $\hat\mu^{-1}(h)$ form a regular foliation, then the quotient space 
\begin{multline*}\hat\mu^{-1}(h)/\on{Stab}(h)\\=\bigl\{f\in M_{\hat\Sigma,V_-\cup V_+}(G);\ f(\mbf{w}_v)\in C^\perp_v, f(a_v)=h_v\ (\forall v\in V)\bigr\}/\on{Stab}(h),\end{multline*}
is a Poisson manifold. 

Moreover, if the $C$-orbits of $M_{\Sigma,V}(G)$ form a regular foliation, then $\hat\mu^{-1}(h)/\on{Stab}(h)$ is naturally isomorphic to a Poisson submanifold of $M_{\Sigma,V}(G)/C$.
\end{thm}
\begin{proof}
We apply Proposition~\ref{prop:hatclean} to prove the first statement. To do so, we only need to show that $C\subseteq \hat C$ is closed. But for each $v\in V$, we have $\hat C_v\cong C_v^\perp\rtimes C_v$ via the map $(g_1,g_2)\to (g_1g_2^{-1},g_2)$; hence $C_v\subseteq \hat C_v$ is closed, which yields the result.

The proof of Proposition~\ref{prop:hatclean} shows that the assumptions of Theorem~\ref{thm:GHred}.2 hold, and if the $C$-orbits of $M_{\Sigma,V}(G)$ form a regular foliation, then the assumptions of Theorem~\ref{thm:GHred}.1 also hold. Hence $\hat\mu^{-1}(h)/\on{Stab}(h)$ is naturally isomorphic to a Poisson submanifold of $M_{\Sigma,V}(G)/C$.
\end{proof}

Theorem \ref{thm:bdrdata} is particularly interesting in the case when $h=\mbf{1}$ is the unit. Since this constrains the holonomies along the paths $a_v$ to be trivial, we can contract these paths to points.

\begin{figure}[h!]

\centering
\def\svgwidth{12cm}
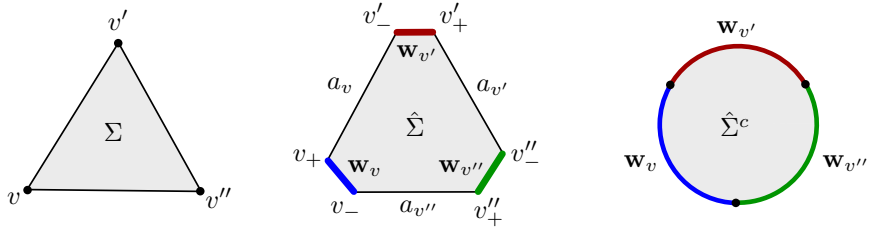
\caption{\label{fig:tilSig}Blow up at the marked points, followed by contraction of $a_v$'s}
\label{fig:contracted}
\end{figure}

\begin{example}\label{ex:SymplDbl}[\cite{Severa98,Severa:2011ug,Severa:2005vla}]
Once again, let  $(\g,\mf a,\mf b)$ be a Manin triple, let $\Sigma,V$ be a disk with four marked points, and let us choose the subgroups $C_v$ as in the picture,
\begin{center}
\begingroup%
  \makeatletter%
  \providecommand\color[2][]{%
    \errmessage{(Inkscape) Color is used for the text in Inkscape, but the package 'color.sty' is not loaded}%
    \renewcommand\color[2][]{}%
  }%
  \providecommand\transparent[1]{%
    \errmessage{(Inkscape) Transparency is used (non-zero) for the text in Inkscape, but the package 'transparent.sty' is not loaded}%
    \renewcommand\transparent[1]{}%
  }%
  \providecommand\rotatebox[2]{#2}%
  \ifx\svgwidth\undefined%
    \setlength{\unitlength}{72.69290161bp}%
    \ifx\svgscale\undefined%
      \relax%
    \else%
      \setlength{\unitlength}{\unitlength * \real{\svgscale}}%
    \fi%
  \else%
    \setlength{\unitlength}{\svgwidth}%
  \fi%
  \global\let\svgwidth\undefined%
  \global\let\svgscale\undefined%
  \makeatother%
  \begin{picture}(1,0.88168309)%
    \put(0,0){\includegraphics[width=\unitlength]{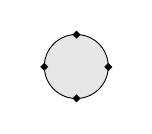}}%
    \put(0.49679312,0.72491664){\color[rgb]{0,0,0}\makebox(0,0)[b]{\smash{$B$}}}%
    \put(0.49679312,0.10035065){\color[rgb]{0,0,0}\makebox(0,0)[b]{\smash{$B$}}}%
    \put(0.77055002,0.4287448){\color[rgb]{0,0,0}\makebox(0,0)[lb]{\smash{$A$}}}%
    \put(0.25149029,0.4287448){\color[rgb]{0,0,0}\makebox(0,0)[rb]{\smash{$A$}}}%
    \put(0.50902737,0.41868066){\color[rgb]{0,0,0}\makebox(0,0)[b]{\smash{$\Sigma$}}}%
  \end{picture}%
\endgroup%

\end{center}
where $A,B\subset G$ are Lie groups integrating $\mf a$ and $\mf b$ with $A\cap B=\{\mbf{1}\}$.
Setting $h=\mbf{1}$ (which is a regular value for $\hat\mu$), we get the constraint on the holonomies pictured in \cref{fig:SymplDbl}.

\begin{figure}[h!]
\begin{center}
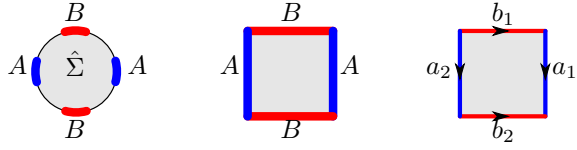
\end{center}
\caption{Setting $h=\mbf{1}$ is effectively the same as contracting the uncolored portions of the boundary in the first picture; this results in the second picture. We label the holonomies as in the third picture. Since $\hat\Sigma$ is contractible, the product of the holonomies around the boundary must be trivial, that is: $a_1b_1=b_2a_2$.}\label{fig:SymplDbl}
\end{figure}

In this case $\on{Stab}(h)=\mbf{1}$ and $\hat\mu^{-1}(\mbf{1})/\on{Stab}(\mbf{1})$ is
$$\hat\mu^{-1}(\mbf{1})=\{(a_1,b_1,a_2,b_2)\in A\times B\times A\times B;\ a_1b_1=b_2a_2\},$$
where the holonomies are as pictured in \cref{fig:SymplDbl}.\footnote{In the case that the discrete group $A\cap B$ is non-trivial, $\hat\mu^{-1}(\mbf{1})$ is a covering space for $M_{\Sigma,V}(G)/C$, and thus inherits a Poisson structure.}

$M_{\Sigma,V}(G)/C$ is the Lu-Weinstein double symplectic groupoid (cf. \cite{Lu:1989vb,thesis-3}). The moduli space for this surface was first identified with the Lu-Weinstein double symplectic groupoid in the work of the second author \cite{Severa98,Severa:2011ug,Severa:2005vla}.

In the case that $G=AB$, we can identify the moduli space with $G$ via
$$g=a_1b_1=b_2a_2.$$
The resulting Poisson structure on $G$ is the so-called Heisenberg double.

\end{example}

\subsection*{Description in terms of flat bundles}
We now give an alternative description of $M_{\Sigma,V}/C$ and its Poisson submanifolds in terms of flat bundles. Given a domain wall $\mbf{w}_v\in \mbf{Wall}$ obtained by blowing up the marked point $v\in V$, we define $C_{\mbf{w}_v}:=C_v$.

\begin{defn}
Suppose that $\mbf{w}\subseteq \partial\hat\Sigma$ is a domain wall. A \emph{flat} reduction of the structure from $G$ to $C_\mbf{w}$ over $\mbf{w}$, is a choice of 
principal $C_\mbf{w}$-subbundle $Q_{\mbf{w}}\subseteq P\rvert_{\mbf{w}}$ which is flat with respect to the connection on $P\rvert_{\mbf{w}}$.\footnote{Equivalently, this is a choice of a flat principal $C_\mbf{w}$-bundle $Q_{\mbf{w}}\to\mbf{w}$ together with an isomorphism of flat bundles $P\rvert_{\mbf{w}}\to Q_{\mbf{w}}\times_{C_\mbf{w}}G$.}
\end{defn}

 \begin{thm}\label{thm:FullModSpc}
 Suppose the $C$-orbits of $M_{\Sigma,V}$ form a regular foliation. Then the moduli space of pairs \begin{equation}\label{eq:PairMod}(P\to\hat\Sigma,\{Q_{\mbf{w}}\to \mbf{w}\}_{\mbf{w}\in \mbf{Wall}})\end{equation}
of flat $G$-bundles over $\hat\Sigma$ equipped with a flat reduction of the structure group from $G$ to $C_\mbf{w}$ over each domain wall $\mbf{w}\in\mbf{Wall}$, carries a natural Poisson structure and can be canonically identified  with $M_{\Sigma,V}/C$ (as a Poisson manifold).
 \end{thm}


We now describe the Poisson submanifold $\hat\mu^{-1}(\mbf{1})/\on{Stab}(\mbf{1})\subseteq M_{\Sigma,V}/C$ as a moduli space. 
For simplicity, we assume\footnote{These are not necessary assumptions, and we encourage the interested reader to describe the generalization.} that $s\in S^2(\g)$ is non-degenerate, and each $C_v\subseteq G$ is Lagrangian (i.e. $C^\perp_v= C_v$ for each $v\in V$). 

Let $\hat\Sigma^c$ be the surface $\hat\Sigma$ with each arc $a_v$ contracted to a point, as in Figure~\ref{fig:contracted}. The surface $\hat\Sigma^c$ is thus obtained from $\Sigma$ by blowing up the marked points and contracting the original boundary arcs. We identify every domain wall $\mbf{w}\subseteq \partial\hat \Sigma$ with its image in $\partial\hat \Sigma^c$. Note that adjacent domain walls intersect non-trivially in $\hat\Sigma^c$.

\begin{defn}\label{def:CompRed}
Suppose $P\to \hat\Sigma^c$ is a principal $G$-bundle, and that $\mbf{w},\mbf{w}'\subseteq \partial\hat\Sigma^c$ are two domain walls. We say that a reduction of the structure group from $G$ to $C_{\mbf{w}}$ over $\mbf{w}$ is compatible with a reduction of the structure group from $G$ to $C_{\mbf{w}'}$ over $\mbf{w}'$ if there exists a common reduction of the structure group from $G$ to $C_{\mbf{w}}\cap C_{\mbf{w}'}$ over the intersection $\mbf{w}\cap \mbf{w}'$ (if it exists, the common reduction of structure is unique up to a unique isomorphism).\footnote{Equivalently, the corresponding $C_\mbf{w}$-subbundle  $Q_{\mbf{w}}\subset P\rvert_{\mbf{w}}$ and the corresponding $C_{\mbf{w}'}$-subbundle $Q_{\mbf{w}'}\subset P\rvert_{\mbf{w}'}$ intersect non-trivially over every point in the intersection $\mbf{w}\cap \mbf{w}'$. Moreover, in the special case that $\mbf{w}=\mbf{w}_v=\mbf{w}'$ and the image of $\mbf{w}_v$ in $\hat\Sigma^c$ forms a closed loop, we insist (naturally) that the two reductions of the structure from $G$ to $C_{\mbf{w}_v}$ over the common endpoint $q(v_+)=q(v_-)$ agree.} 
\end{defn}

We have:

\begin{thm}\label{thm:CompCornData}
Suppose that the graph of $\hat\mu$ composes cleanly with $\mbf{1}\in G^V$, and the $\on{Stab}(\mbf{1})$-orbits of $\hat\mu^{-1}(\mbf{1})$ form a regular foliation, then the moduli space of pairs 
\begin{equation}\label{eq:CompPairs}(P\to\hat\Sigma^c,\{Q_{\mbf{w}}\to \mbf{w}\}_{\mbf{w}\in \mbf{Wall}}),\end{equation}
consisting of a flat $G$-bundles over $\hat\Sigma^c$ equipped with mutually compatible flat reductions of the structure group from $G$ to $C_\mbf{w}$ over each domain wall $\mbf{w}\in\mbf{Wall}$, is a Poisson manifold canonically isomorphic to $\hat\mu^{-1}(\mbf{1})/\on{Stab}(\mbf{1})\subseteq M_{\Sigma,V}/C$.

\end{thm}
\begin{rem}
The moduli space described in Theorem~\ref{thm:CompCornData} is a symplectic manifold (cf. \cite{LiBland:2013ue}).
\end{rem}

\section{Surfaces with domain walls}
Consider a finite family $(\Sigma_d, V_d)$ indexed by a set, $\mbf{Dom}$.
Once again, we construct $\hat\Sigma_d$ by blowing up $\Sigma_d$ at every $v\in V_d$. We let $\mbf{w}_v\subset \hat\Sigma_d$ denote the preimage of $v\in V_d$ and denote the boundary points  of $\mbf{w}_v$ by $v_-$ and $v_+$ (see \cref{fig:Sigmahat}). We let $$\widetilde{\mbf{Wall}}=\underset{{d\in\mbf{Dom},v\in V_d}}{\bigsqcup} \mbf{w}_v,$$ and choose an orientation for each connected component of $\widetilde{\mbf{Wall}}$, defining the map $\on{Sign}:\pi_0(\widetilde{\mbf{Wall}})\to \{-,+\}$ by
$$\on{Sign}(\mbf{w}):=\begin{cases}
+&\text{if }\mbf{w}\subseteq\Sigma_d,\text{ and the inclusion }\\&\mbf{w}\to \partial \Sigma_d\text{ is orientation preserving},\\
-&\text{otherwise,}
\end{cases}$$
for any $\mbf{w}\in \pi_0(\widetilde{\mbf{Wall}})$. 

Now choose an orientation preserving action of $\mbb{Z}_2$ on $\widetilde{\mbf{Wall}}$, which we extend trivially to an action on $\sqcup_{d\in\mbf{Dom}}\hat\Sigma_d$. We let $$\Sigma =  \left(\sqcup_{d\in\mbf{Dom}}\hat\Sigma_d\right)/\mbb{Z}_2,$$
and $\mbf{Wall}=\pi_0(\widetilde{\mbf{Wall}})/\mbb{Z}_2$. As before, we refer to each connected component $\mbf{w}\in \mbf{Wall}$ as a \emph{domain wall}. 

For each $d\in\mbf{Dom}$ choose a Lie group $G_d$ and an element $s_d\in (S^2\g_d)^{G_d}$ - we describe this choice as coloring the domain $\Sigma_d$ with $(G_d,s_d)$. Then
$$M=\prod_d M_{\Sigma_d,V_d}(G_d)$$
is a quasi-Poisson $G$-manifold, where
$$G=\prod_d G_d^{V_d}.$$

Let us now construct a reducing subgroup of $G$ in the following way. Taking $V=\sqcup_{d\in\mbf{Dom}}V_d$, we have a natural identification  $V\cong  \pi_0(\widetilde{\mbf{Wall}})$. The orbits of the induced action of $\mbb{Z}_2$ on $V$ are either pairs (orbit of cardinality 2), or singletons (orbits of cardinality 1), and the orbit space is naturally identified with the set of domain walls $V/\mbb{Z}_2\cong \mbf{Wall}$.

 For every singleton-orbit $\{v\}$, for which $v\in V_d\subseteq V$, we choose a subgroup 
$$C_{\{v\}}\subseteq G_d$$
which is coisotropic with respect to $s_d$.
The singleton-orbit, $\{v\}\in V/\mbb{Z}_2\cong \mbf{Wall}$, is naturally identified with a domain wall $\mbf{w}\in\mbf{Wall}$, and we describe the choice of $C_{\{v\}}$ as \emph{coloring the domain wall} $\mbf{w}$ with $C_\mbf{w}:=C_{\{v\}}$. For future use we define $s_\mbf{w}=\on{Sign}(v)s_d$.

For every pair-orbit $\mc{O}=\{v,v'\}$, for which $v\in V_d$, $v'\in V_{d'}$ we choose a subgroup 
$$C_{\{v,v'\}}\subseteq G_d\times G_{d'}$$
which is coisotropic with respect to 
$\on{Sign}(v)s_{d}+\on{Sign}(v')s_{d'}$.
The pair-orbit, $\{v,v'\}\in V/\mbb{Z}_2\cong \mbf{Wall}$, is naturally identified with a domain wall $\mbf{w}\in\mbf{Wall}$, and we describe the choice of $C_{\{v,v'\}}$ as \emph{coloring the domain wall} $\mbf{w}$ with $C_\mbf{w}:=C_{\{v,v'\}}$. As before, we define $s_\mbf{w}=\on{Sign}(v)s_{d}+\on{Sign}(v')s_{d'}$

Following Wehrheim and Woodward \cite{Wehrheim:2010cg,Wehrheim:2010fa} we shall call the resulting surface, $\Sigma$, together with the chosen coloring of each domain $\Sigma_d$ with $(G_d,s_d)$, and each domain wall $\mbf{w}\in\mbf{Wall}$ with the coisotropic subgroup $C_\mbf{w}$, a \emph{quilted surface}.

  The product
$$C=\prod_{\mbf{w}\in\mbf{Wall}}C_\mbf{w}=\left[\left(\prod_{\{v\}\in V/\mbb{Z}_2}C_v\right)\times \left(\prod_{\{v,v'\}\in V/\mbb{Z}_2}C_{\{v,v'\}}\right)\right]\subseteq \prod_{d\in\mbf{Dom}} G_d^{V_d}$$
is a reducing subgroup of $G$. Indeed, it is coisotropic with respect to $$s:=\sum_{d\in\mbf{Dom},v\in V_d}\on{Sign}(v)s_d,$$ as in \cref{rem:CanNegates}.

As a result (by Theorem \ref{thm:Red}), 
$$M/C=\big(\prod_d M_{\Sigma_d,V_d}(G_d)\big)/C$$
is a Poisson manifold (provided the equivalence relation is regular). The manifold $M/C$ can be interpreted as a moduli space classifying compatible flat bundles over the quilted surface; for brevity, we will not describe the corresponding moduli problem in detail (as we did in Theorem~\ref{thm:FullModSpc} for a single domain), though we encourage the interested reader to do so. Instead, we will now use the ``hat'' construction from Section~\ref{sec:HatConstr} to give an explicit geometric description of $M/C$ in terms of holonomy representations of the fundamental groupoid (as would be obtained from such flat bundles).


We begin by describing the induced space $\hat M$ corresponding to $M=\prod_d M_{\Sigma_d,V_d}(G_d)$, the action of $\hat C$ on $\hat M$, and the map $\hat\mu:\hat M\to G$ in more detail. 
As before, we let $\mf{c}_\mbf{w}^\perp=s_\mbf{w}^\sharp \ann(\mf{c}_\mbf{w})$, where $\mf{c}_\mbf{w}$ is the Lie algebra corresponding to $C_\mbf{w}$. Once again, $\mf{c}_\mbf{w}^\perp\subseteq \mf{c}_\mbf{w}$ is an ideal, and we let $C_\mbf{w}^\perp\subseteq C_\mbf{w}$ be the corresponding connected Lie group (which is normalized by $C_\mbf{w}$). We take $$\hat C_\mbf{w}= \{(g_{\mbf{w}_-},g_{\mbf{w}_+})\in C_{\mbf w}\times C_{\mbf w}\mid\ g_{\mbf{w}_+}g_{\mbf{w}_-}^{-1}\in C_{\mbf w}^\perp\}.$$
We have $$\hat C=\prod_{\mbf{w}\in \mbf{Wall}}\hat C_\mbf{w}\subseteq \prod_{\mbf{w}\in\mbf{Wall}}C_\mbf{w}\times C_\mbf{w}= \prod_{d\in \mbf{Dom}}G^{V_d}\times G^{V_d}.$$


\begin{figure}[h!]
\centering
\def\svgwidth{7cm}
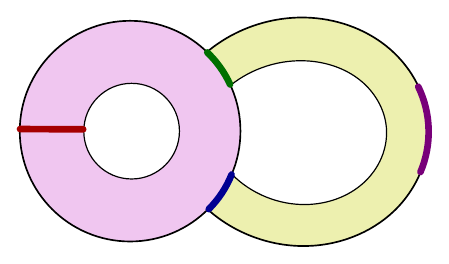
\caption{A quilted surface with 2 domains and 4 domain walls}
\end{figure}

Now we relate the induced space $\hat M$ to the the quilted surface. An element $f\in\hat M$ is a collection of  elements $f_d\in M_{\hat\Sigma_d, {V_d}_-\cup {V_d}_+}(G_d)$ ($d\in\mbf{Dom}$), satisfying the condition:
$$f(\mbf w)\in C^\perp_\mbf{w},$$
for every domain wall $\mbf{w}\in\mbf{Wall}$.
Here $f(\mbf w)$ is the holonomy (if $\mbf w$ is on the boundary of $\Sigma$) or the pair of holonomies (if $\mbf w$ is inside $\Sigma$) along $\mbf w$.
%
Every element $(g_{\mbf{w}_-},g_{\mbf{w}_+})\in\hat C_{\mbf w}$ acts on $\hat M$ by acting at the initial and final boundary points of $\mbf{w}$,  ($\mbf{w}_-$, and $\mbf{w}_+$, respectively).\footnote{Recall that the orientation of $\mbf{w}$ as a domain wall was chosen arbitrarily, and may not coincide with the orientation of the boundaries of the domains. For example,
 if $\mbf{w}\in \mbf{Wall}\cong V/\mbb{Z}_2$ corresponds to the singleton-orbit $\{v\}\in V_d\subseteq V$, then we have $$(\mbf{w}_-,\mbf{w}_+)=\begin{cases}
(v_-,v_+)&\text{ if }\on{Sign}(v)=+,\text{ and }\\
(v_+,v_-)&\text{ if }\on{Sign}(v)=-.\\
\end{cases}$$}
This defines the action of $\hat C$ on $\hat M$.

Finally, the map 
$$\hat\mu:\hat M\to G=\prod_{d\in\mbf{Dom}}G_d^{V_d}$$
is the collection of holonomies along the arcs $a_v$, $v\in V_d$ (i.e.\ along the boundary arcs of $\Sigma$ which are not domain walls).

The isomorphism
$$\hat M/\hat C\cong M/C$$
can be seen via the embedding $M\hookrightarrow\hat M$ given by the constraint $f(\mbf w)=\mbf{1}$ for all $\mbf w \in\mbf{Wall}$ (effectively contracting every $\mbf w$ to a point).
We write
\begin{equation}\label{eq:ModSpcQuilSurf}\mc{M}_\Sigma:=\hat M/\hat C\end{equation}

Theorem~\ref{thm:mmred} and Proposition~\ref{prop:hatclean} now give the following result.

\begin{subtheorem}{thm}\label{thm:quiltedAll}
\begin{thm}\label{thm:quilted}
Let $h\in G$. If $h$ is in a clean position relative to $\hat\mu$ and the $\on{Stab(h)}$-orbits of $\hat\mu^{-1}(h)$ form a regular equivalence relation, then the quotient space
\begin{equation}\label{eq:quilted}
N:=\{f\in\hat M;\ f(a_v)=h_v\}/\on{Stab(h)}
\end{equation}
is a Poisson manifold.

Moreover, if the $C$-orbits of $M$ form a regular equivalence relation, then $N$ is naturally isomorphic to a Poisson submanifold of $M/C$.
\end{thm}

The theorem is particularly interesting when $h=\mbf{1}$, as then we can contract the arcs $a_v$. 
More generally we might want to only contract a subset of the arcs $a_v$. For instance, suppose we choose a subset $V_{triv}\subseteq V$ and want to contract those arcs $a_v$ with $v\in V_{triv}$.
 We can do this as follows: Let $$H=\{h\in \prod_{d\in\mbf{Dom}}G_d^{V_d};\  h_{v}=\mbf{1}\text{ for all } v\in V_{triv}\},$$
and 
$$\hat H=\{g\in \prod_{d\in\mbf{Dom}}G_d^{{V_d}_-\cup {V_d}_+};\  g_{v_+}=g_{\sigma( v)_-}\text{ for all } v\in V_{triv}\}.$$
Then for any $g\in \prod_{d}G_d^{{V_d}_-\cup {V_d}_+}$, and $h\in H$, 
\begin{equation}\label{eq:ginhatH}g\cdot h\in H\Rightarrow g\in \hat H,\end{equation}
where the action on the left hand side is the twisted action \eqref{eq:tautwist}.
Thus the stabilizer groupoid for $H$ is $\on{Stab}(H)=\left(\hat C\cap\hat H\right)\ltimes H$ (cf. Remark~\ref{rem:Stab}).

Let $\Sigma_{V_{triv}}$ be the surface obtained from $\Sigma$ by contracting each arc $a_v\subseteq \Sigma$ (for $v\in V_{triv}$) to a point.

Thus Proposition~\ref{prop:hatclean} implies:
\begin{thm}\label{cor:PartTriv} If $H$ is in a clean position relative to $\hat \mu$, $S=\hat C \cdot H$ is a submanifold of $\prod_{d}G_d^{V_d}$, and the $\hat C\cap \hat H$-action on $\hat\mu^{-1}(H)$ defines a regular equivalence relation, then the moduli space
$$\mc{M}_{\Sigma_{V_{triv}}}:=\hat\mu^{-1}(H)/(\hat C\cap \hat H)$$ is
a Poisson manifold. 

Moreover, if the action of $C$ on $M$ defines a regular equivalence relation, then  $\mc{M}_{\Sigma_{V_{triv}}}$ is
naturally isomorphic to a Poisson submanifold of $M/C$.
\end{thm}
\end{subtheorem}
Explicitly, we have 
\begin{multline*}\hat\mu^{-1}(H)=\{f\in\hat M;\ f(a_v)=\mbf{1}\text{ for all }v\in V_{triv}\}
\\=\left\{\left.f\in \prod_{d\in\mbf{Dom}} M_{\hat\Sigma_d, {V_d}_-\cup {V_d}_+}(G_d)\right\rvert \begin{split}f(\mbf{w})&\in C^\perp_\mbf{w}\text{ for every } \mbf{w}\in\mbf{Wall},\text{ and }\\
 f(a_v)&=1\text{ for all }v\in V_{triv}
 \end{split}
\right\}
\end{multline*}

In the sequel, we will refer to 
$$\hat C\cap\hat H = \left\{g\in \prod_d G_d^{{V_d}_-\cup {V_d}_+}\left\lvert \begin{split}
 g_{\mbf{w}_-},g_{\mbf{w}_+}&\in C_{\mbf{w}}\\
 g_{\mbf{w}_+}g_{\mbf{w}_-}^{-1}&\in C_{\mbf{w}}^\perp\\
 g_{v_+}&=g_{\sigma(v)_-} \text{ for all } v\in V_{triv}
 \end{split}\right.\right\}$$
 as the \emph{group of residual gauge transformations}. Here
  $g_{\mbf{w}}:=g_v$ when the domain wall $\mbf{w}\in \mbf{Wall}\cong V/\mbb{Z}_2$ corresponds to the singleton-orbit $\{v\}\subset V$, and $g_{\mbf{w}}:=(g_v,g_{v'})$ when the domain wall $\mbf{w}$ corresponds to the pair-orbit $\{v,v'\}\subset V$.

We leave the reformulation of Theorem \ref{thm:quiltedAll} in terms of moduli spaces of flat bundles (in the spirit of Theorems~\ref{thm:FullModSpc} and \ref{thm:CompCornData}) to the reader.

\begin{rem}[Quilted surfaces with residual marked points]\label{rmk:qPoisQuilt}
	One may also leave the marked points on certain boundary components of the domains $\Sigma_d$ ($d\in\mbf{Dom}$) unreduced (or unsewn); the resulting moduli spaces are quasi-Poisson rather than Poisson. For example, suppose that we decompose the set $B:=\sqcup_{d}\partial\Sigma_d$ of all boundary components into two closed subsets $B=B^{sew}\sqcup B^{res}$. We will leave those marked points $V^{res}:=V\cap B^{res}$ in $B^{res}$ unreduced (here $V=\sqcup_d V_d$). 
	
	As before, le $\hat\Sigma_d$ denote the surface obtained by blowing up $\Sigma_d$ at $V_d^{sew}:=V_d\cap B^{sew}$, $\widetilde{\mbf{Wall}}$ the union of the exceptional divisors and $\widetilde{\mbf{Wall}}\to \mbf{Wall}$ the quotient obtained by sewing along the domain walls. 
	
	We color each domain, $\Sigma_d$ with a pair $(G_d,s_d)$, and each domain wall, $\mbf{w}\in\mbf{Wall}$, with a coisotropic subgroup $C_\mbf{w}$. We take
\begin{align*}
G&=\prod_{d\in\mbf{Dom}}G_d^{V_d^{sew}}&H&=\prod_{d\in\mbf{Dom}}G_d^{V_d\cap V^{res}},&\text{ and }C&=\prod_{\mbf{w}\in\mbf{Wall}}C_{\mbf{w}}\subseteq G
\end{align*}

Once again, we let $\hat M$ denote the induced manifold, 
$$\hat M = \bigg\{f\in \prod_d M_{\hat\Sigma_d,{V_d}^{res}\cup  {V_d}_-^{sew}\cup {V_d}_+^{sew}}(G_d)\bigg\lvert f(\mbf{w})\in C^\perp_\mbf{w}\text{ for each }\mbf{w}\in \mbf{Wall}\bigg\}$$
 
 Then Proposition~\ref{prop:hatclean} implies that the moduli space
 $$\{f\in\hat M;\ f(a_v)=\mbf{1}\text{ for each }v\in V^{sew}\}/\on{Stab}(\mbf{1})$$ is a quasi-Poisson $H$-manifold, where $$\on{Stab}(\mbf{1})=\{\hat c\in \hat C;\  \hat c\cdot \mbf{1}=\mbf{1}\}.$$
\end{rem}

\section{Spin networks}
In this section, we reinterpret the quasi-Poisson structure on marked surfaces as well as the Poisson structure on quilted surfaces in terms of spin networks. 
\begin{rem}
Spin networks were first introduced by Penrose \cite{Penrose:1971vl} (see also \cite{Baez:1996hh}). Poisson brackets of spin networks were studied by Roche and Szenes \cite{Roche:2000ws}, following work by Goldman \cite{Goldman:1986eh,Goldman:1984hr},  Andersen, Mattes and Reshetikhin \cite{Andersen:1996ur}.
\end{rem}
\subsection{Spin networks on a marked surface}
Let $(\Sigma,V)$ denote a marked surface, as in \cref{sec:QPsMs}.
\begin{defn}
A \emph{graph diagram in $(\Sigma,V)$} consists of the following data:
\begin{itemize}
\item a directed graph, $\Gamma$ with edges $E_\Gamma$ and vertices $V_\Gamma$,
\item a subset of vertices $V_\Gamma^{anch}\subseteq V_\Gamma$, called `anchor points', and
\item a smooth map $\iota:\Gamma\to \Sigma$ (it is understood that the domain of this map is really the geometric realization, $\lvert\Gamma\rvert$, of $\Gamma$) such that 
$$\iota^{-1}(V)=V_\Gamma^{anch}.$$
\end{itemize}
Often we will abuse notation and denote the graph diagram simply by $\Gamma$.

A morphism between graph diagrams $\iota:\Gamma\to\Sigma$ and $\iota':\Gamma'\to\Sigma$ is a map between the geometric realizations, $\mu:\lvert\Gamma\rvert\to\lvert\Gamma'\rvert$, such that $\iota=\iota'\circ\mu$ and $\mu(V_\Gamma)\subseteq V_{\Gamma'}$.

A homotopy between graph diagrams $\iota_0:\Gamma\to\Sigma$ and $\iota_1:\Gamma\to\Sigma$ is a one parameter family of graph diagrams $\iota_t:\Gamma\to\Sigma$, $t\in[0,1]$.
\end{defn}
Suppose $\Gamma$ is a graph diagram in $(\Sigma,V)$. The Lie group $G^{V_\Gamma}$ acts on $G^{E_\Gamma}$ by
$$(g\cdot g')_e=g_{\mbf{t}(e)}g'_e g_{\mbf{s}(e)}^{-1},$$
where $g\in G^{V_\Gamma}$, $g'\in G^{E_\Gamma}$, and $e\in E_\Gamma$.
We define
$$\mc{A}_\Gamma(G):= C^\infty\big(G^{E_\Gamma}\big)^{G^{V_{\Gamma}^{int}}},$$
where $V_\Gamma^{int}=V_\Gamma\smallsetminus V_\Gamma^{anch}$.
Note that there is a residual action of $G^{V_\Gamma^{anch}}$ on $\mc{A}_\Gamma(G)$. For any $v\in V_{\Gamma}^{anch}$, we let $\rho_v:G\times \mc{A}_\Gamma(G)\to \mc{A}_\Gamma(G)$ denote the action of the $v^{th}$-factor.
\begin{rem}\label{rem:HomSpaceDescr}
Since $\Gamma$ is a graph, $G^{E_\Gamma}=\on{Hom}\big(\Pi_1(\lvert\Gamma\rvert,V_\Gamma),G\big)$. Therefore
$$\mc{A}_\Gamma(G)= C^\infty\bigg(\on{Hom}\big(\Pi_1(\lvert\Gamma\rvert,V_\Gamma),G\big)\bigg)^{G^{V_{\Gamma}^{int}}}.$$
\end{rem}

By functoriality, a morphism $\mu:\lvert\Gamma\rvert\to\lvert\Gamma'\rvert$ of graph diagrams defines a morphism 
$$\mu^!:\on{Hom}\big(\Pi_1(\lvert\Gamma'\rvert,V_{\Gamma'}),G\big)\to \on{Hom}\big(\Pi_1(\lvert\Gamma\rvert,V_\Gamma),G\big),$$
which is equivariant with respect to the map $$\mu_V^!:G^{V_{\Gamma'}}\to G^{V_{\Gamma}}$$ defined by
$\mu^!(g)_v=g_{\mu(v)}$, for any $g\in G^{V_{\Gamma'}}$. 
Since $\mu_V^!$ maps the normal subgroup $G^{V_{\Gamma'}^{int}}$ into $G^{V_{\Gamma}^{int}}$, pull-back along $\mu^!$ defines an equivariant morphism of algebras
$\mu_*:\mc{A}_\Gamma\to\mc{A}_{\Gamma'}$. We may summarize this as:

\begin{lem}

The assignment $(\iota:\Gamma\to \Sigma) \mapsto \mc{A}_{\Gamma}$ is a functor from graph diagrams to algebras endowed with an action of $G^V$, where $g\in G^V$ acts via
$$g\cdot f= \big(\prod_{v\in V_\Gamma^{anch}}\rho_v(g_{\iota(v)})\big)f.$$
\end{lem}

\begin{defn}
A \emph{spin network} in $(\Sigma,V)$ is a pair $(\Gamma,f)$, where $\Gamma\to\Sigma$ is a graph diagram, and $f\in\mc{A}_\Gamma(G)$.

We say that spin networks $(\Gamma_0,f_0)$ and $(\Gamma_1,f_1)$ are \emph{homotopic} if the underlying graph diagrams are homotopic and $f_0=f_1\in \mc{A}_{\Gamma_0}\equiv\mc{A}_{\Gamma_1}$ (to identify the algebras, we use the fact that the definition of $\mc{A}_\Gamma$ only depends on the underlying graph, and not the map $\Gamma\to \Sigma$).

A morphism of spin networks $\mu:(\Gamma,f)\to (\Gamma',f')$ is a morphism of graph diagrams $\mu:\Gamma\to\Gamma'$ such that $f'=\mu_*f$. 
\end{defn}
We consider two spin networks to be equivalent if they are related by a chain of homotopies and morphisms (or the formal inverses of morphisms). For a spin network $(\Gamma,f)$, we let $[\Gamma,f]$ denote the corresponding equivalence class. Define $\on{SpinNet}_{(\Sigma,V)}(G)$ to be the set of equivalence classes of spin networks. 

\begin{lem}\label{lem:AlgofSpinNet}
$\on{SpinNet}_{(\Sigma,V)}(G)$ is a $G^V$-algebra where scalar multiplication, addition, and multiplication are defined as
\begin{subequations}\label[pluralequation]{eq:SpinNetOp}
\begin{align}
\lambda\cdot[\Gamma,f]&=[\Gamma,\lambda f]\\
[\Gamma,f]+[\Gamma',f']&=[\Gamma\cup\Gamma',f\oplus f']\\
[\Gamma,f]\cdot[\Gamma',f']&=[\Gamma\cup\Gamma',f\otimes f']
\end{align}
\end{subequations}
\end{lem}
\cref{lem:AlgofSpinNet} will follow from \cref{prop:EvAlgId}.

\subsection{Spin networks and functions on the moduli space}\label{sec:SpNetMrSFunct}
Suppose $(\Gamma,f)$ is a spin network in $(\Sigma,V)$. Then the we may push $f$ along the map $\iota:\Gamma\to\Sigma$ to define a function $\on{ev}(\Gamma,f)$ on the moduli space $M_{\Sigma,V}(G)$, as we shall now describe.

For a finite set of points $X=\{x_i\}\in \Sigma\smallsetminus V$, we let 
$$M_{\Sigma,V,X}(G):=\Hom\big(\Pi_1(\Sigma,V\cup X),G\big)$$
Notice that 
\begin{equation*}\label{eq:AddExtraVert}M_{\Sigma,V}(G)=M_{\Sigma,V,X}(G)/G^X.\end{equation*}
The map $$\iota:(\lvert\Gamma\rvert,V_{\Gamma})\to (\Sigma,V\cup V_{\Gamma}^{int})$$ yields a $G^V\times G^{V_\Gamma^{int}}$-equivariant map
\begin{equation}\label{eq:tauEx}\iota^!:M_{\Sigma,V,V_\Gamma^{int}}(G)\to \Hom\big(\Pi_1(\lvert\Gamma\rvert,V_{\Gamma}),G\big).\end{equation}

Therefore, the function 
$$\on{ev}(\Gamma,f):={\iota^!}^*f\in C^\infty\big(M_{\Sigma,V,V_\Gamma^{int}}(G)\big)$$ 
is $G^{V_\Gamma^{int}}$-invariant. Hence it descends to define a function on the moduli space
$$M_{\Sigma,V}(G)\cong M_{\Sigma,V,V_\Gamma^{int}}(G)/G^{V_\Gamma^{int}}.$$
Moreover, the map $$\on{ev}(\Gamma,\cdot):\mc{A}_\Gamma(G)\to C^\infty(M_{\Sigma,V}(G))$$ is $G^V$-equivariant.

Notice that $\on{ev}(\Gamma_0,f_0)=\on{ev}(\Gamma_1,f_1)$ whenever $(\Gamma_0,f_0)$ and $(\Gamma_1,f_1)$ are homotopic. Moreover, if 
$$\mu:(\iota:\lvert\Gamma\rvert\to\Sigma)\to(\iota':\lvert\Gamma'\rvert\to\Sigma)$$ is a morphism of graph diagrams, then the following diagram of equivariant morphisms commute
$$\begin{tikzpicture}
\mmat{m}{M_{\Sigma,V,V_{\Gamma'}^{int}}(G)&\\ 
\Hom\big(\Pi_1(\lvert\Gamma'\rvert,V_{\Gamma'}),G\big)&\Hom\big(\Pi_1(\lvert\Gamma\rvert,V_{\Gamma}),G\big)\\};
\draw[->] (m-1-1) edge node {$\iota'^!$} (m-2-1);
\path[->] (m-1-1) edge node {$\iota^!$} (m-2-2);
\draw[->] (m-2-1) edge node[swap] {$\mu^!$} (m-2-2);
\end{tikzpicture}$$
Therefore,  ${\iota^!}^*( f)={\iota'^!}^*\circ\mu_*(f)$ for any $f\in \mc{A}_\Gamma$, i.e.
$\on{ev}(\Gamma',\mu_*f)=\on{ev}(\Gamma,f)$. So $\on{ev}$ descends to a map on the set of equivalence classes of spin networks, $\on{SpinNet}_{\Sigma,V}(G)$.

\begin{prop}\label{prop:EvAlgId}
The map $$\on{ev}:\on{SpinNet}_{\Sigma,V}(G)\to C^\infty\big(M_{\Sigma,V}(G)\big)$$ is an isomorphism of $G^V$-algebras.
\end{prop}
\begin{proof}
With only a slight modification, the statement follows from the proofs of the corresponding statements for (unmarked) surfaces in \cite{Baez:1996hh,Roche:2000ws}, but we outline it here for completeness.

Let $\widetilde{\on{SpinNet}_{\Sigma,V}(G)}$ be the set of all spin networks in $(\Sigma,V)$ (i.e. we do not identify equivalent spin networks). 
We may define the operations in \cref{eq:SpinNetOp} directly on $\widetilde{\on{SpinNet}_{\Sigma,V}(G)}$ (we do not claim that they satisfy the axioms of an algebra on this set).
Since $(\Gamma,f)\to \on{ev}(\Gamma,f)$ is induced by the map of spaces \labelcref{eq:tauEx}, it follows that
$\on{ev}$ intertwines the operations of scalar multiplication, addition, and multiplication. 

Next we show that $\on{ev}:\widetilde{\on{SpinNet}_{\Sigma,V}(G)}\to C^\infty\big(M_{\Sigma,V}(G)\big)$ is surjective. Let $\iota_{skel}:\lvert\Gamma_{skel}\rvert\to \Sigma$ be an embedded graph diagram with a single anchor point at every marked point, for which there exists a deformation retract $r:\Sigma\to \iota_{skel}(\lvert\Gamma_{skel}\rvert)$. 
Then the map defined in \cref{eq:tauEx},
$$\Hom\big(\Pi_1(\Sigma,V\cup V_\Gamma^{int}),G\big)\to \Hom\big(\Pi_1(\lvert\Gamma_{skel}\rvert,V_{\Gamma_{skel}}),G\big)$$
 is a diffeomorphism. It follows that $\on{ev}:\mc{A}_{\Gamma_{skel}}\to C^\infty\big(M_{\Sigma,V}(G)\big)$ is an isomorphism.

It remains to show that that $\on{ev}(\Gamma,f)=\on{ev}(\Gamma',f')$ only if the two spin networks are equivalent. We may assume that there exist maps between the geometric realizations $\mu:\lvert\Gamma\rvert\to\lvert\Gamma_{skel}\rvert$ (just compose the map $\lvert\Gamma\rvert\to\Sigma$ with the retract $\Sigma\to \lvert\Gamma_{skel}\rvert$) and likewise $\mu':\lvert\Gamma'\rvert\to\lvert\Gamma_{skel}\rvert$. Since $\on{ev}:\mc{A}_{\Gamma_{skel}}\to C^\infty\big(M_{\Sigma,V}(G)\big)$ is an isomorphism, $\on{ev}(\Gamma,f)=\on{ev}(\Gamma',f')$ only if $\mu_*f=\mu'_*f'$, i.e. $(\Gamma,f)$ and $(\Gamma',f')$ are equivalent.
\end{proof}

\subsection{The quasi-Poisson bracket on $\on{SpinNet}_{\Sigma,V}(G)$}\label{sec:QPoisBrakSpin}

\begin{figure}[h!]
\begin{center}
\begin{subfigure}[t]{.45\linewidth}
\centering
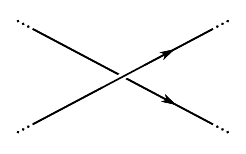
\caption{The edge $e\in E_\Gamma$ and $e'\in E_{\Gamma'}$ intersect at $A$.}
\end{subfigure}
\begin{subfigure}[t]{.45\linewidth}
\centering
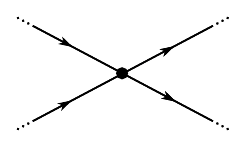
\caption{After subdividing the edges $e$ and $e'$ at $A$ and then the newly created vertices, the newly created edges are labelled as depicted.}
\end{subfigure}
\end{center}
\caption{\label{fig:GraphJoin}}
\end{figure}

Suppose that $\iota:\Gamma\to \Sigma$ and $\iota':\Gamma'\to \Sigma$ are graph diagrams such that the maps $\iota$ and $\iota'$ are transverse to each other.\footnote{That is, the restriction of $\iota$ and $\iota'$ to any two edges are transverse, and $\iota(V_\Gamma^{int})\cap\iota'(V_{\Gamma'}^{int})=\emptyset.$} Let $(\Gamma\times_\Sigma\Gamma')^{anch}=(V_\Gamma^{anch})\times_{\Sigma} (V_{\Gamma'}^{anch})$ and 
$$(\Gamma\times_\Sigma\Gamma')^{int}=\big(\Gamma\times_\Sigma\Gamma'\big)\smallsetminus (\Gamma\times_\Sigma\Gamma')^{anch}.$$
Given $A\in (\Gamma\times_\Sigma\Gamma')^{int}$, let $e\in E_\Gamma$ and $e'\in E_{\Gamma'}$ be the two edges intersecting at $A$.
Let $\Gamma\cup_{A}\Gamma'$ be the graph obtained from $\Gamma\cup\Gamma'$ by subdividing the edges $e$ and $e'$ at $A$ and then merging the newly created vertices. It is clear that $$\iota\cup\iota':\Gamma\cup_{A}\Gamma'\to\Sigma$$ is a graph diagram.
 We let $e_A$ and $e'_A$ denote the newly created edges running from $\mbf{s}(e)$ and $\mbf{s}(e')$ (respectively) to $A$ and 
  $\tilde e_A$ and $\tilde e'_A$ denote the newly created edges running from $A$ to $\mbf{t}(e)$ and $\mbf{t}(e')$ (respectively), as depicted in \cref{fig:GraphJoin}. Note that the orientations of the new edges are inherited from the orientations of the original edge.

We define the map
 $$\varrho:(\g\times\g)\times G^{E_{(\Gamma\cup_{A}\Gamma')}}\to T\big(G^{E_{\Gamma}}\big)\times T\big(G^{E_{\Gamma'}}\big)$$ by
 $$\varrho(\xi,\eta,g)_{e_0}:=\begin{cases}
(L_{g_{\tilde e_A}}R_{g_{e_A}})_*\xi &\text{if } e_0=e,\\
(L_{g_{\tilde e'_A}}R_{g_{e'_A}})_*\eta &\text{if } e_0=e',\\
g_{e_0}&\text{otherwise.}
\end{cases}$$
Where $L_g,R_g:G\to G$ denote Left,Right multiplication by $g\in G$, and we have identified $\g$ with the tangent space at the identity.
The map $\varrho$ is $G^{V_{(\Gamma\cup_{A}\Gamma')}}$-equivariant, where
the action is defined on $\g\times\g$ by
$$g\cdot(\xi,\eta):=(\Ad_{g_A}\xi,\Ad_{g_A}\eta),\quad g\in G^{V_{(\Gamma\cup_{A}\Gamma')}},\quad \xi,\eta\in\g.$$
The universal property of the tensor product implies that $\varrho$ extends to a map
 $$\varrho:(\g\otimes\g)\times G^{E_{(\Gamma\cup_{A}\Gamma')}}\to T\big(G^{E_{\Gamma}}\big)\otimes T\big(G^{E_{\Gamma'}}\big)\subseteq \otimes^2 T \big( G^{E_{(\Gamma\cup\Gamma')}}\big)$$
  We define 
 $$\Psi_s:=\varrho(s,\cdot):G^{E_{(\Gamma\cup_{A}\Gamma')}}\to \otimes^2 T \big( G^{E_{(\Gamma\cup\Gamma')}}\big).$$
Since $s$ is $G$ invariant, it follows that $\Psi_s$ is 
$G^{V_{(\Gamma\cup_{A}\Gamma')}}$-equivariant. 
Therefore 
$$\big(f\otimes f'\to \Psi_s^*(\d f\otimes \d f')\big):\mc{A}_{\Gamma}(G)\otimes \mc{A}_{\Gamma'}(G)\to \mc{A}_{(\Gamma\cup_{A}\Gamma')}(G)$$
is a $G^{V_{\Gamma}^{anch}}\times G^{V_{\Gamma'}^{anch}}$-equivariant linear map.


\begin{prop}\label{prop:qPBrkSpNet}
Let $\on{SpinNet}_{\Sigma,V}(G)$ be endowed with the bracket
\begin{multline}\label{eq:PBrkSpNetMrkPr}\big\{[\Gamma,f],[\Gamma',f']\big\}:=\sum_{A\in(\Gamma\times_\Sigma\Gamma')^{int}}\on{sign}(A)\big[\Gamma\cup_{A}\Gamma',\Psi_s^*(\d f\otimes \d f')\big]\\
+\sum_{A=(v,v')\in (\Gamma\times_\Sigma\Gamma')^{anch}}\frac{1}{2}\on{sign}(A)\big[\Gamma\cup\Gamma',(\d f\otimes \d f')\big((\rho_v\otimes\rho_{v'})s\big)\big],
\end{multline}
where $\on{sign}(\cdot)=\pm1$ is computed as pictured in \cref{fig:Signp}, and $\Gamma$ and $\Gamma'$ are assumed to be transverse graph diagrams.
Then \begin{equation}\label{eq:EvIntBrk}\pi\big(\d\on{ev}(\Gamma,f),\d\on{ev}(\Gamma',f')\big)=\on{ev}\big\{[\Gamma,f],[\Gamma',f']\big\}.\end{equation}
%
\end{prop}
\begin{proof}
Suppose first that $V_\Gamma^{int}=\emptyset= V_{\Gamma'}^{int}$.
Then $\mc{A}_\Gamma(G)=C^\infty(G^{E_\Gamma})$, and 
$\d\on{ev}(\Gamma,f)=\prod_{e\in E_\Gamma}\hol_{e}^* \d f$. 
 Thus
\begin{equation}\label{eq:dSpinNetFunct}\d\on{ev}(\Gamma,f)=\sum_{e\in E_\Gamma}\big((L_{\hol_e})_e^* \d f\big)(\hol_{e}^*\theta^L)\end{equation}
where $(L_{g})_e,(R_{g})_e:G^{E_\Gamma}\to G^{E_\Gamma}$ denotes the left,right multiplication by $g\in G$ on the $e\in E_{\Gamma}$-th factor of $G^{E_\Gamma}$.
Substituting \cref{eq:dSpinNetFunct} (and the corresponding expression for $\d\on{ev}(\Gamma',f')$) into \cref{eq:pi_hpair}  results in the equality
\begin{multline*}
\pi\big(\d\on{ev}(\Gamma,f),\d\on{ev}(\Gamma',f')\big)\\
=\frac{1}{2}\sum_{e\in E_\Gamma,e'\in E_{\Gamma'}} \big((L_{\hol_e})_e^* \d f\big)\otimes  \big((L_{\hol_{e'}})_{e'}^* \d f'\big) (\Ad_{\hol_{(e,e')}}\otimes 1)s
\end{multline*}
Expanding $(e,e')$ in the last line using \cref{eq:IntPairing}, and simplifying the resulting expression using the equalities $\hol_{e}=\hol_{\tilde e_A}\hol_{e_A}$ and $\hol_{e'}=\hol_{\tilde e'_A}\hol_{e'_A}$, yields
\begin{multline*}\sum_{A\in e\cap e', e\in E_\Gamma,e'\in E_{\Gamma'}}\frac{1}{2} \lambda(A)\on{sign}(A)\big( \d f\otimes   \d f'\big)\\\big(\big((L_{\hol_{\tilde e_A}})_e (R_{\hol_{ e_A}})_e\otimes (L_{\hol_{\tilde e'_A}})_{e'} (R_{\hol_{ e'_A}})_{e'}\big)s\big).
\end{multline*}
If $A\in \partial\Sigma$, then $\lambda(A)=1$, and 
$$(L_{\hol_{\tilde e_A}})_e (R_{\hol_{ e_A}})_e=\rho_v,\quad (L_{\hol_{\tilde e'_A}})_e (R_{\hol_{ e'_A}})_e=\rho_{v'},$$
where $A=(v,v')\in (\Gamma\times_\Sigma\Gamma')^{anch}$. On the other hand, if $A\not\in\partial\Sigma$, then $\lambda(A)=2$, and 
$$\big( \d f\otimes   \d f'\big)\big(\big((L_{\hol_{\tilde e_A}})_e (R_{\hol_{ e_A}})_e\otimes (L_{\hol_{\tilde e'_A}})_{e'} (R_{\hol_{ e'_A}})_{e'}\big)s\big)=\Psi_s^*(\d f\otimes \d f').$$
Therefore \cref{eq:EvIntBrk} holds in this case.

\begin{figure}[h!]
\centering
\def\svgwidth{\linewidth}
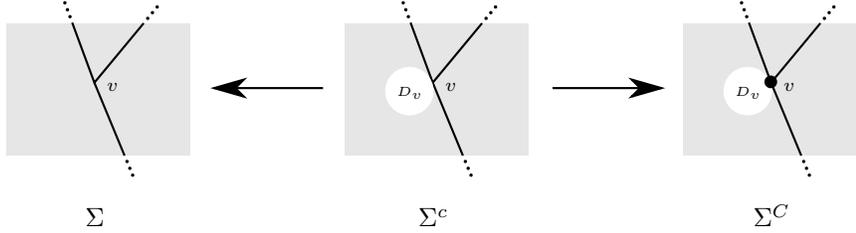
\caption{For each vertex $v$, we delete a small open disk, $D_v$, from $\Sigma$ such that its boundary $\partial D_v$ intersects the graph precisely at $v$.\label{fig:SigmaC}}
\end{figure}

Suppose now that $\Gamma$ and $\Gamma'$ are arbitrary graph diagrams (which are transverse to each other). Let $\Sigma^c$ be the marked surface obtained from $\Sigma$ as follows: for each $v\in V_{\Gamma}^{int}\cup V_{\Gamma'}^{int}$ delete a small open disk $D_v\subset\Sigma\smallsetminus{\Gamma\cup\Gamma'}$ such that $$\partial D_v\cap(\Gamma\cup\Gamma')=v.$$
Let $\Sigma^C$ be the marked surface obtained by $\Sigma^c$ by adding a marked point at $v$ (cf. \cref{fig:SigmaC}). 

If we view $(\Gamma,f)$ and $(\Gamma',f')$ as spin networks in $\Sigma^C$, all their vertices are anchor points. Therefore \cref{eq:PBrkSpNetMrkPr} holds for these spin networks in $\Sigma^C$. Finally, applying \cref{cor:emb} to the embeddings $$(\Sigma^c,V)\to (\Sigma,V)$$ and $$(\Sigma^c,V)\to (\Sigma^C,V\cup V_{\Gamma}^{int}\cup V_{\Gamma'}^{int})$$ shows that \cref{eq:PBrkSpNetMrkPr} holds for these spin networks in $\Sigma$.
\end{proof}

\subsection{Spin networks on quilted surfaces}\label{sec:SpinNet}
Let $\Sigma$ be a quilted surface.
\begin{defn}
A \emph{graph diagram in $\Sigma$} consists of the following data:
\begin{itemize}
\item a directed graph, $\Gamma$ with edges $E_\Gamma$ and vertices $V_\Gamma$,
\item a subset of vertices $V_\Gamma^{anch}\subseteq V_\Gamma$, called `anchor points', and
\item a smooth map $\iota:\Gamma\to \Sigma$ (it is understood that this map is really from the geometric realization) such that 
$$\iota^{-1}(\mbf{Wall})=V_\Gamma^{anch}.$$
\end{itemize}
Often we will abuse notation and denote the graph diagram simply by $\Gamma$.

A morphism between graph diagrams $\iota:\Gamma\to\Sigma$ and $\iota':\Gamma'\to\Sigma$ is a map between the geometric realizations of the graphs, $\mu:\lvert\Gamma\rvert\to\lvert\Gamma'\rvert$, such that $\iota=\iota'\circ\mu$ and $\mu(V_\Gamma)\subseteq V_{\Gamma'}$.

A homotopy between graph diagrams $\iota_0:\Gamma\to\Sigma$ and $\iota_1:\Gamma\to\Sigma$ is  one parameter family of graph diagrams $\iota_t:\Gamma\to\Sigma$, $t\in[0,1]$.
\end{defn}

Suppose $\iota:\Gamma\to \Sigma$ is a graph diagram. We denote by $\tilde\iota$ the natural map 
$V_\Gamma^{anch}\to\mbf{Wall}$. 
Let $\Gamma_d=\iota^{-1}(\Sigma_d)$,\footnote{I.e. the subgraph of $\Gamma$ consisting of edges mapped into $\Sigma_d$} then $\Gamma_d\to \Sigma_d$ is a graph diagram in the marked surface $(\Sigma_d,V_d)$, (here we have composed with the blow down map $\hat\Sigma_d\to \Sigma_d$). 
Define
$$\mc{A}_{(\iota:\Gamma\to\Sigma)}:=C^\infty\bigg(\prod_{d\in\mbf{Dom}}G_d^{E_{\Gamma_d}}\bigg)^{C_\Gamma^{int}\times C_\Gamma^{anch}},$$
where $$C_\Gamma^{int}:=\prod_{d\in\mbf{Dom}}G_d^{V_{\Gamma_d}^{int}},\quad C_\Gamma^{anch}:=\prod_{v\in V_{\Gamma}^{anch}}C_{\tilde\iota(v)},$$ 
with 
$V_{\Gamma_d}^{int}:=V_{\Gamma_d}\smallsetminus V_{\Gamma_d}^{anch}$.
Notice that
\begin{equation}\label{eq:FuncSpcSemiTensDesc}\mc{A}_{(\iota:\Gamma\to\Sigma)}=\big(\hat\otimes_{d\in\mbf{Dom}}\mc{A}_{\Gamma_d}(G_d)\big)^{C_\Gamma^{anch}}.\end{equation}

To simplify notation, we will often abbreviate $\mc{A}_{(\iota:\Gamma\to\Sigma)}$ to $\mc{A}_\Gamma$.

As in the case of spin networks in marked surfaces, we have the following lemma.
\begin{lem}

\begin{itemize}
\item The algebras $\mc{A}_{(\iota_0:\Gamma\to\Sigma)}$ and $\mc{A}_{(\iota_1:\Gamma\to\Sigma)}$ are canonically identified for homotopic graph diagrams $\iota_0:\Gamma\to\Sigma$ and $\iota_1:\Gamma\to\Sigma$.
\item The assignment $(\iota:\Gamma\to\Sigma)\mapsto \mc{A}_{(\iota:\Gamma\to\Sigma)}$ is a functor from graph diagrams to algebras.
\end{itemize}
\end{lem}

\begin{defn}
A \emph{spin network in $\Sigma$} is a pair $(\Gamma,f)$, where $\Gamma\to \Sigma$ is a graph diagram, and $f\in \mc{A}_\Gamma$.

We say that spin networks $(\Gamma_0,f_0)$ and $(\Gamma_1,f_1)$ are \emph{homotopic} if the underlying graph diagrams are homotopic and $f_0=f_1\in \mc{A}_{\Gamma_0}\equiv\mc{A}_{\Gamma_1}$. 

A morphism of spin networks $\mu:(\Gamma,f)\to (\Gamma',f')$ is a morphism of graph diagrams $\mu:\lvert\Gamma\rvert\to\lvert\Gamma'\rvert$ such that $f'=\mu_*f$. 
\end{defn}

Let $$\mc{M}_\Sigma:=\big(\prod_{d\in\mbf{Dom}}M_{\Sigma_d,V_d}(G_d)\big)/C$$ denote the moduli space for the quilted surface, $\Sigma$. Then 
$$C^\infty(\mc{M}_\Sigma)=C^\infty\big(\prod_{d\in\mbf{Dom}}M_{\Sigma_d,V_d}(G_d)\big)^C=\big(\hat\otimes_{d\in\mbf{Dom}}C^\infty\big(M_{\Sigma_d,V_d}(G_d)\big)\big)^C$$
Since the map
$$\on{ev}:\hat\otimes_{d\in\mbf{Dom}}\mc{A}_{\Gamma_d}(G_d)\to \hat\otimes_{d\in\mbf{Dom}}C^\infty\big(M_{\Sigma_d,V_d}(G_d)\big)$$
is $\prod_{d\in\mbf{Dom}}G_d^{V_d}$-equivariant it restricts to a map of algebras
$$\on{ev}:\mc{A}_\Gamma\to C^\infty(\mc{M}_\Sigma).$$

As before, we consider two spin networks to be equivalent if they are related by a chain of homotopies and morphisms, and we define $\on{SpinNet}_\Sigma$ to be the set of equivalence classes of spin networks in the quilted surface, $\Sigma$. The same arguments as in \cref{sec:SpNetMrSFunct} shows that $\on{ev}$ descends to a map of equivalence classes:
$$\on{ev}:\on{SpinNet}_\Sigma\to C^\infty(\mc{M}_\Sigma).$$
Moreover, the analogues of \cref{lem:AlgofSpinNet}, \cref{prop:EvAlgId}, and \cref{prop:qPBrkSpNet} hold for $\on{SpinNet}_\Sigma$. 

To be precise, suppose that $\iota:\Gamma\to\Sigma$ and $\iota':\Gamma'\to\Sigma$ are two graph diagrams which are transverse. That is, $\iota(V_\Gamma)\cap\iota'(V_{\Gamma'})=\emptyset$ and the restrictions of $\iota$ and $\iota'$ to the edges are transverse. Suppose $A\in \Gamma\times_\Sigma\Gamma'$ and let $d_A\in\mbf{Dom}$ denote the domain such that $A\in \Sigma_{d_A}$. We define the graph diagram $\iota\cup_A\iota':\Gamma\cup_{A}\Gamma'\to \Sigma$ as in \cref{sec:QPoisBrakSpin}, and we define
$$\Psi:\prod_{d\in\mbf{Dom}}G_{d}^{E_{(\Gamma\cup_{\mbf{p}}\Gamma')_{d}}}\to \otimes^2 T\big(\prod_{d\in\mbf{Dom}}G_{d}^{E_{(\Gamma\cup\Gamma')_{d}}}\big)$$
to be the sum of $\Psi_{s_{d_A}}$ on the $d_A$-th factor with the zero sections on the other factors.

\begin{thm}\label{thm:PoisAlgSpinNet}
The map 
$$\on{ev}:\on{SpinNet}_\Sigma\to C^\infty(\mc{M}_\Sigma)$$
is an isomorphism of Poisson algebras, where scalar multiplication, addition, and multiplication are defined on $\on{SpinNet}_\Sigma$ by \cref{eq:SpinNetOp}, and the Poisson bracket is defined by
\begin{equation}\label{eq:PoissSpinNetBrack}\{[\Gamma,f],[\Gamma',f']\}:=\sum_{A\in \Gamma\times_\Sigma\Gamma'}\on{sign}(A)[\Gamma\cup_{A}\Gamma',\Psi^*(\d f\otimes \d f')],\end{equation}
where $\Gamma$ and $\Gamma'$ are assumed to be transverse graph diagrams.
\end{thm}
\begin{proof}
The proof that $\on{ev}:\on{SpinNet}_\Sigma\to C^\infty(\mc{M}_\Sigma)$ is an isomorphism of algebras is entirely analogous to that of \cref{prop:EvAlgId}, and so we omit it.

As explained in \cref{sec:QPoisBrakSpin}, $\Psi$ is equivariant with respect to the action of 
$\prod_{d\in\mbf{Dom}}G_d^{V_{(\Gamma\cup_A\Gamma')_d}}.$
 In particular, $\Psi^*(\d f\otimes \d f')\in \mc{A}_{(\Gamma\cup_A\Gamma')}$.
 
\Cref{eq:FuncSpcSemiTensDesc} identifies $\on{SpinNet}_\Sigma$ with the subalgebra of 
$C$-invariant elements of $$\hat\otimes_{d\in\mbf{Dom}}\on{SpinNet}_{\Sigma_d,V_d}(G_d).$$
Moreover, the following diagram commutes:
$$\begin{tikzpicture}
\mmat{m}{\on{SpinNet}_\Sigma & C^\infty(\mc{M}_\Sigma)\\
\hat\otimes_{d\in\mbf{Dom}}\on{SpinNet}_{\Sigma_d,V_d}(G_d) & C^\infty\big(\prod_{d\in\mbf{Dom}}M_{\Sigma_d,V_d}(G_d)\big)\\};
\draw[->] (m-1-1) edge node{$\on{ev}$} (m-1-2);
\draw[->] (m-1-1) edge (m-2-1);
\draw[->] (m-1-2) edge (m-2-2);
\draw[->] (m-2-1) edge node{$\on{ev}$} (m-2-2);
\end{tikzpicture}$$
Hence \cref{thm:Red,prop:qPBrkSpNet} imply that \cref{eq:PBrkSpNetMrkPr}
 defines a Poisson bracket on $\on{SpinNet}_\Sigma$ for which $\on{ev}$ is a morphism of Poisson algebras. Explicitly, this bracket is
\begin{multline*}\big\{[\Gamma,f],[\Gamma',f']\big\}:=\sum_{d\in\mbf{Dom}}\bigg(\sum_{A\in(\Gamma\times_\Sigma\Gamma')_d^{int}}\on{sign}(A)\big[\Gamma\cup_{A}\Gamma',\Psi^*(\d f\otimes \d f')\big]\\
+\sum_{A=(p,p')\in (\Gamma_d\times_{\Sigma_d}\Gamma'_d)^{anch}}\frac{1}{2}\on{sign}(A)\big[\Gamma\cup\Gamma',(\d f\otimes \d f')\big((\rho_p\otimes\rho_{p'})s_d\big)\big]\bigg).
\end{multline*}
The first term simplifies to yield \cref{eq:PoissSpinNetBrack}, so we need only show that the second term vanishes.

\begin{figure}[h]
\begin{center}
\begin{subfigure}[t]{.49\linewidth}
\centering
\def\svgwidth{\linewidth}
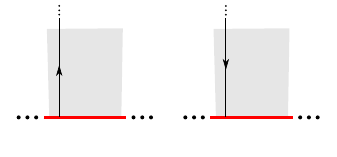
\caption{$\on{Sign}(p,d)=\pm1$ depends on the direction of the edge in $\Gamma$ leaving $p$.}
\end{subfigure}
\begin{subfigure}[t]{.49\linewidth}
\centering
\def\svgwidth{\linewidth}
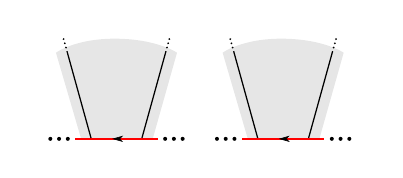
\caption{
$\on{Sign}(p,p',\mbf{w},d)=\pm1$ depends on the order of $p$ and $p'$, where we give $\mbf{w}\subset\partial\hat\Sigma_d$ the boundary orientation.}
\end{subfigure}
\begin{subfigure}[t]{.49\linewidth}
\centering
\def\svgwidth{\linewidth}
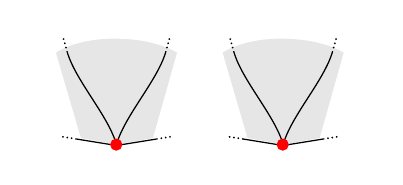
\caption{
The blow down map $\hat\Sigma_d\to\Sigma_d$ contracts $\mbf{w}$ to a point. Thus the blow down equates $\epsilon(p,p',\mbf{w},d)$ with $\on{sign}(A)$, where the latter is computed as in \cref{fig:Signp}}
\end{subfigure}
\end{center}
\caption{
\label{fig:Signs}}
\end{figure}

In terms of graph diagrams in the quilted surface, $\Sigma$, the second term can be rewritten as
$$\frac{1}{2}\sum_{\mbf{w}\in\mbf{Wall}}\;\sum_{p\in\mbf{w}\cap\Gamma,p'\in\mbf{w}\cap\Gamma'}\;\sum_{d\in\mbf{Dom}(\mbf{w})}\epsilon(p,p',\mbf{w},d)\big[\Gamma\cup\Gamma',(\d f\otimes \d f')\big((\rho_p\otimes\rho_{p'})s_d\big)\big]$$
where 
$$\mbf{Dom}(\mbf{w})=\{d\in\mbf{Dom};\  \mbf{w}\subset\partial\hat\Sigma_d\}.$$
Both 
$$\epsilon(p,p',\mbf{w},d)=\on{sign}(p,d)\on{sign}(p',d)\on{sign}(p,p',\mbf{w},d),$$ and $\on{sign}(p,d)=\pm1$ and $\on{sign}(p,p',\mbf{w},d)=\pm1$
are computed as in \cref{fig:Signs}. Note that $\epsilon(p,p',\mbf{w},d)=\on{sign}(A)$, where the left hand side is computed in the blow up, $\hat\Sigma_d$, and the right hand side is computed for $A=(p,p')$ in the blow down, $\Sigma_d$.

Now both $f$ and $f'$ are $C_{\mbf{w}}$-invariant, that is 
$$\bigoplus_{d\in\mbf{Dom}(\mbf{w})}\rho_p^*\d f\in\mf{c}_{\mbf{w}}^\perp,\quad
\bigoplus_{d\in\mbf{Dom}(\mbf{w})}\rho_{p'}^*\d f'\in\mf{c}_{\mbf{w}}^\perp$$
Since $C_\mbf{w}$ was chosen to be coisotropic, $\sum_{d\in\mbf{Dom}(\mbf{w})}\epsilon(p,p',\mbf{w},d)s_d$ vanishes on $\mf{c}_{\mbf{w}}^\perp$. This completes the proof.
\end{proof}

Suppose $\Sigma$ and $\Sigma'$ are two quilted surfaces. A map $\Sigma\to\Sigma'$ is called an \emph{embedding of quilted surfaces} if 
\begin{itemize}
\item it is an embedding of surfaces $\Sigma\to \Sigma'$,
\item it maps each domain, $\Sigma_d$, of $\Sigma$ into a domain, $\Sigma'_{d'}$, of $\Sigma'$ which is colored with the same structure group (i.e. $G_d\equiv G_{d'}$), and
\item it maps each domain wall, $\mbf{w}$, of $\Sigma$ into a domain wall, $\mbf{w}'$, of $\Sigma'$ which is colored with the same coisotropic relation (i.e. $C_{\mbf{w}} \equiv C_{\mbf{w}'}$).
\end{itemize}

\begin{cor}\label{cor:embedding}
An embedding of quilted surfaces $\Sigma\to\Sigma'$ induces a Poisson morphism 
$$\mc{M}_{\Sigma'}\to\mc{M}_\Sigma$$
of the corresponding moduli spaces.
\end{cor}
\begin{proof}
This follows directly from the functoriality of the assignment  $\Sigma\to \on{SpinNet}_{\Sigma}$ of the Poisson algebra of spin networks to quilted surfaces and \cref{thm:PoisAlgSpinNet}.

Alternatively, it follows from \cref{cor:emb} and the definition of the Poisson structure on $\mc{M}_\Sigma$ given in \cref{eq:ModSpcQuilSurf}.
\end{proof}

\section{Colorful examples}

\begin{example}[Poisson Lie groups]\label{ex:PoissLie2}
Suppose that $(\g,\mf{a},\mf{b})$ is a Manin triple and consider the quilted surface $\Sigma$ pictured below
\begin{center}
\begingroup%
  \makeatletter%
  \providecommand\color[2][]{%
    \errmessage{(Inkscape) Color is used for the text in Inkscape, but the package 'color.sty' is not loaded}%
    \renewcommand\color[2][]{}%
  }%
  \providecommand\transparent[1]{%
    \errmessage{(Inkscape) Transparency is used (non-zero) for the text in Inkscape, but the package 'transparent.sty' is not loaded}%
    \renewcommand\transparent[1]{}%
  }%
  \providecommand\rotatebox[2]{#2}%
  \ifx\svgwidth\undefined%
    \setlength{\unitlength}{107.45510254bp}%
    \ifx\svgscale\undefined%
      \relax%
    \else%
      \setlength{\unitlength}{\unitlength * \real{\svgscale}}%
    \fi%
  \else%
    \setlength{\unitlength}{\svgwidth}%
  \fi%
  \global\let\svgwidth\undefined%
  \global\let\svgscale\undefined%
  \makeatother%
  \begin{picture}(1,0.89646526)%
    \put(0,0){\includegraphics[width=\unitlength]{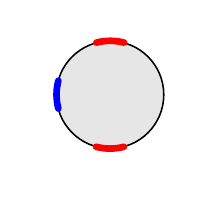}}%
    \put(0.48386092,0.0847116){\color[rgb]{0,0,0}\makebox(0,0)[b]{\smash{$B$}}}%
    \put(0.20529777,0.45763216){\color[rgb]{0,0,0}\makebox(0,0)[rb]{\smash{$A$}}}%
    \put(0.31048103,0.25491191){\color[rgb]{0,0,0}\makebox(0,0)[rb]{\smash{$\mbf{1}$}}}%
    \put(0.48386092,0.75475889){\color[rgb]{0,0,0}\makebox(0,0)[b]{\smash{$B$}}}%
    \put(0.31048103,0.65694028){\color[rgb]{0,0,0}\makebox(0,0)[rb]{\smash{$\mbf{1}$}}}%
  \end{picture}%
\endgroup%

\end{center}
where $A,B\subset G$ are Lie groups integrating $\mf{a}$ and $\mf{b}$ such that $A\cap B=\mathbf{1}$.
Constraining the holonomies along the uncolored boundary arcs marked by $\mbf{1}$ to be trivial is effectively the same as contracting those arcs to points, which results in the first image below.
\begin{center}
\begingroup%
  \makeatletter%
  \providecommand\color[2][]{%
    \errmessage{(Inkscape) Color is used for the text in Inkscape, but the package 'color.sty' is not loaded}%
    \renewcommand\color[2][]{}%
  }%
  \providecommand\transparent[1]{%
    \errmessage{(Inkscape) Transparency is used (non-zero) for the text in Inkscape, but the package 'transparent.sty' is not loaded}%
    \renewcommand\transparent[1]{}%
  }%
  \providecommand\rotatebox[2]{#2}%
  \ifx\svgwidth\undefined%
    \setlength{\unitlength}{227.45605469bp}%
    \ifx\svgscale\undefined%
      \relax%
    \else%
      \setlength{\unitlength}{\unitlength * \real{\svgscale}}%
    \fi%
  \else%
    \setlength{\unitlength}{\svgwidth}%
  \fi%
  \global\let\svgwidth\undefined%
  \global\let\svgscale\undefined%
  \makeatother%
  \begin{picture}(1,0.36446893)%
    \put(0,0){\includegraphics[width=\unitlength]{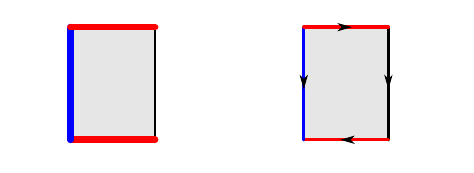}}%
    \put(0.7327972,0.33237344){\color[rgb]{0,0,0}\makebox(0,0)[b]{\smash{$b_1$}}}%
    \put(0.7327972,0.03298525){\color[rgb]{0,0,0}\makebox(0,0)[b]{\smash{$b_2$}}}%
    \put(0.83971593,0.19257909){\color[rgb]{0,0,0}\makebox(0,0)[lb]{\smash{$g$}}}%
    \put(0.61752715,0.19257909){\color[rgb]{0,0,0}\makebox(0,0)[rb]{\smash{$a$}}}%
    \put(0.24039439,0.02595093){\color[rgb]{0,0,0}\makebox(0,0)[b]{\smash{$B$}}}%
    \put(0.34731312,0.19257909){\color[rgb]{0,0,0}\makebox(0,0)[lb]{\smash{$G$}}}%
    \put(0.12512434,0.19257909){\color[rgb]{0,0,0}\makebox(0,0)[rb]{\smash{$A$}}}%
    \put(0.24039439,0.33546126){\color[rgb]{0,0,0}\makebox(0,0)[b]{\smash{$B$}}}%
  \end{picture}%
\endgroup%

\end{center}
Following \cref{cor:PartTriv}, the moduli space is
$$\{(b_1,g,b_2,a)\in B\times G\times B \times A;\  b_2gb_1=a\}/(\hat H\cap \hat C),$$
where the holonomies are as pictured above.
Since $A\cap B=\mbf{1}$, the group $(\hat H\cap \hat C)\cong B\times B$ of residual gauge transformations  acts non-trivially only at the right two vertices by 
$$(b,b')\cdot(b_1,g,b_2,a)=(bb_1,b'gb^{-1},b_2b'^{-1},a).$$
Thus the moduli space is identified with $A$. The Poisson structure on $A$ is the Poisson Lie structure. This example should be compared with \cref{ex:PoisLie1}.

Consider the embedding of quilted surfaces pictured below:
\begin{center}
\begingroup%
  \makeatletter%
  \providecommand\color[2][]{%
    \errmessage{(Inkscape) Color is used for the text in Inkscape, but the package 'color.sty' is not loaded}%
    \renewcommand\color[2][]{}%
  }%
  \providecommand\transparent[1]{%
    \errmessage{(Inkscape) Transparency is used (non-zero) for the text in Inkscape, but the package 'transparent.sty' is not loaded}%
    \renewcommand\transparent[1]{}%
  }%
  \providecommand\rotatebox[2]{#2}%
  \ifx\svgwidth\undefined%
    \setlength{\unitlength}{147.40119629bp}%
    \ifx\svgscale\undefined%
      \relax%
    \else%
      \setlength{\unitlength}{\unitlength * \real{\svgscale}}%
    \fi%
  \else%
    \setlength{\unitlength}{\svgwidth}%
  \fi%
  \global\let\svgwidth\undefined%
  \global\let\svgscale\undefined%
  \makeatother%
  \begin{picture}(1,0.34871632)%
    \put(0,0){\includegraphics[width=\unitlength]{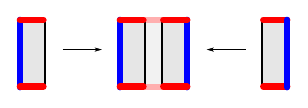}}%
    \put(0.26686553,0.20480329){\color[rgb]{0,0,0}\makebox(0,0)[b]{\smash{$t$}}}%
    \put(0.74133706,0.20480329){\color[rgb]{0,0,0}\makebox(0,0)[b]{\smash{$s$}}}%
  \end{picture}%
\endgroup%

\end{center}
As explained above, the moduli space corresponding to the quilted surfaces depicted on the left and right is the Poisson Lie group, $A$. Meanwhile, as explained in \cref{ex:SymplDbl}, the moduli space for the quilted surface depicted in the middle is the symplectic groupoid
$$\{(a_1,b_1,a_2,b_2)\in A\times B\times A\times B;\ a_1b_1=b_2a_2\},$$
 integrating the Poisson Lie structure on $A$. The embedding of quilted surfaces induces the map
 $$t\times s:(a_1,b_1,a_2,b_2)\to (a_2,a_1).$$
 The first component of the map is the Lie groupoid target map, while the second component is the Lie groupoid source map. \cref{cor:embedding} shows that this map is a Poisson morphism.

The multiplication for the Poisson Lie group $A\cong\mc{M}_\Sigma$ can be understood by considering the following sequence of embeddings:
\begin{center}
\def\svgwidth{\linewidth}
\begingroup%
  \makeatletter%
  \providecommand\color[2][]{%
    \errmessage{(Inkscape) Color is used for the text in Inkscape, but the package 'color.sty' is not loaded}%
    \renewcommand\color[2][]{}%
  }%
  \providecommand\transparent[1]{%
    \errmessage{(Inkscape) Transparency is used (non-zero) for the text in Inkscape, but the package 'transparent.sty' is not loaded}%
    \renewcommand\transparent[1]{}%
  }%
  \providecommand\rotatebox[2]{#2}%
  \ifx\svgwidth\undefined%
    \setlength{\unitlength}{490.76020508bp}%
    \ifx\svgscale\undefined%
      \relax%
    \else%
      \setlength{\unitlength}{\unitlength * \real{\svgscale}}%
    \fi%
  \else%
    \setlength{\unitlength}{\svgwidth}%
  \fi%
  \global\let\svgwidth\undefined%
  \global\let\svgscale\undefined%
  \makeatother%
  \begin{picture}(1,0.15455567)%
    \put(0,0){\includegraphics[width=\unitlength]{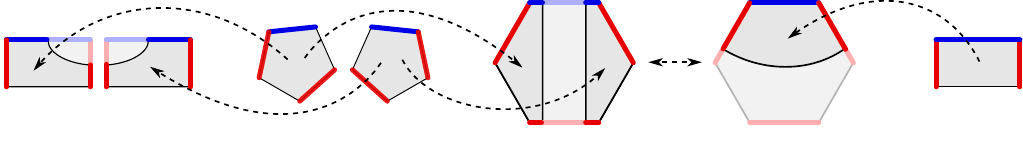}}%
    \put(0.66122378,0.09977918){\color[rgb]{0,0,0}\makebox(0,0)[b]{\smash{$=$}}}%
    \put(0.55180737,0.00336691){\color[rgb]{0,0,0}\makebox(0,0)[b]{\smash{$\tilde\Sigma'$}}}%
    \put(0.76698374,0.00336691){\color[rgb]{0,0,0}\makebox(0,0)[b]{\smash{$\tilde\Sigma'$}}}%
    \put(0.95647342,0.00336691){\color[rgb]{0,0,0}\makebox(0,0)[b]{\smash{$\Sigma$}}}%
    \put(0.33595768,0.00336691){\color[rgb]{0,0,0}\makebox(0,0)[b]{\smash{$\tilde\Sigma\sqcup\tilde\Sigma$}}}%
    \put(0.09606915,0.00336691){\color[rgb]{0,0,0}\makebox(0,0)[b]{\smash{$\Sigma\sqcup\Sigma$}}}%
  \end{picture}%
\endgroup%

\end{center}
which induces the following sequence of Poisson morphisms:
$$\mc{M}_{\Sigma\sqcup\Sigma}\xrightarrow{\iota} \mc{M}_{\tilde\Sigma\sqcup\tilde\Sigma}\xleftarrow{\varphi} \mc{M}_{\tilde\Sigma'}\xrightarrow{=} \mc{M}_{\tilde\Sigma'}\xrightarrow{m}\mc{M}_\Sigma$$ 
The map $\varphi$ is an embedding which contains the image of $\iota$. Composing this sequence, $\mu:=m\circ \varphi^{-1}\circ \iota$, yields a multiplication map 
$$\mu: \mc{M}_{\Sigma}\times\mc{M}_{\Sigma}=\mc{M}_{\Sigma\sqcup\Sigma}\to \mc{M}_\Sigma.$$
Associativity can be verified via another sequence of embeddings of quilted surfaces (we omit the details).
%
%
Under the canonical identification $\mc{M}_{\Sigma}\cong A$ with the Lie group $A$, the map $\mu$ becomes the usual multiplication. However, as a consequence of this construction, we see that multiplication is a Poisson map.

\end{example}

\begin{example}[Double Poisson Lie group]\label{ex:DblPoissLie}
Suppose once again that $(\g,\mf{a},\mf{b})$ is a Manin triple, and that $A,B\subset G$ are Lie groups integrating $\mf{a}$ and $\mf{b}$ such that $A\cap B=\mathbf1$. Then, $G$ is a Poisson Lie group (cf. \cite{Drinfeld83}, see also \cite{thesis-3}, for example), we may describe this Poisson structure in terms of a moduli space as follows: Let $\Sigma$ be the quilted surface pictured below left,
\begin{center}
\begingroup%
  \makeatletter%
  \providecommand\color[2][]{%
    \errmessage{(Inkscape) Color is used for the text in Inkscape, but the package 'color.sty' is not loaded}%
    \renewcommand\color[2][]{}%
  }%
  \providecommand\transparent[1]{%
    \errmessage{(Inkscape) Transparency is used (non-zero) for the text in Inkscape, but the package 'transparent.sty' is not loaded}%
    \renewcommand\transparent[1]{}%
  }%
  \providecommand\rotatebox[2]{#2}%
  \ifx\svgwidth\undefined%
    \setlength{\unitlength}{283.33515625bp}%
    \ifx\svgscale\undefined%
      \relax%
    \else%
      \setlength{\unitlength}{\unitlength * \real{\svgscale}}%
    \fi%
  \else%
    \setlength{\unitlength}{\svgwidth}%
  \fi%
  \global\let\svgwidth\undefined%
  \global\let\svgscale\undefined%
  \makeatother%
  \begin{picture}(1,0.23904345)%
    \put(0,0){\includegraphics[width=\unitlength]{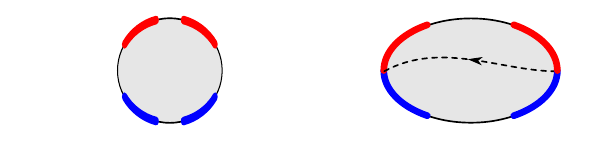}}%
    \put(0.64152709,0.17988638){\color[rgb]{0,0,0}\makebox(0,0)[b]{\smash{$B$}}}%
    \put(0.64152713,0.04716411){\color[rgb]{0,0,0}\makebox(0,0)[b]{\smash{$A$}}}%
    \put(0.79024001,0.1503037){\color[rgb]{0,0,0}\makebox(0,0)[b]{\smash{$g$}}}%
    \put(0.94646634,0.17988638){\color[rgb]{0,0,0}\makebox(0,0)[b]{\smash{$B$}}}%
    \put(0.94646634,0.04716411){\color[rgb]{0,0,0}\makebox(0,0)[b]{\smash{$A$}}}%
    \put(0.38962798,0.11587881){\color[rgb]{0,0,0}\makebox(0,0)[lb]{\smash{$\mbf{1}$}}}%
    \put(0.18068818,0.1158788){\color[rgb]{0,0,0}\makebox(0,0)[rb]{\smash{$\mbf{1}$}}}%
    \put(0.1841183,0.17988638){\color[rgb]{0,0,0}\makebox(0,0)[b]{\smash{$B$}}}%
    \put(0.1841183,0.04716411){\color[rgb]{0,0,0}\makebox(0,0)[b]{\smash{$A$}}}%
    \put(0.39305814,0.17988638){\color[rgb]{0,0,0}\makebox(0,0)[b]{\smash{$B$}}}%
    \put(0.39305814,0.04716411){\color[rgb]{0,0,0}\makebox(0,0)[b]{\smash{$A$}}}%
  \end{picture}%
\endgroup%

\end{center}
Contracting the uncolored boundary arcs marked by $\mbf{1}$ results in the image on the right. Since $A\cap B=\mbf{1}$, the group $(\hat H\cap \hat C)$ of residual gauge transformations is trivial (cf. \cref{cor:PartTriv}).
Consequently, computing the holonomy along the dotted arc pictured in the second image, we may identify the moduli space with $G$. The resulting Poisson structure on $G$ is Drinfel'd's Double Poisson Lie structure.

Applying \cref{cor:embedding} to the following embedding of quilted surfaces 
\begin{center}
\begingroup%
  \makeatletter%
  \providecommand\color[2][]{%
    \errmessage{(Inkscape) Color is used for the text in Inkscape, but the package 'color.sty' is not loaded}%
    \renewcommand\color[2][]{}%
  }%
  \providecommand\transparent[1]{%
    \errmessage{(Inkscape) Transparency is used (non-zero) for the text in Inkscape, but the package 'transparent.sty' is not loaded}%
    \renewcommand\transparent[1]{}%
  }%
  \providecommand\rotatebox[2]{#2}%
  \ifx\svgwidth\undefined%
    \setlength{\unitlength}{187.82496338bp}%
    \ifx\svgscale\undefined%
      \relax%
    \else%
      \setlength{\unitlength}{\unitlength * \real{\svgscale}}%
    \fi%
  \else%
    \setlength{\unitlength}{\svgwidth}%
  \fi%
  \global\let\svgwidth\undefined%
  \global\let\svgscale\undefined%
  \makeatother%
  \begin{picture}(1,0.39001237)%
    \put(0,0){\includegraphics[width=\unitlength]{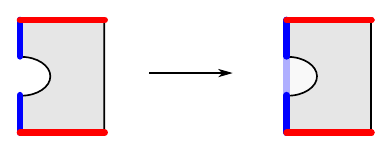}}%
  \end{picture}%
\endgroup%

\end{center}
yields the morphism of Poisson Lie groups $A\to G$.
\end{example}

\begin{example}[Symplectic double groupoid integrating Drinfel'd's double \cite{Severa98,Severa:2011ug}]\label{ex:symplDblDbl}
Suppose once again that $(\g,\mf{a},\mf{b})$ is a Manin triple, and that $A,B\subset G$ are Lie groups integrating $\mf{a}$ and $\mf{b}$ such that the product map $A\times B\to  G$ is a diffeomorphism. We may identify the symplectic double groupoid integrating the Poisson Lie structure on $G$ (described in the previous example) with a moduli space, as follows: Let $\Sigma$ be the quilted surface pictured below, where we have already contracted all the uncolored boundary arcs.
\begin{center}
\begingroup%
  \makeatletter%
  \providecommand\color[2][]{%
    \errmessage{(Inkscape) Color is used for the text in Inkscape, but the package 'color.sty' is not loaded}%
    \renewcommand\color[2][]{}%
  }%
  \providecommand\transparent[1]{%
    \errmessage{(Inkscape) Transparency is used (non-zero) for the text in Inkscape, but the package 'transparent.sty' is not loaded}%
    \renewcommand\transparent[1]{}%
  }%
  \providecommand\rotatebox[2]{#2}%
  \ifx\svgwidth\undefined%
    \setlength{\unitlength}{287.77614746bp}%
    \ifx\svgscale\undefined%
      \relax%
    \else%
      \setlength{\unitlength}{\unitlength * \real{\svgscale}}%
    \fi%
  \else%
    \setlength{\unitlength}{\svgwidth}%
  \fi%
  \global\let\svgwidth\undefined%
  \global\let\svgscale\undefined%
  \makeatother%
  \begin{picture}(1,0.38586472)%
    \put(0,0){\includegraphics[width=\unitlength]{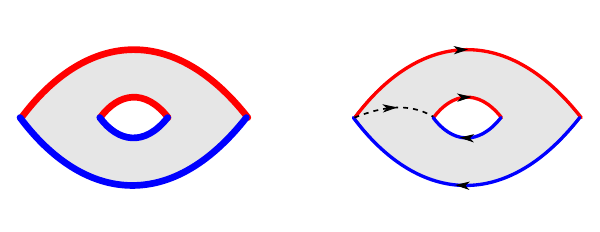}}%
    \put(0.22340294,0.33694432){\color[rgb]{0,0,0}\makebox(0,0)[b]{\smash{$B$}}}%
    \put(0.22340294,0.03188513){\color[rgb]{0,0,0}\makebox(0,0)[b]{\smash{$A$}}}%
    \put(0.22340294,0.23686653){\color[rgb]{0,0,0}\makebox(0,0)[b]{\smash{$B$}}}%
    \put(0.22340294,0.12084317){\color[rgb]{0,0,0}\makebox(0,0)[b]{\smash{$A$}}}%
    \put(0.77939064,0.33694432){\color[rgb]{0,0,0}\makebox(0,0)[b]{\smash{$b_1$}}}%
    \put(0.77939064,0.03188513){\color[rgb]{0,0,0}\makebox(0,0)[b]{\smash{$a_1$}}}%
    \put(0.77939064,0.23686653){\color[rgb]{0,0,0}\makebox(0,0)[b]{\smash{$b_2$}}}%
    \put(0.77939064,0.12084317){\color[rgb]{0,0,0}\makebox(0,0)[b]{\smash{$a_2$}}}%
    \put(0.66517833,0.17790569){\color[rgb]{0,0,0}\makebox(0,0)[b]{\smash{$g$}}}%
  \end{picture}%
\endgroup%

\end{center}
The holonomies along the arcs pictured above identify the moduli space with
$$\{(a_1,b_1,a_2,b_2,g)\in A\times B\times A\times B\times G;\  ga_1b_1=a_2b_2g\},$$
where once again, the group of residual gauge transformations is trivial.

Applying \cref{cor:embedding} to the following embedding of quilted surfaces 
\begin{center}
\def\svgwidth{\linewidth}
\begingroup%
  \makeatletter%
  \providecommand\color[2][]{%
    \errmessage{(Inkscape) Color is used for the text in Inkscape, but the package 'color.sty' is not loaded}%
    \renewcommand\color[2][]{}%
  }%
  \providecommand\transparent[1]{%
    \errmessage{(Inkscape) Transparency is used (non-zero) for the text in Inkscape, but the package 'transparent.sty' is not loaded}%
    \renewcommand\transparent[1]{}%
  }%
  \providecommand\rotatebox[2]{#2}%
  \ifx\svgwidth\undefined%
    \setlength{\unitlength}{352.57614746bp}%
    \ifx\svgscale\undefined%
      \relax%
    \else%
      \setlength{\unitlength}{\unitlength * \real{\svgscale}}%
    \fi%
  \else%
    \setlength{\unitlength}{\svgwidth}%
  \fi%
  \global\let\svgwidth\undefined%
  \global\let\svgscale\undefined%
  \makeatother%
  \begin{picture}(1,0.23935012)%
    \put(0,0){\includegraphics[width=\unitlength]{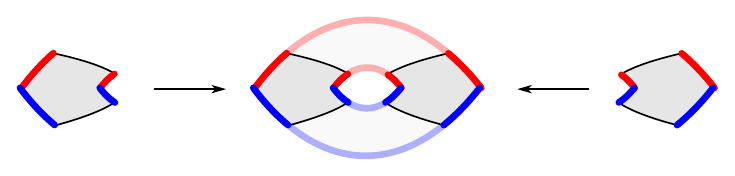}}%
    \put(0.25434485,0.13158511){\color[rgb]{0,0,0}\makebox(0,0)[b]{\smash{$t$}}}%
    \put(0.75352779,0.13158511){\color[rgb]{0,0,0}\makebox(0,0)[b]{\smash{$s$}}}%
  \end{picture}%
\endgroup%

\end{center}
yields the Poisson morphism 
$$t\times s:(a_1,b_1,a_2,b_2,g)\to (g,b_2gb_1^{-1})$$
whose components are the target/source map (respectively) onto $G$ endowed with the Poisson Lie structure described in \cref{ex:DblPoissLie}.

This moduli space was first studied in \cite{Severa98,Severa:2011ug}.

\end{example}

\begin{example}[Lu-Yakimov Poisson homogeneous spaces]\label{ex:luYak}
Suppose once again that $(\g,\mf{a},\mf{b})$ is a Manin triple, and that $A,B\subset G$ are Lie groups integrating $\mf{a}$ and $\mf{b}$ such that $A\cap B=\mbf{1}$. Let $C\subseteq G$ be a closed subgroup whose Lie subalgebra $\mf{c}\subseteq \g$ is coisotropic. Lu and Yakimov \cite{Lu06} describe a Poisson structure on $G/C$, which we may identify with the moduli space for the following quilted surface: 
\begin{center}
\begingroup%
  \makeatletter%
  \providecommand\color[2][]{%
    \errmessage{(Inkscape) Color is used for the text in Inkscape, but the package 'color.sty' is not loaded}%
    \renewcommand\color[2][]{}%
  }%
  \providecommand\transparent[1]{%
    \errmessage{(Inkscape) Transparency is used (non-zero) for the text in Inkscape, but the package 'transparent.sty' is not loaded}%
    \renewcommand\transparent[1]{}%
  }%
  \providecommand\rotatebox[2]{#2}%
  \ifx\svgwidth\undefined%
    \setlength{\unitlength}{128.29272461bp}%
    \ifx\svgscale\undefined%
      \relax%
    \else%
      \setlength{\unitlength}{\unitlength * \real{\svgscale}}%
    \fi%
  \else%
    \setlength{\unitlength}{\svgwidth}%
  \fi%
  \global\let\svgwidth\undefined%
  \global\let\svgscale\undefined%
  \makeatother%
  \begin{picture}(1,0.64531439)%
    \put(0,0){\includegraphics[width=\unitlength]{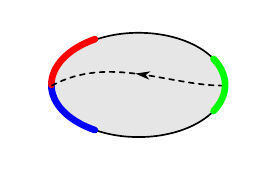}}%
    \put(0.17304547,0.460866){\color[rgb]{0,0,0}\makebox(0,0)[b]{\smash{$B$}}}%
    \put(0.17304554,0.16774817){\color[rgb]{0,0,0}\makebox(0,0)[b]{\smash{$A$}}}%
    \put(0.91312551,0.30194668){\color[rgb]{0,0,0}\makebox(0,0)[b]{\smash{$C$}}}%
    \put(0.50147872,0.3955325){\color[rgb]{0,0,0}\makebox(0,0)[b]{\smash{$g$}}}%
  \end{picture}%
\endgroup%

\end{center}
Computing the holonomy along the dotted arc yields an element $g\in G$, but the group $C$ of residual gauge transformations acts by right multiplication on this element.
Thus, following \cref{cor:PartTriv}, the moduli space $G/C$ carries a Poisson structure.

The symplectic groupoid integrating the Poisson structure on $G/C$ is the moduli space corresponding to the quilted surface pictured below:
\begin{center}
\begingroup%
  \makeatletter%
  \providecommand\color[2][]{%
    \errmessage{(Inkscape) Color is used for the text in Inkscape, but the package 'color.sty' is not loaded}%
    \renewcommand\color[2][]{}%
  }%
  \providecommand\transparent[1]{%
    \errmessage{(Inkscape) Transparency is used (non-zero) for the text in Inkscape, but the package 'transparent.sty' is not loaded}%
    \renewcommand\transparent[1]{}%
  }%
  \providecommand\rotatebox[2]{#2}%
  \ifx\svgwidth\undefined%
    \setlength{\unitlength}{287.77614746bp}%
    \ifx\svgscale\undefined%
      \relax%
    \else%
      \setlength{\unitlength}{\unitlength * \real{\svgscale}}%
    \fi%
  \else%
    \setlength{\unitlength}{\svgwidth}%
  \fi%
  \global\let\svgwidth\undefined%
  \global\let\svgscale\undefined%
  \makeatother%
  \begin{picture}(1,0.38586472)%
    \put(0,0){\includegraphics[width=\unitlength]{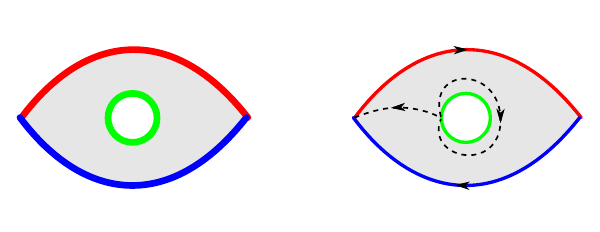}}%
    \put(0.22340294,0.33694432){\color[rgb]{0,0,0}\makebox(0,0)[b]{\smash{$B$}}}%
    \put(0.22340294,0.03188513){\color[rgb]{0,0,0}\makebox(0,0)[b]{\smash{$A$}}}%
    \put(0.77939064,0.33694432){\color[rgb]{0,0,0}\makebox(0,0)[b]{\smash{$b$}}}%
    \put(0.77939064,0.03188513){\color[rgb]{0,0,0}\makebox(0,0)[b]{\smash{$a$}}}%
    \put(0.84193035,0.18215671){\color[rgb]{0,0,0}\makebox(0,0)[lb]{\smash{$c$}}}%
    \put(0.66517833,0.17790569){\color[rgb]{0,0,0}\makebox(0,0)[b]{\smash{$g$}}}%
    \put(0.27482289,0.18215671){\color[rgb]{0,0,0}\makebox(0,0)[lb]{\smash{$C$}}}%
  \end{picture}%
\endgroup%

\end{center}
The holonomies along the arcs pictured above identify the moduli space with
$$\{(a,b,c,g)\in A\times B\times C^\perp\times G;\  abg=gc\}/C$$
where $c'\in C$ acts by
$$c'\cdot(a,b,c,g)= (a,b,c'cc'^{-1},gc'^{-1}).$$

Applying \cref{cor:embedding} to the embedding of quilted surfaces pictured below
\begin{center}
\begingroup%
  \makeatletter%
  \providecommand\color[2][]{%
    \errmessage{(Inkscape) Color is used for the text in Inkscape, but the package 'color.sty' is not loaded}%
    \renewcommand\color[2][]{}%
  }%
  \providecommand\transparent[1]{%
    \errmessage{(Inkscape) Transparency is used (non-zero) for the text in Inkscape, but the package 'transparent.sty' is not loaded}%
    \renewcommand\transparent[1]{}%
  }%
  \providecommand\rotatebox[2]{#2}%
  \ifx\svgwidth\undefined%
    \setlength{\unitlength}{320.57614746bp}%
    \ifx\svgscale\undefined%
      \relax%
    \else%
      \setlength{\unitlength}{\unitlength * \real{\svgscale}}%
    \fi%
  \else%
    \setlength{\unitlength}{\svgwidth}%
  \fi%
  \global\let\svgwidth\undefined%
  \global\let\svgscale\undefined%
  \makeatother%
  \begin{picture}(1,0.26324211)%
    \put(0,0){\includegraphics[width=\unitlength]{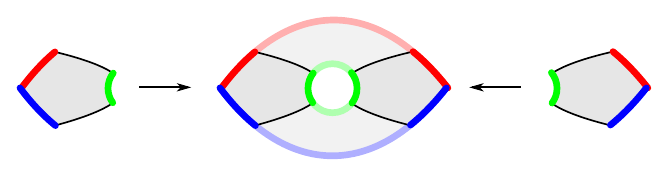}}%
  \end{picture}%
\endgroup%

\end{center}
yields the following source and target maps onto $G/C$:
$$t\times s:(a,b,c,g)\to (g,bg).$$

That the Lu-Yakimov space, $G/C$, carries a Poisson action of the Drinfeld Double, $G$, can be seen from the following sequence of embeddings (cf. Example~\ref{ex:PoissLie2}):
\begin{center}
\def\svgwidth{\linewidth}
\begingroup%
  \makeatletter%
  \providecommand\color[2][]{%
    \errmessage{(Inkscape) Color is used for the text in Inkscape, but the package 'color.sty' is not loaded}%
    \renewcommand\color[2][]{}%
  }%
  \providecommand\transparent[1]{%
    \errmessage{(Inkscape) Transparency is used (non-zero) for the text in Inkscape, but the package 'transparent.sty' is not loaded}%
    \renewcommand\transparent[1]{}%
  }%
  \providecommand\rotatebox[2]{#2}%
  \ifx\svgwidth\undefined%
    \setlength{\unitlength}{790.54555664bp}%
    \ifx\svgscale\undefined%
      \relax%
    \else%
      \setlength{\unitlength}{\unitlength * \real{\svgscale}}%
    \fi%
  \else%
    \setlength{\unitlength}{\svgwidth}%
  \fi%
  \global\let\svgwidth\undefined%
  \global\let\svgscale\undefined%
  \makeatother%
  \begin{picture}(1,0.170322)%
    \put(0,0){\includegraphics[width=\unitlength]{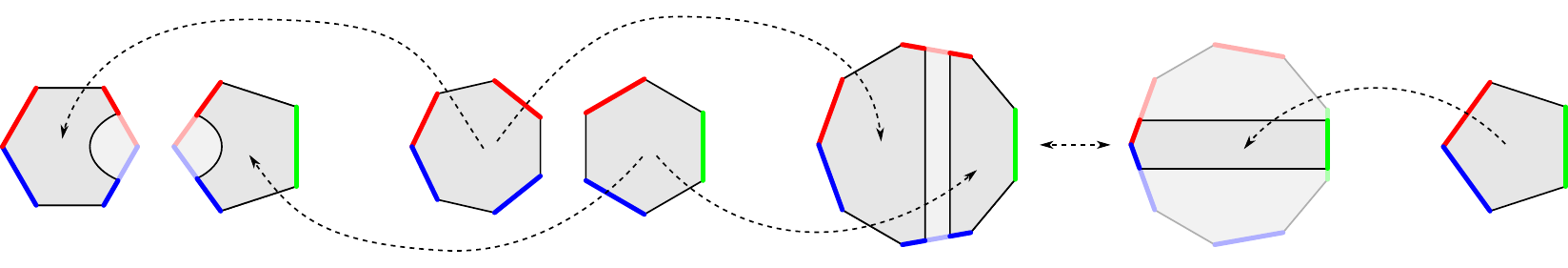}}%
    \put(0.68595035,0.08383988){\color[rgb]{0,0,0}\makebox(0,0)[b]{\smash{$=$}}}%
  \end{picture}%
\endgroup%

\end{center}
\end{example}

\begin{example}[Poisson Homogeneous spaces]
Suppose once again that $(\g,\mf{a},\mf{b})$ is a Manin triple, and that $A,B\subset G$ are Lie groups integrating $\mf{a}$ and $\mf{b}$ such that $A\cap B=\mbf{1}$. As explained in Example~\ref{ex:PoissLie2}, $A$ carries the structure of a Poisson Lie group. 
	Following Drinfel'd's classification \cite{Drinfeld:1993il}, pointed Poisson homogeneous spaces for $A$ with connected stabilizers are classified, up to isomorphism, by Lagrangian subalgebras $\mf{c}\subseteq \g$.\footnote{Note that not every Lagrangian subalgebra $\mf{c}$ corresponds to a Poisson homogeneous space, but the correspondences become one-to-one when one works with local lie groups and homogeneous spaces (cf. \cite{Drinfeld:1993il,Evens:2001ue,Evens:2006kk}). \\
	Alternatively (if one doesn't wish to work with pointed spaces), Poisson homogeneous spaces for $A$ with connected stabilizers are classified, up to isomorphism, by $A$-orbits of Lagrangian subalgebras $\mf{c}\subseteq \g$.} Let $C\subset G$ be a subgroup whose corresponding subalgebra $\mf{c}\subseteq \g$ is Lagrangian. Let $K=A\cap C$. Consider the quilted surface pictured below:
	\begin{center}
\begingroup%
  \makeatletter%
  \providecommand\color[2][]{%
    \errmessage{(Inkscape) Color is used for the text in Inkscape, but the package 'color.sty' is not loaded}%
    \renewcommand\color[2][]{}%
  }%
  \providecommand\transparent[1]{%
    \errmessage{(Inkscape) Transparency is used (non-zero) for the text in Inkscape, but the package 'transparent.sty' is not loaded}%
    \renewcommand\transparent[1]{}%
  }%
  \providecommand\rotatebox[2]{#2}%
  \ifx\svgwidth\undefined%
    \setlength{\unitlength}{227.45605469bp}%
    \ifx\svgscale\undefined%
      \relax%
    \else%
      \setlength{\unitlength}{\unitlength * \real{\svgscale}}%
    \fi%
  \else%
    \setlength{\unitlength}{\svgwidth}%
  \fi%
  \global\let\svgwidth\undefined%
  \global\let\svgscale\undefined%
  \makeatother%
  \begin{picture}(1,0.36446893)%
    \put(0,0){\includegraphics[width=\unitlength]{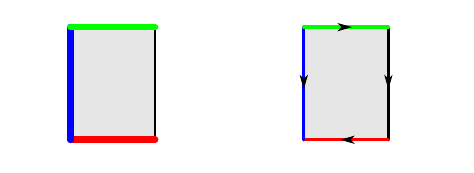}}%
    \put(0.7327972,0.33237344){\color[rgb]{0,0,0}\makebox(0,0)[b]{\smash{$c$}}}%
    \put(0.7327972,0.03298525){\color[rgb]{0,0,0}\makebox(0,0)[b]{\smash{$b$}}}%
    \put(0.83971593,0.19257909){\color[rgb]{0,0,0}\makebox(0,0)[lb]{\smash{$g$}}}%
    \put(0.61752715,0.19257909){\color[rgb]{0,0,0}\makebox(0,0)[rb]{\smash{$a$}}}%
    \put(0.24039439,0.02595093){\color[rgb]{0,0,0}\makebox(0,0)[b]{\smash{$B$}}}%
    \put(0.34731312,0.19257909){\color[rgb]{0,0,0}\makebox(0,0)[lb]{\smash{$G$}}}%
    \put(0.12512434,0.19257909){\color[rgb]{0,0,0}\makebox(0,0)[rb]{\smash{$A$}}}%
    \put(0.24039439,0.33546126){\color[rgb]{0,0,0}\makebox(0,0)[b]{\smash{$C$}}}%
  \end{picture}%
\endgroup%

\end{center}
Following \cref{cor:PartTriv}, the moduli space is
$$\{(c,g,b,a)\in C\times G\times B \times A;\  bgc=a\}/\sim,$$
where the holonomies are as pictured above, and the quotient is by the group of residual gauge transformations.
Since $A\cap B=\mbf{1}$, the group of residual gauge transformations is $K\times C\times B$  acting  by 
$$(k,c',b')\cdot(c,g,b,a)=(c'ck^{-1},b'gc'^{-1},bb'^{-1},ak^{-1}).$$
Thus the moduli space is identified with $A/K$. 

The Poisson Lie group $A$ acts on $A/K$ via the following sequence of embeddings:
\begin{center}
\def\svgwidth{\linewidth}
\begingroup%
  \makeatletter%
  \providecommand\color[2][]{%
    \errmessage{(Inkscape) Color is used for the text in Inkscape, but the package 'color.sty' is not loaded}%
    \renewcommand\color[2][]{}%
  }%
  \providecommand\transparent[1]{%
    \errmessage{(Inkscape) Transparency is used (non-zero) for the text in Inkscape, but the package 'transparent.sty' is not loaded}%
    \renewcommand\transparent[1]{}%
  }%
  \providecommand\rotatebox[2]{#2}%
  \ifx\svgwidth\undefined%
    \setlength{\unitlength}{490.76020508bp}%
    \ifx\svgscale\undefined%
      \relax%
    \else%
      \setlength{\unitlength}{\unitlength * \real{\svgscale}}%
    \fi%
  \else%
    \setlength{\unitlength}{\svgwidth}%
  \fi%
  \global\let\svgwidth\undefined%
  \global\let\svgscale\undefined%
  \makeatother%
  \begin{picture}(1,0.15455567)%
    \put(0,0){\includegraphics[width=\unitlength]{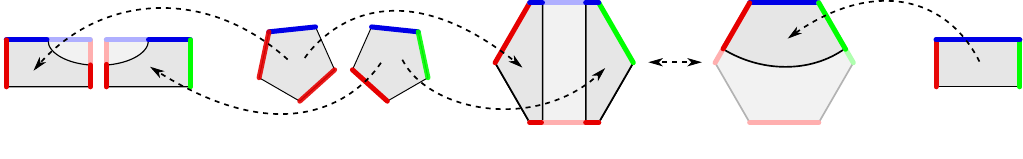}}%
    \put(0.66122378,0.09977918){\color[rgb]{0,0,0}\makebox(0,0)[b]{\smash{$=$}}}%
  \end{picture}%
\endgroup%

\end{center}

To compute the bivector field on $A/K$, we work with the following skeleton:
	\begin{center}
\begingroup%
  \makeatletter%
  \providecommand\color[2][]{%
    \errmessage{(Inkscape) Color is used for the text in Inkscape, but the package 'color.sty' is not loaded}%
    \renewcommand\color[2][]{}%
  }%
  \providecommand\transparent[1]{%
    \errmessage{(Inkscape) Transparency is used (non-zero) for the text in Inkscape, but the package 'transparent.sty' is not loaded}%
    \renewcommand\transparent[1]{}%
  }%
  \providecommand\rotatebox[2]{#2}%
  \ifx\svgwidth\undefined%
    \setlength{\unitlength}{95.49769897bp}%
    \ifx\svgscale\undefined%
      \relax%
    \else%
      \setlength{\unitlength}{\unitlength * \real{\svgscale}}%
    \fi%
  \else%
    \setlength{\unitlength}{\svgwidth}%
  \fi%
  \global\let\svgwidth\undefined%
  \global\let\svgscale\undefined%
  \makeatother%
  \begin{picture}(1,0.97904021)%
    \put(0,0){\includegraphics[width=\unitlength]{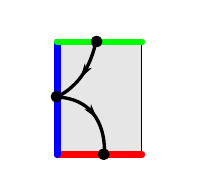}}%
    \put(0.50612855,0.09762017){\color[rgb]{0,0,0}\makebox(0,0)[b]{\smash{$B$}}}%
    \put(0.76078717,0.49449451){\color[rgb]{0,0,0}\makebox(0,0)[lb]{\smash{$G$}}}%
    \put(0.23157878,0.49449451){\color[rgb]{0,0,0}\makebox(0,0)[rb]{\smash{$A$}}}%
    \put(0.50612855,0.83481071){\color[rgb]{0,0,0}\makebox(0,0)[b]{\smash{$C$}}}%
    \put(0.52337514,0.35268936){\color[rgb]{0,0,0}\makebox(0,0)[lb]{\smash{$g_1$}}}%
    \put(0.4478279,0.58982413){\color[rgb]{0,0,0}\makebox(0,0)[lb]{\smash{$g_2$}}}%
  \end{picture}%
\endgroup%

\end{center}
Let $g_1=b_1a_1\in BA,$ and $g_2=a_2c_2\in AC$ denote the holonomies along the corresponding edges, then we may choose the following (skeletal) coordinates for $A/K$:
$$\{(ba_1,a_2c)\in BA\times AC\}/\sim\quad
\xrightarrow{(ba_1,a_2c)\to (a_1a_2)K} A/K,$$
where the equivalence relation on the left hand side is induced by the action \begin{equation}\label{eq:PHequiv}(a',b',c')\cdot (ba_1,a_2c)=(b'ba_1a'^{-1},a'a_2cc'{-1}),\quad (a',b',c')\in A\times B\times C\end{equation}
In these coordinates, the bivector is  the restriction of \eqref{eq:FR}:
$$\pi=-\frac{1}{2}\sum_i \xi_i^L(1)\wedge {\eta^i}^R(2)+{\eta^i}^L(1)\wedge {\xi_i}^R(2)$$ to $BA\times AC\subset G^2$, 
where, for $\zeta\in \g$, $\zeta^L(k)$ and $\zeta^R(k)$ denote the corresponding left and right-invariant vector fields which are tangent to the $k^{th}$ factor of $G^2$ ($k=1,2$), and $\{\xi_i\}\subset\mf{a}$, and $\{\eta^i\}\subset \mf{b}$ are dual bases. At the point $(1,1)\in G^2$, 
$$\pi_{(1,1)}\sim -\frac{1}{2}\sum_i\xi_i(2)\wedge {\eta^i}(2)$$
modulo the equivalence relation \eqref{eq:PHequiv} (since at $(1,1)$, $\eta^i\in \mf{b}\Rightarrow{\eta^i}^L(1)\sim 0$ and $\xi_i\in\mf{a}\Rightarrow\xi_i^L(1)\sim \xi_i^R(2)$).

 Now, choose a vector space complement $\mf{p}\subseteq \mf{b}$ to $\ann(\mf{k})\subseteq \mf{b}\cong \mf{a}^*$ (where $\mf{k}\subseteq \mf{a}$ is the Lie algebra of $K\subseteq A$). Then $\mf{a}=\mf{k}\oplus\ann(\mf{p})$, and $\mf{b}=\mf{p}\oplus\ann(\mf{k})$. There exists a unique $\varpi\in \wedge^2 (\mf{a}/\mf{k})\cong\wedge^2\ann(\mf{p})$ such that 
$$\mf{c}=\{(\xi+\varpi^\sharp(\eta)+\eta)\in \g\mid \xi\in\mf{k},\eta\in \ann(\mf{k})\},$$
where $\varpi^\sharp(\eta)= (\eta\otimes \on{id})(\varpi)$ and we have used the identification $\varpi\in\wedge^2(\mf{a}/\mf{k})\cong\wedge^2\big(\ann(\mf{k})^*\big)$.

Let $\{\xi_i^1\}\subset \mf{k}$ and $\{\eta^i_1\}\subset\mf{p}$ be bases in duality, and similarly with $\{\xi_i^0\}\subset \ann(\mf{p})$ and $\{\eta^i_0\}\subset\ann(\mf{k})$. Then
$$\pi_{(1,1)}\sim -\frac{1}{2}\sum_i\xi_i(2)\wedge {\eta^i}(2)=-\frac{1}{2}\sum_i{\xi_i^0}(2)\wedge {\eta^i_0}(2)+{\xi_i^1}(2)\wedge {\eta^i_1}(2)$$
Modulo \eqref{eq:PHequiv}, we have 
$$\pi_{(1,1)}\sim \frac{1}{2}\sum_i{\xi_i^0}(2)\wedge {\varpi^\sharp(\eta^i_0)}(2)$$
(since $\xi_i^1\in \mf{c}\Rightarrow{\xi_i^1}(2)\sim 0$ and $\varpi^\sharp(\eta^i_0)+\eta_0^i\in \mf{c}\Rightarrow {\eta^i_0}(2)\sim -{\varpi^\sharp(\eta^i_0)}(2)$).
Thus, at the point $K\in A/K$, the bivector is just $$\varpi\in \wedge^2 (\mf{a}/\mf{k})\cong\wedge^2T_K(A/K),$$
which shows that the pointed Poisson homogeneous space $(A/K,\ast=K)$ corresponds to the Lagrangian Lie subalgebra $\mf{c}\subset \g$ under Drinfel'd's classification. 
	
\end{example}

\begin{example}[Drinfel'd Polyubles \cite{Fock:1999wz}]\label{ex:poly}
Suppose once again that $(\g,\mf{a},\mf{b})$ is a Manin triple, and that $A,B\subset G$ are Lie groups integrating $\mf{a}$ and $\mf{b}$ such that $A\cap B=\mbf{1}$. The \emph{$n^{th}$ Drinfel'd polyuble of $A$} (introduced in \cite{Fock:1999wz} for factorizable Poisson Lie groups, $A$, and generalized to arbitrary Poisson Lie groups by Jiang Hua-Lu and Victor Mouquin \cite{Lu:FBHentwA}) is a Lie-Poisson structure on 
$$\mc{P}_n(A):=\begin{cases}A\times G^k &\text{when $n=2k+1$ is odd, and}\\ G^{k+1}&\text{when $n=2k+2$ is even.}\end{cases}$$
The Poisson manifold $\mc{P}_n(A)$ is naturally identified with the moduli space $\mc{M}_{\Sigma_n}$ for a quilted surface $\Sigma_n$, as we now describe.\footnote{For the factorizable case, this is essentially explained in \cite{Fock:1999wz}, but using a slightly different language.} When $n$ is even, the quilted surface $\Sigma_n$ is the $k=\frac{n-2}{2}$-fold punctured disk pictured below left,

\noindent%
\begin{minipage}{\linewidth}
\makebox[\linewidth]{
\begingroup%
  \makeatletter%
  \providecommand\color[2][]{%
    \errmessage{(Inkscape) Color is used for the text in Inkscape, but the package 'color.sty' is not loaded}%
    \renewcommand\color[2][]{}%
  }%
  \providecommand\transparent[1]{%
    \errmessage{(Inkscape) Transparency is used (non-zero) for the text in Inkscape, but the package 'transparent.sty' is not loaded}%
    \renewcommand\transparent[1]{}%
  }%
  \providecommand\rotatebox[2]{#2}%
  \ifx\svgwidth\undefined%
    \setlength{\unitlength}{255.83703613bp}%
    \ifx\svgscale\undefined%
      \relax%
    \else%
      \setlength{\unitlength}{\unitlength * \real{\svgscale}}%
    \fi%
  \else%
    \setlength{\unitlength}{\svgwidth}%
  \fi%
  \global\let\svgwidth\undefined%
  \global\let\svgscale\undefined%
  \makeatother%
  \begin{picture}(1,0.48142815)%
    \put(0,0){\includegraphics[width=\unitlength]{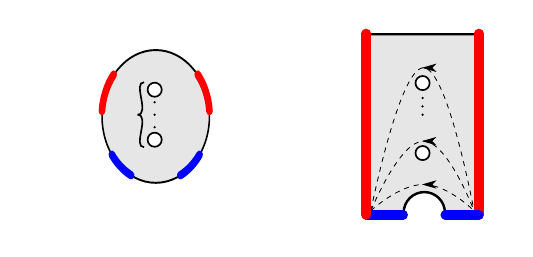}}%
    \put(0.39941937,0.21856426){\color[rgb]{0,0,0}\makebox(0,0)[lb]{\smash{$\mbf{1}$}}}%
    \put(0.40784692,0.29980236){\color[rgb]{0,0,0}\makebox(0,0)[lb]{\smash{$B$}}}%
    \put(0.17627383,0.29980236){\color[rgb]{0,0,0}\makebox(0,0)[rb]{\smash{$B$}}}%
    \put(0.18469093,0.21856426){\color[rgb]{0,0,0}\makebox(0,0)[rb]{\smash{$\mbf{1}$}}}%
    \put(0.67219471,0.23726255){\color[rgb]{0,0,0}\makebox(0,0)[rb]{\smash{$B$}}}%
    \put(0.91334248,0.23726255){\color[rgb]{0,0,0}\makebox(0,0)[lb]{\smash{$B$}}}%
    \put(0.24728218,0.26107703){\color[rgb]{0,0,0}\makebox(0,0)[rb]{\smash{$k$}}}%
    \put(0.79276869,0.42759081){\color[rgb]{0,0,0}\makebox(0,0)[b]{\smash{$G$}}}%
    \put(0.7239749,0.03643919){\color[rgb]{0,0,0}\makebox(0,0)[b]{\smash{$A$}}}%
    \put(0.79411394,0.15739727){\color[rgb]{0,0,0}\makebox(0,0)[b]{\smash{$g_0$}}}%
    \put(0.79252745,0.23543386){\color[rgb]{0,0,0}\makebox(0,0)[b]{\smash{$g_1$}}}%
    \put(0.79335471,0.37675746){\color[rgb]{0,0,0}\makebox(0,0)[b]{\smash{$g_k$}}}%
    \put(0.20643391,0.15526484){\color[rgb]{0,0,0}\makebox(0,0)[rb]{\smash{$A$}}}%
    \put(0.38472047,0.15526484){\color[rgb]{0,0,0}\makebox(0,0)[lb]{\smash{$A$}}}%
    \put(0.86781646,0.03643919){\color[rgb]{0,0,0}\makebox(0,0)[b]{\smash{$A$}}}%
  \end{picture}%
\endgroup%
}
\captionof{figure}{}\label{fig:Polyuble1}
\end{minipage}
contracting the uncolored boundary arcs marked by $\mbf{1}$ results in the surface pictured above right. Since $A\cap B=\mbf{1}$, the group $(\hat H\cap \hat C)$ of residual gauge transformations is trivial (cf. \cref{cor:PartTriv}). Thus, computing the holonomies along the dotted arcs canonically identifies $\mc{M}_{\Sigma_n}$ with $$\{(g_0,g_1,\dots,g_k)\in G^{k+1}\}=\mc{P}_n(A).$$
The resulting Poisson structure on $\mc{P}_n(A)$ is the Drinfel'd Polyuble Poisson Lie structure, as explained in \cite{Fock:1999wz} (cf. Remark~\ref{rem:FR}). In particular, for $n=2$, the quilted surface pictured in Figure~\ref{fig:Polyuble1} is just the one from Example~\ref{ex:DblPoissLie} corresponding to the Drinfel'd Double. 
 The double symplectic groupoid integrating $\mc{P}_n(A)$ is the moduli space for the quilted surface from Example~\ref{ex:symplDblDbl}, modified by attaching $k$ handles in the interior.

On the other hand, when $n$ is odd, then the quilted surface $\Sigma_n$ is the $k=\frac{n-1}{2}$-fold punctured disk pictured below left,
\begin{center}
\begingroup%
  \makeatletter%
  \providecommand\color[2][]{%
    \errmessage{(Inkscape) Color is used for the text in Inkscape, but the package 'color.sty' is not loaded}%
    \renewcommand\color[2][]{}%
  }%
  \providecommand\transparent[1]{%
    \errmessage{(Inkscape) Transparency is used (non-zero) for the text in Inkscape, but the package 'transparent.sty' is not loaded}%
    \renewcommand\transparent[1]{}%
  }%
  \providecommand\rotatebox[2]{#2}%
  \ifx\svgwidth\undefined%
    \setlength{\unitlength}{255.83703613bp}%
    \ifx\svgscale\undefined%
      \relax%
    \else%
      \setlength{\unitlength}{\unitlength * \real{\svgscale}}%
    \fi%
  \else%
    \setlength{\unitlength}{\svgwidth}%
  \fi%
  \global\let\svgwidth\undefined%
  \global\let\svgscale\undefined%
  \makeatother%
  \begin{picture}(1,0.48142815)%
    \put(0,0){\includegraphics[width=\unitlength]{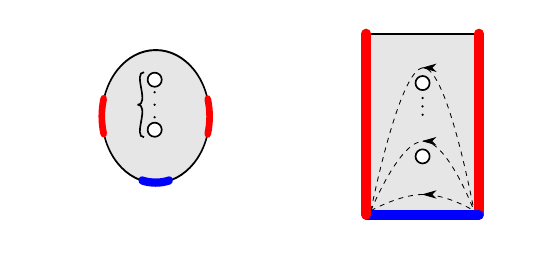}}%
    \put(0.38065742,0.16227843){\color[rgb]{0,0,0}\makebox(0,0)[lb]{\smash{$\mbf{1}$}}}%
    \put(0.40784692,0.24977051){\color[rgb]{0,0,0}\makebox(0,0)[lb]{\smash{$B$}}}%
    \put(0.17627383,0.24977051){\color[rgb]{0,0,0}\makebox(0,0)[rb]{\smash{$B$}}}%
    \put(0.20345288,0.16227843){\color[rgb]{0,0,0}\makebox(0,0)[rb]{\smash{$\mbf{1}$}}}%
    \put(0.67219471,0.23726255){\color[rgb]{0,0,0}\makebox(0,0)[rb]{\smash{$B$}}}%
    \put(0.91334248,0.23726255){\color[rgb]{0,0,0}\makebox(0,0)[lb]{\smash{$B$}}}%
    \put(0.24728218,0.27983898){\color[rgb]{0,0,0}\makebox(0,0)[rb]{\smash{$k$}}}%
    \put(0.79276869,0.42759081){\color[rgb]{0,0,0}\makebox(0,0)[b]{\smash{$G$}}}%
    \put(0.79276869,0.03643919){\color[rgb]{0,0,0}\makebox(0,0)[b]{\smash{$A$}}}%
    \put(0.82263093,0.13238135){\color[rgb]{0,0,0}\makebox(0,0)[rb]{\smash{$a$}}}%
    \put(0.79252745,0.23543386){\color[rgb]{0,0,0}\makebox(0,0)[b]{\smash{$g_1$}}}%
    \put(0.79335471,0.37675746){\color[rgb]{0,0,0}\makebox(0,0)[b]{\smash{$g_k$}}}%
    \put(0.2924502,0.09272502){\color[rgb]{0,0,0}\makebox(0,0)[b]{\smash{$A$}}}%
  \end{picture}%
\endgroup%

\end{center} contracting the uncolored boundary arcs marked by $\mbf{1}$ results in the surface pictured above right. As before, the group of residual gauge transformations is trivial, and computing the holonomies along the dotted arcs canonically identifies $\mc{M}_{\Sigma_n}$ with $$\{(a_0,g_1,\dots,g_k)\in A\times G^k\}=\mc{P}_n(A).$$
The resulting Poisson structure on $\mc{P}_n(A)$ is the Drinfel'd Polyuble Poisson Lie structure, as explained in \cite{Fock:1999wz} (cf. Remark~\ref{rem:FR}).  The double symplectic groupoid integrating $\mc{P}_n(A)$ is the moduli space for the quilted surface pictured in Figure~\ref{fig:SymplDbl}, modified by attaching $k$ handles in the interior.

As mentioned in  \cite{Fock:1999wz}, for $n$ either even or odd, the natural embedding of $\Sigma_n$ into the (unpunctured) disk from Example~\ref{ex:PoissLie2} yields a Poisson morphism $A\to \mc{P}_n(A)$ between the corresponding moduli spaces (cf. Example~\ref{ex:DblPoissLie} for the case $n=2$).

To take a different perspective on this example, we recall that the $n^{th}$ Drinfeld' polyuble is the Poisson Lie group whose corresponding Manin triple is $(\g_n,\mf{a}_n,\mf{b}_n)$, where 

\begin{align}
\label{eq:DrPolygn}\g_n&:=\begin{cases}
\g\oplus(\bar\g\oplus\g)^{\frac{n-1}{2}},&\text{if $n$ is odd},\\
(\g\oplus\bar\g)^{\frac{n}{2}},&\text{if $n$ is even},
\end{cases}\\
\notag\mf{a}_n&:=\begin{cases}
\mf{a}\oplus(\g_\Delta)^{\frac{n-1}{2}},&\text{if $n$ is odd},\\
(\g_\Delta)^{\frac{n}{2}},&\text{if $n$ is even},
\end{cases}\\
\notag\mf{b}_n&:=\begin{cases}
(\g_\Delta)^{\frac{n-1}{2}}\oplus\mf{b},&\text{if $n$ is odd},\\
\mf{a}\oplus(\g_\Delta)^{\frac{n-2}{2}}\oplus\mf{b},&\text{if $n$ is even},
\end{cases}
\end{align}

Where $\g_\Delta=\{(\xi,\xi)\}\subseteq\g\oplus\bar\g$ is the diagonal subalgebra, as in Example~\ref{ex:FuseRedIsRed}.

Thus, after first replacing the Manin triple $(\g,\mf{a},\mf{b})$ with $(\g_n,\mf{a}_n,\mf{b}_n)$, we may
 identify the $n^{th}$ Drinfel'd polyuble with the moduli space for the quilted surface from Example~\ref{ex:PoissLie2}:

\noindent%
\begin{minipage}{\linewidth}
\makebox[\linewidth]{
  \begin{tabular}{p{.3\textwidth}p{.55\textwidth}}
\begin{center}
\begingroup%
  \makeatletter%
  \providecommand\color[2][]{%
    \errmessage{(Inkscape) Color is used for the text in Inkscape, but the package 'color.sty' is not loaded}%
    \renewcommand\color[2][]{}%
  }%
  \providecommand\transparent[1]{%
    \errmessage{(Inkscape) Transparency is used (non-zero) for the text in Inkscape, but the package 'transparent.sty' is not loaded}%
    \renewcommand\transparent[1]{}%
  }%
  \providecommand\rotatebox[2]{#2}%
  \ifx\svgwidth\undefined%
    \setlength{\unitlength}{89.67727051bp}%
    \ifx\svgscale\undefined%
      \relax%
    \else%
      \setlength{\unitlength}{\unitlength * \real{\svgscale}}%
    \fi%
  \else%
    \setlength{\unitlength}{\svgwidth}%
  \fi%
  \global\let\svgwidth\undefined%
  \global\let\svgscale\undefined%
  \makeatother%
  \begin{picture}(1,0.8820082)%
    \put(0,0){\includegraphics[width=\unitlength]{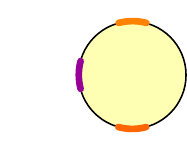}}%
    \put(0.69896567,0.01474736){\color[rgb]{0,0,0}\makebox(0,0)[b]{\smash{$B_n$}}}%
    \put(0.3651795,0.46159656){\color[rgb]{0,0,0}\makebox(0,0)[rb]{\smash{$A_n$}}}%
    \put(0.49121454,0.21868857){\color[rgb]{0,0,0}\makebox(0,0)[rb]{\smash{$\mbf{1}$}}}%
    \put(0.69896567,0.81762639){\color[rgb]{0,0,0}\makebox(0,0)[b]{\smash{$B_n$}}}%
    \put(0.49121454,0.70041599){\color[rgb]{0,0,0}\makebox(0,0)[rb]{\smash{$\mbf{1}$}}}%
    \put(0.69896567,0.46079126){\color[rgb]{0,0,0}\makebox(0,0)[b]{\smash{$G_n$}}}%
  \end{picture}%
\endgroup%
\end{center}&
\begin{center}$\begin{aligned}
G_n&:=G^n\\
A_n&:=\begin{cases}
A\times (G_\Delta)^{\frac{n-1}{2}},&\text{if $n$ is odd},\\
(G_\Delta)^{\frac{n}{2}},&\text{if $n$ is even},
\end{cases}\\
B_n&:=\begin{cases}
(G_\Delta)^{\frac{n-1}{2}}\times B,&\text{if $n$ is odd},\\
A\times(G_\Delta)^{\frac{n-2}{2}}\times B,&\text{if $n$ is even},
\end{cases}
\end{aligned}$\end{center}
\end{tabular}
}
\captionof{figure}{}\label{fig:Polyuble2}
\end{minipage}
where $G_\Delta$ is the diagonal subgroup, as in Example~\ref{ex:FuseRedIsRed}. In particular, using the multiplicative structure described in Example~\ref{ex:PoissLie2}, this shows that $\mc{P}_n$ is a Poisson Lie group, as claimed.

We now relate the $G_n$-colored quilted surface pictured in Figure~\ref{fig:Polyuble2} to the $G$-colored quilted surface pictured in Figure~\ref{fig:Polyuble1}. For brevity, we will consider only the case where $n$ is even (the odd case is completely analogous). 

Let $(\Sigma,V)$ be the underlying domain for the quilted surface pictured in Figure~\ref{fig:Polyuble2} (i.e. the disk with three marked points on the boundary).
First, notice that coloring a domain $\Sigma$ with $G_n=G^n$ is equivalent to coloring $n$ copies of the domain with $G$. 

\begin{center}
\begingroup%
  \makeatletter%
  \providecommand\color[2][]{%
    \errmessage{(Inkscape) Color is used for the text in Inkscape, but the package 'color.sty' is not loaded}%
    \renewcommand\color[2][]{}%
  }%
  \providecommand\transparent[1]{%
    \errmessage{(Inkscape) Transparency is used (non-zero) for the text in Inkscape, but the package 'transparent.sty' is not loaded}%
    \renewcommand\transparent[1]{}%
  }%
  \providecommand\rotatebox[2]{#2}%
  \ifx\svgwidth\undefined%
    \setlength{\unitlength}{200.54492188bp}%
    \ifx\svgscale\undefined%
      \relax%
    \else%
      \setlength{\unitlength}{\unitlength * \real{\svgscale}}%
    \fi%
  \else%
    \setlength{\unitlength}{\svgwidth}%
  \fi%
  \global\let\svgwidth\undefined%
  \global\let\svgscale\undefined%
  \makeatother%
  \begin{picture}(1,0.40056143)%
    \put(0,0){\includegraphics[width=\unitlength]{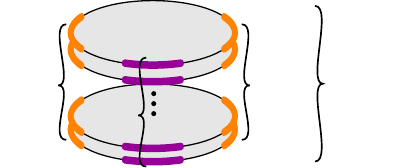}}%
    \put(0.78630491,0.19594158){\color[rgb]{0,0,0}\makebox(0,0)[lb]{\smash{$n$ layers of $G$-colored sheets}}}%
    \put(0.11795285,0.1869132){\color[rgb]{0,0,0}\makebox(0,0)[rb]{\smash{$B_n$}}}%
    \put(0.61260508,0.1869132){\color[rgb]{0,0,0}\makebox(0,0)[lb]{\smash{$B_n$}}}%
    \put(0.32538767,0.11510884){\color[rgb]{0,0,0}\makebox(0,0)[rb]{\smash{$A_n$}}}%
    \put(0.36372923,0.31269856){\color[rgb]{0,0,0}\makebox(0,0)[b]{\smash{$G$}}}%
  \end{picture}%
\endgroup%

\end{center}

Next, we use the fact (cf. Example~\ref{ex:FuseRedIsRed}) that coloring a pair of marked points by $G_\Delta$ corresponds to fusing those two marked points and then reducing by $G$ (effectively erasing the marked point); that is, coloring a domain wall by $G_\Delta$ corresponds to sewing the surface together along that domain wall:
\begin{center}
\begingroup%
  \makeatletter%
  \providecommand\color[2][]{%
    \errmessage{(Inkscape) Color is used for the text in Inkscape, but the package 'color.sty' is not loaded}%
    \renewcommand\color[2][]{}%
  }%
  \providecommand\transparent[1]{%
    \errmessage{(Inkscape) Transparency is used (non-zero) for the text in Inkscape, but the package 'transparent.sty' is not loaded}%
    \renewcommand\transparent[1]{}%
  }%
  \providecommand\rotatebox[2]{#2}%
  \ifx\svgwidth\undefined%
    \setlength{\unitlength}{280.00681152bp}%
    \ifx\svgscale\undefined%
      \relax%
    \else%
      \setlength{\unitlength}{\unitlength * \real{\svgscale}}%
    \fi%
  \else%
    \setlength{\unitlength}{\svgwidth}%
  \fi%
  \global\let\svgwidth\undefined%
  \global\let\svgscale\undefined%
  \makeatother%
  \begin{picture}(1,0.33210345)%
    \put(0,0){\includegraphics[width=\unitlength]{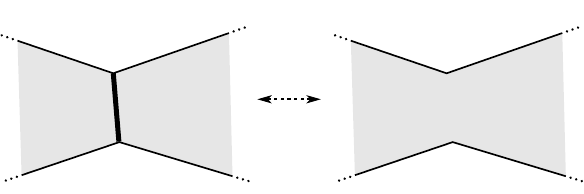}}%
    \put(0.19769747,0.28726079){\color[rgb]{0,0,0}\makebox(0,0)[b]{\smash{$G_\Delta$-colored
Domain wall.}}}%
    \put(0.76911214,0.28154665){\color[rgb]{0,0,0}\makebox(0,0)[b]{\smash{No domain wall.}}}%
    \put(0.497852,0.17728506){\color[rgb]{0,0,0}\makebox(0,0)[b]{\smash{$=$}}}%
  \end{picture}%
\endgroup%

\end{center}
 Since $A_n= G_\Delta^{\frac{n}{2}}$ is just $\frac{n}{2}$ copies of $G_\Delta$, reducing by $A_n$ corresponds to sewing the alternating pairs of sheets together at that marked point, as pictured below (for the case $n=6$): 

\begin{minipage}{.9\linewidth}
\makebox[.9\linewidth]{
\begingroup%
  \makeatletter%
  \providecommand\color[2][]{%
    \errmessage{(Inkscape) Color is used for the text in Inkscape, but the package 'color.sty' is not loaded}%
    \renewcommand\color[2][]{}%
  }%
  \providecommand\transparent[1]{%
    \errmessage{(Inkscape) Transparency is used (non-zero) for the text in Inkscape, but the package 'transparent.sty' is not loaded}%
    \renewcommand\transparent[1]{}%
  }%
  \providecommand\rotatebox[2]{#2}%
  \ifx\svgwidth\undefined%
    \setlength{\unitlength}{200.54492188bp}%
    \ifx\svgscale\undefined%
      \relax%
    \else%
      \setlength{\unitlength}{\unitlength * \real{\svgscale}}%
    \fi%
  \else%
    \setlength{\unitlength}{\svgwidth}%
  \fi%
  \global\let\svgwidth\undefined%
  \global\let\svgscale\undefined%
  \makeatother%
  \begin{picture}(1,0.40056143)%
    \put(0,0){\includegraphics[width=\unitlength]{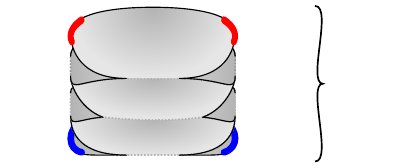}}%
    \put(0.78630491,0.19594158){\color[rgb]{0,0,0}\makebox(0,0)[lb]{\smash{$6$ layers of $G$-colored sheets}}}%
    \put(0.15784416,0.22680452){\color[rgb]{0,0,0}\makebox(0,0)[rb]{\smash{$G_\Delta$}}}%
    \put(0.57271377,0.22680452){\color[rgb]{0,0,0}\makebox(0,0)[lb]{\smash{$G_\Delta$}}}%
    \put(0.36372922,0.20898113){\color[rgb]{0,0,0}\makebox(0,0)[b]{\smash{$G_\Delta$}}}%
    \put(0.36372922,0.01750283){\color[rgb]{0,0,0}\makebox(0,0)[b]{\smash{$G_\Delta$}}}%
    \put(0.57271377,0.32254366){\color[rgb]{0,0,0}\makebox(0,0)[lb]{\smash{$B$}}}%
    \put(0.57271377,0.13904363){\color[rgb]{0,0,0}\makebox(0,0)[lb]{\smash{$G_\Delta$}}}%
    \put(0.57271377,0.04330448){\color[rgb]{0,0,0}\makebox(0,0)[lb]{\smash{$A$}}}%
    \put(0.15784416,0.32254366){\color[rgb]{0,0,0}\makebox(0,0)[rb]{\smash{$B$}}}%
    \put(0.15784416,0.13904363){\color[rgb]{0,0,0}\makebox(0,0)[rb]{\smash{$G_\Delta$}}}%
    \put(0.15784416,0.04330448){\color[rgb]{0,0,0}\makebox(0,0)[rb]{\smash{$A$}}}%
    \put(0.36372922,0.10526372){\color[rgb]{0,0,0}\makebox(0,0)[b]{\smash{$G_\Delta$}}}%
  \end{picture}%
\endgroup%
}
\captionof{figure}{We sew the sheets together along the domain walls coloured by $G_\Delta$, so that the orientations coincide when one crosses the $G_\Delta$-coloured domain wall. Note that the orientations of  alternating layers are reversed - this corresponds to the fact that alternating layers (domains) are coloured with $\bar\g$ rather than $\g$ (cf. \eqref{eq:DrPolygn}).}\label{fig:sewing}
\end{minipage}

Similarly, since $B_n= A\times(G_\Delta)^{\frac{n-2}{2}}\times B$ is just $\frac{n-2}{2}$ copies of $G_\Delta$ (bounded on one end by a copy of $A$, and on the other end by a copy of $B$) reducing by $B_n$ corresponds to sewing alternating pairs of the interior sheets together at that marked point (as pictured above), and then coloring the marked points on the initial (respectively, final) sheet by $A$ (respectively $B$). Unfolding the resulting quilted surface to lay it flat yields the quilted surface pictured in Figure~\ref{fig:Polyuble1}.

\end{example}

\begin{example}[Fission spaces \cite{Boalch:2011vt}]\label{ex:FissionSpace}
Suppose that $s\in S^2(\g)^\g$ is non-degenerate, $\g=\mf{u}_+\oplus\h\oplus\mf{u}_-$ as a vector space (but not as a Lie algebra), where $\mf{b}_\pm:=\h\oplus\mf{u}_\pm\subseteq \g$ are coisotropic subalgebras satisfying $\mf{b}_\pm^\perp=\mf{u}_\pm$. Suppose further that the Lie subalgebras $\mf{u}_\pm,\mf{b}_\pm,\h$ all integrate to closed subgroups $U_\pm,B_\pm,H\subseteq G$ such that $H=B_+\cap B_-$. The metric on $\g$ descends to a non-degenerate invariant metric on $\h\subseteq\g$, and
\begin{equation}\label{eq:BoalchCois}C_\pm:=\{(c,\tilde c)\in H\times G;\  \tilde c c^{-1}\in U_\pm\}\end{equation}
is a coisotropic subgroup of $H\times G$.

\begin{figure}[h!]
\begin{center}
\def\svgwidth{\linewidth}
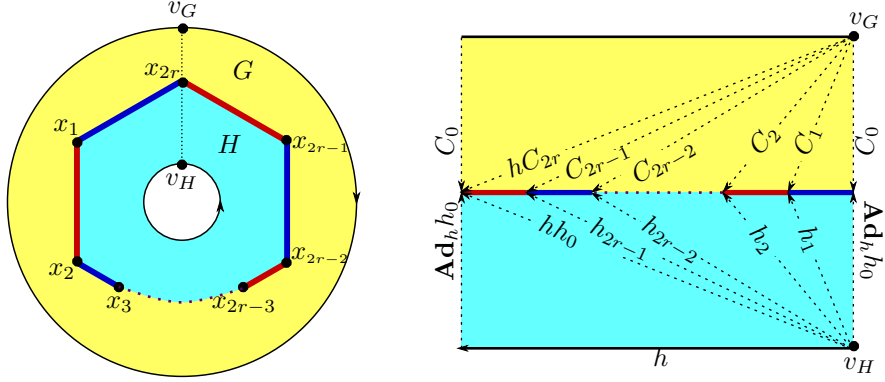
\caption{\label{fig:PBoalchDrawCut} On the surface pictured above, the structure group in the yellow domain is $G$ while the structure group in the blue domain is $H$. Along the boundary of the two domains, blue edges are colored with $C_+$ while the red edges are colored with $C_-$. The dotted segment of the boundary between the two annuli is meant to indicate an alternating sequence of blue and red edges.
 Cutting along the vertical dotted line in the first picture yields the second picture. Acting by $H$ at the points $x_1,\dots, x_{2r}$ allows us to set the holonomies $h_0,\dots,h_{2r-1}$ to the identity.}
\end{center}
\end{figure}

Consider the quilted surface pictured in \cref{fig:PBoalchDrawCut}. Computing holonomies along the dashed lines yields
 \begin{multline*}\{(h,h_0,\dots, h_{2r-1};C_0,C_1,\dots, C_{2r})\in H^{2r+1}\times G^{2r+1}\\;\  h_{2i+1}^{-1}C_{2i+1} C_{2i}^{-1}h_{2i}\in U_+\text{ and }h_{2i}^{-1}C_{2i} C_{2i-1}^{-1}h_{2i-1}\in U_-,\}.\end{multline*}
 Meanwhile, since $B_+\cap B_-=G$, the group of residual gauge transformations is $\prod_{x_1,\dots, x_{2r}}H$, acting at the appropriate points on the quilted surface. 
 Thus, up to a gauge transformation, we may assume that $h_0=h_1=\dots=h_{2r-1}=\mbf{1}$. Setting $S_i=C_{i} C_{i-1}^{-1}$, we see that
 that the moduli space
can be identified with 
$$_G\mc{A}_H^r:=\{(h;S_{2r},\dots,S_1;C_0)\in H\times (U_-\times U_+)^r\times G\}.$$

As explained in Remark~\ref{rmk:qPoisQuilt}, $_G\mc{A}_H^r$ is a quasi-Poisson $G\times H$ manifold, where $g\in G$ and $k\in H$ act at the marked points $v_G$ and $v_H$, respectively:
$$(g,k)\cdot(h;S_{2r},\dots,S_1;C_0)=(khk^{-1},kS_{2r}k^{-1},\dots, kS_{1}k^{-1},kC_0g^{-1}),$$
The holonomy along the boundary components,
$$(h;S_{2r},\dots,S_1;C_0)\to (C_0^{-1}hS_{2r}\cdots S_1C_0,h^{-1}),$$
defines a moment map on $_G\mc{A}_H^r$.

This quasi-Hamiltonian $G\times H$-space was first discovered by Boalch \cite{Boalch:2011vt,Boalch:2007ty,Boalch:2009tn}, who used it to study  meromorphic connections on Riemann surfaces.

\end{example}

\begin{example}[Bott-Samelson variety \cite{Balazs:2013uk}]
Suppose that $s_\g\in S^2(\g)^\g$ is non-degenerate, $\g=\mf{u}_+\oplus\h\oplus\mf{u}_-$ as a vector space (but not as a Lie algebra), where $\mf{b}_\pm:=\mf{u}_\pm\oplus\h\subseteq \g$ are coisotropic subalgebras satisfying $\mf{b}_\pm^\perp=\mf{u}_\pm$. Let $s_\h\in S^2(\h)^\h$ be the projection of $s_\g$ along $$\g=\h\oplus(\mf{u}_+\oplus\mf{u}_-)\to\h.$$ Suppose further that the Lie subalgebras $\mf{b}_\pm,\h$ integrate to closed subgroups $B_\pm,H\subseteq G$ such that $H=B_+\cap B_-$.

We take $\tilde \g=\g\oplus\h$, and define 
$$\hat{\mf{b}}_\pm=\{(\xi\pm\eta,\eta)\mid \xi\in\mf{u}_\pm\text{ and }\eta\in \h\}\subseteq \tilde\g$$
Then $(\tilde\g,\hat{\mf{b}}_+,\hat{\mf{b}}_-)$ forms a Manin triple with respect to the non-degenerate copairing $s:=s_\g\oplus-s_\h \in S^2(\tilde\g)^{\tilde\g}$.

Suppose that the Lie subalgebras $\hat{\mf{b}}_\pm$ integrate to closed subgroups $\hat B_\pm\subseteq \tilde G:=G\times H$ which project isomorphically to $B_\pm$ along the map $G\times H\to G$; and suppose for simplicity t that $\hat B_+\cap \hat B_-=\mbf{1}$.

Let $(\tilde \g_n,\mf{b}^+_n,\mf{b}^-_n)$ denote the Manin triple corresponding to the $n$-th Drinfel'd Polyuble, $\mc{P}_n(\hat B_+)$, as described in Example~\ref{ex:poly}, and let $(\tilde G_n, B^+_n,B^-_n)$ denote the corresponding triple of Lie groups (cf. Figure~\ref{fig:Polyuble2}). 

Notice that $\tilde B_+=B_+\times H\subseteq \tilde G$ is a coisotropic subgroup.
We consider the Lu-Yakimov Poisson homogeneous space, $\mc{M}_{big}=\tilde G^n/\tilde B_+^n\cong G^n/B_+^n$, seen as the moduli space for the following quilted surface with $\tilde G^n$-colored domain (cf. Example~\ref{ex:luYak}):
\begin{center}
\begingroup%
  \makeatletter%
  \providecommand\color[2][]{%
    \errmessage{(Inkscape) Color is used for the text in Inkscape, but the package 'color.sty' is not loaded}%
    \renewcommand\color[2][]{}%
  }%
  \providecommand\transparent[1]{%
    \errmessage{(Inkscape) Transparency is used (non-zero) for the text in Inkscape, but the package 'transparent.sty' is not loaded}%
    \renewcommand\transparent[1]{}%
  }%
  \providecommand\rotatebox[2]{#2}%
  \ifx\svgwidth\undefined%
    \setlength{\unitlength}{128.29272461bp}%
    \ifx\svgscale\undefined%
      \relax%
    \else%
      \setlength{\unitlength}{\unitlength * \real{\svgscale}}%
    \fi%
  \else%
    \setlength{\unitlength}{\svgwidth}%
  \fi%
  \global\let\svgwidth\undefined%
  \global\let\svgscale\undefined%
  \makeatother%
  \begin{picture}(1,0.64531439)%
    \put(0,0){\includegraphics[width=\unitlength]{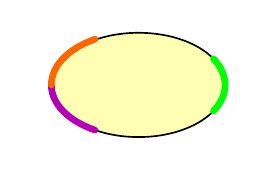}}%
    \put(0.22293146,0.460866){\color[rgb]{0,0,0}\makebox(0,0)[rb]{\smash{$B_n^-$}}}%
    \put(0.22293146,0.16774817){\color[rgb]{0,0,0}\makebox(0,0)[rb]{\smash{$B_n^+$}}}%
    \put(0.86851013,0.30194668){\color[rgb]{0,0,0}\makebox(0,0)[lb]{\smash{$\tilde B_+^n$}}}%
    \put(0.51395027,0.30823214){\color[rgb]{0,0,0}\makebox(0,0)[b]{\smash{$\tilde G^n$}}}%
  \end{picture}%
\endgroup%

\end{center}

Recall that a $\tilde G^n$-colored domain $\Sigma$ is equivalent to $n$ copies of $\Sigma$, colored by $\tilde G$. `Unfolding' the $n$-layered quilted surface pictured above (as in Example~\ref{ex:poly}), we obtain the following $\tilde G$-colored quilted surface with $n$ $\tilde B_+$-colored domain walls:
\noindent%
\begin{minipage}{\linewidth}
\makebox[\linewidth]{
\begingroup%
  \makeatletter%
  \providecommand\color[2][]{%
    \errmessage{(Inkscape) Color is used for the text in Inkscape, but the package 'color.sty' is not loaded}%
    \renewcommand\color[2][]{}%
  }%
  \providecommand\transparent[1]{%
    \errmessage{(Inkscape) Transparency is used (non-zero) for the text in Inkscape, but the package 'transparent.sty' is not loaded}%
    \renewcommand\transparent[1]{}%
  }%
  \providecommand\rotatebox[2]{#2}%
  \ifx\svgwidth\undefined%
    \setlength{\unitlength}{194.72602539bp}%
    \ifx\svgscale\undefined%
      \relax%
    \else%
      \setlength{\unitlength}{\unitlength * \real{\svgscale}}%
    \fi%
  \else%
    \setlength{\unitlength}{\svgwidth}%
  \fi%
  \global\let\svgwidth\undefined%
  \global\let\svgscale\undefined%
  \makeatother%
  \begin{picture}(1,0.7208405)%
    \put(0,0){\includegraphics[width=\unitlength]{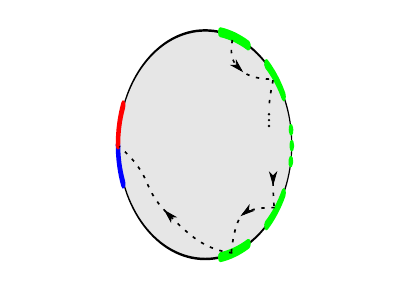}}%
    \put(0.27829564,0.42732594){\color[rgb]{0,0,0}\makebox(0,0)[rb]{\smash{$\hat B_-$}}}%
    \put(0.27829564,0.28350909){\color[rgb]{0,0,0}\makebox(0,0)[rb]{\smash{$\hat B_+$}}}%
    \put(0.59817403,0.64814221){\color[rgb]{0,0,0}\makebox(0,0)[lb]{\smash{$\tilde B_+$}}}%
    \put(0.69981862,0.53734953){\color[rgb]{0,0,0}\makebox(0,0)[lb]{\smash{$\tilde B_+$}}}%
    \put(0.59817403,0.04832511){\color[rgb]{0,0,0}\makebox(0,0)[lb]{\smash{$\tilde B_+$}}}%
    \put(0.69981862,0.18403261){\color[rgb]{0,0,0}\makebox(0,0)[lb]{\smash{$\tilde B_+$}}}%
    \put(0.41905042,0.21855373){\color[rgb]{0,0,0}\makebox(0,0)[lb]{\smash{$g_1$}}}%
    \put(0.58524358,0.20656921){\color[rgb]{0,0,0}\makebox(0,0)[rb]{\smash{$g_2$}}}%
    \put(0.5651484,0.53562776){\color[rgb]{0,0,0}\makebox(0,0)[rb]{\smash{$g_n$}}}%
  \end{picture}%
\endgroup%
}
\captionof{figure}{}\label{fig:LuYakimovPolyU}
\end{minipage}
Computing the holonomies along the dotted arcs identifies the corresponding moduli space with $$\mc{M}_{big}:=\overset{n-1}{\overbrace{\tilde G\times_{\tilde B_+}\cdots\times_{\tilde B_+}\tilde G}}\times_{\tilde B_+}\tilde G/\tilde B_+.$$
As in Example~\ref{ex:luYak}, $\mc{M}_{big}$ has a Poisson action of the Drinfel'd Polyuble $\mc{P}_{2n}(\hat B_+)$.


We now take a closer look at this moduli space. Let $\Sigma'$ denote the disk with $n+2$ marked points $V=\{-1,0,\dots,n\}$: 
\begin{center}
\begingroup%
  \makeatletter%
  \providecommand\color[2][]{%
    \errmessage{(Inkscape) Color is used for the text in Inkscape, but the package 'color.sty' is not loaded}%
    \renewcommand\color[2][]{}%
  }%
  \providecommand\transparent[1]{%
    \errmessage{(Inkscape) Transparency is used (non-zero) for the text in Inkscape, but the package 'transparent.sty' is not loaded}%
    \renewcommand\transparent[1]{}%
  }%
  \providecommand\rotatebox[2]{#2}%
  \ifx\svgwidth\undefined%
    \setlength{\unitlength}{123.8197998bp}%
    \ifx\svgscale\undefined%
      \relax%
    \else%
      \setlength{\unitlength}{\unitlength * \real{\svgscale}}%
    \fi%
  \else%
    \setlength{\unitlength}{\svgwidth}%
  \fi%
  \global\let\svgwidth\undefined%
  \global\let\svgscale\undefined%
  \makeatother%
  \begin{picture}(1,0.7730014)%
    \put(0,0){\includegraphics[width=\unitlength]{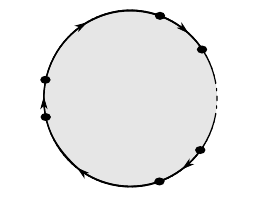}}%
    \put(0.14792792,0.46747819){\color[rgb]{0,0,0}\makebox(0,0)[rb]{\smash{$-1$}}}%
    \put(0.14792792,0.29341338){\color[rgb]{0,0,0}\makebox(0,0)[rb]{\smash{$0$}}}%
    \put(0.6509869,0.73473713){\color[rgb]{0,0,0}\makebox(0,0)[lb]{\smash{$n$}}}%
    \put(0.81083904,0.60064225){\color[rgb]{0,0,0}\makebox(0,0)[lb]{\smash{$n-1$}}}%
    \put(0.6509869,0.00876486){\color[rgb]{0,0,0}\makebox(0,0)[lb]{\smash{$1$}}}%
    \put(0.81083904,0.17301476){\color[rgb]{0,0,0}\makebox(0,0)[lb]{\smash{$2$}}}%
    \put(0.32195183,0.13563896){\color[rgb]{0,0,0}\makebox(0,0)[lb]{\smash{$p_1$}}}%
    \put(0.19083419,0.3820499){\color[rgb]{0,0,0}\makebox(0,0)[lb]{\smash{$p_0$}}}%
    \put(0.27771582,0.59075567){\color[rgb]{0,0,0}\makebox(0,0)[lb]{\smash{$p_{-1}$}}}%
    \put(0.70723937,0.61562283){\color[rgb]{0,0,0}\makebox(0,0)[rb]{\smash{$p_n$}}}%
    \put(0.70271797,0.13862556){\color[rgb]{0,0,0}\makebox(0,0)[rb]{\smash{$p_2$}}}%
  \end{picture}%
\endgroup%

\end{center}
Computing holonomies along the boundary arcs identifies 
$$M_{\Sigma',V}(\tilde G)\cong\{(p_0,\dots,p_n)\in G^{n+1}\},$$
a quasi-Poisson $\tilde G^{n+2}$ space with the $\tilde G^{n+2}$ action given by 
$$(g_{-1}',g_0',\dots,g_n')\cdot(p_0,\dots,p_n)=\big(g_{-1}'p_0{g_0'}^{-1},\dots, g_{n-1}'p_n{g_n'}^{-1}\big),$$
and moment map $\mu:M_{\Sigma',V}(\tilde G)\to \tilde G^{n+2}$ given by 
$$\mu(p_0,\dots,p_n)=(p_{-1},p_0,\dots,p_n)\quad (\text{where } p_{-1}:=(p_0\cdots p_n)^{-1}).$$
Let $$S_{big}=\tilde G\times(\hat B_-\cdot \hat B_+)\times\tilde G^n\subseteq \tilde G^{n+2},$$
and 
$$C=\hat B_-\times \hat B_+\times \tilde B_+^n\subseteq \tilde G^{n+2}.$$
The moduli space $\mc{M}_{big}$ for the quilted surface picture in Figure~\ref{fig:LuYakimovPolyU} is the moment map quotient 
$$\mc{M}_{big}=\mu^{-1}(S_{big})/C$$

Now choose subgroups $B_+\subseteq P_1,\dots,P_n\subseteq G$ between $B_+$ and the ambient group $G$. With $\tilde P_i:=P_i\times H\subseteq \tilde G$, we define 
$$S_{small}=\tilde G\times(\hat B_-\cdot \hat B_+)\times \tilde P_1\times\cdots\times \tilde P_n\subseteq S_{big}$$
(note that $S_{small}$ is a $\hat C$-stable submanifold of $\tilde G^{n+2}$).
Computing the holonomies\footnote{More specifically $$\mu^{-1}(S_{small})=\{(p_0,p_1,\dots,p_n)\mid  p_0=b_-b_+\in \hat B_-\cdot \hat B_+,\text{ and } p_i\in \tilde P_i\}$$ and $(\hat b_-,\hat b_+,\tilde b_1,\dots, \tilde b_n)\in C$ acts by 
$$(\hat b_-,\hat b_+,\tilde b_1,\dots, \tilde b_n)\cdot(b_-b_+;p_1,\dots,p_n)=(\hat b_-b_-b_+\hat b_+^{-1};\hat b_+p_1\tilde b_1^{-1},\tilde b_1 p_2\tilde b_2^{-1},\dots, \tilde b_{n-1}p_n\tilde b_{n}^{-1}).$$
Thus the map (which computes the holonomies pictured in Figure~\ref{fig:LuYakimovPolyU})
$$(b_-b_+,p_1,\dots,p_n)\to [g_1,\dots,g_n]:= [b_+p_1,p_2,p_3,\dots,p_n]\in \tilde P_1\times_{\tilde B_+}\cdots \times_{\tilde B_+}\tilde P_n/\tilde B_+$$ 
descends to an isomorphism on $\mu^{-1}(S_{small})/C$.}
  along the dotted arcs in Figure~\ref{fig:LuYakimovPolyU} identifies the moduli space $\mc{M}_{small}:=\mu^{-1}(S_{small})/C$ with the Bott-Samelson variety 
$$P_1\times_{ B_+}P_2\times_{B_+}\cdots \times_{ B_+} P_n/ B_+\cong\tilde P_1\times_{\tilde B_+}\cdots \times_{\tilde B_+}\tilde P_n/\tilde B_+\subseteq \mc{M}_{big}.$$

In particular, we recover the result of Jiang-Hua Lu \cite{Lu:FBHentwA} that the Bott-Samelson variety $\mc{M}_{small}$ is a Poisson submanifold of $\mc{M}_{big}$.

There are $n$ and distinct Poisson morphisms from the Bott-Samelson variety $\mc{M}_{small}$ to the Shubert variety $G/B_+\cong \tilde G/\tilde B_+$ (which carries the Lu-Yakimov Poisson structure, cf. Example~\ref{ex:luYak}),
$$q_k:[g_1,\dots,g_n]\to [g_1\cdots g_k],\quad 1\leq k\leq n,$$
corresponding to the embeddings (we picture the case $k=2$):
\begin{center}
   \def\svgwidth{.5\linewidth}
\begingroup%
  \makeatletter%
  \providecommand\color[2][]{%
    \errmessage{(Inkscape) Color is used for the text in Inkscape, but the package 'color.sty' is not loaded}%
    \renewcommand\color[2][]{}%
  }%
  \providecommand\transparent[1]{%
    \errmessage{(Inkscape) Transparency is used (non-zero) for the text in Inkscape, but the package 'transparent.sty' is not loaded}%
    \renewcommand\transparent[1]{}%
  }%
  \providecommand\rotatebox[2]{#2}%
  \ifx\svgwidth\undefined%
    \setlength{\unitlength}{240.17102051bp}%
    \ifx\svgscale\undefined%
      \relax%
    \else%
      \setlength{\unitlength}{\unitlength * \real{\svgscale}}%
    \fi%
  \else%
    \setlength{\unitlength}{\svgwidth}%
  \fi%
  \global\let\svgwidth\undefined%
  \global\let\svgscale\undefined%
  \makeatother%
  \begin{picture}(1,0.52986723)%
    \put(0,0){\includegraphics[width=\unitlength]{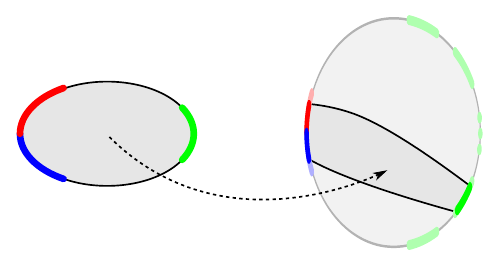}}%
  \end{picture}%
\endgroup%

\end{center}

\end{example}

\appendix
\section{Quilted surfaces and moduli spaces of flat bundles}
In this section we prove Theorems~\ref{thm:FullModSpc}~and~\ref{thm:CompCornData}, and identify the more general Poisson manifolds $\hat\mu^{-1}(h)/\on{Stab}(h)\subseteq M_{\Sigma,V}/C$ with moduli spaces.

\begin{proof}[Proof of Theorem~\ref{thm:FullModSpc}]
First, we recall that $M_{\Sigma,V}$ is the moduli space of pairs \eqref{eq:FramedBndls}
$$(P\to\Sigma; \{\hat v\in P_v\}_{v\in V}),$$ 
of flat $G$-bundles\footnote{Recall that we always take our $G$-bundles to be principal.} over $\Sigma$ equipped with a framing over $V\subset \partial \Sigma$. Of course this moduli space only depends on the homotopy type of the pair $(\Sigma, V)$. Since $(\Sigma, V)$ has the same homotopy type as $(\hat\Sigma,\cup_{v\in V}\mbf{w}_v)$, $M_{\Sigma,V}$ is equivalently described as the moduli space of pairs
$$(P\to\hat\Sigma; \{\hat{\mbf{w}}\in C^\infty_{flat}(\mbf{w}, P\rvert_{\mbf{w}})\}_{\mbf{w}\in\mbf{Wall}}),$$ 
of flat $G$-bundles over $\hat\Sigma$ which are equipped with a flat framing over each domain wall $\mbf{w}\in\mbf{Wall}$. 

Likewise, $M_{\Sigma,V}/C$ is canonically isomorphic to the moduli space of pairs
$$(P\to \Sigma, \{Q_v\to v\}_{v\in V})$$
of flat $G$-bundles over $\Sigma$ which are equipped with reduction of the structure group from $G$ to $C_v$ over each marked point $v\in V$.\footnote{That is, a choice of a principal $C_v$-subbundle $Q_v\subseteq P_v$ over each marked point $v\in V$.} As above, we can also canonically identify $M_{\Sigma,V}/C$ with the moduli space of pairs
$$(P\to\hat\Sigma,\{Q_{\mbf{w}}\to \mbf{w}\}_{\mbf{w}\in\mbf{Wall}})$$
of flat $G$-bundles over $\hat\Sigma$ equipped with a flat reduction of the structure group from $G$ to $C_\mbf{w}$ over each domain wall $\mbf{w}\in\mbf{Wall}$
 (as in the description we gave in the introduction). Indeed the isomorphism between these two moduli spaces is given by sending $(P\to \Sigma, \{Q_v\to v\}_{v\in V})$ to its pullback $(p^*P\to \hat\Sigma, \{p^*Q_v\to \mbf{w}_v\}_{v\in V})$ along the blowdown map $p:\hat\Sigma\to \Sigma$.
 Thus we have proven Theorem~\ref{thm:FullModSpc}.
 \end{proof}

 \begin{proof}[Proof of Theorem~\ref{thm:CompCornData}]
 To describe the Poisson submanifolds $\hat\mu^{-1}(\mbf{1})/\on{Stab}(\mbf{1})\subseteq M_{\Sigma,V}/C$ as moduli spaces, we need to first identify $\hat M$ with a moduli space. By assumption, $s\in S^2(\g)$ is non-degenerate, and each $C_v\subseteq G$ is Lagrangian (i.e. $C^\perp_v= C_v$ for each $v\in V$). Thus $\hat C_v=C_v\times C_v\subseteq G\times G$. (We will treat the general case later).

As a first step, consider the moduli space of triples 
\begin{equation}\label{eq:QuadMod}(P\to \hat\Sigma, \{Q_{\mbf{w}}\to \mbf{w}\}_{\mbf{w}\in \mbf{Wall}}, \{x_v\in Q_{\mbf{w}_v}\rvert_{v}\}_{v\in V_-\cup V_+}),\end{equation}
i.e a pair of the form \eqref{eq:PairMod}, but equipped with an additional framing of the $C_v$-bundle $Q_{\mbf{w}_v}$ at either end point $v_-,v_+\subset \mbf{w}_v$. Since $x_v\in Q_{\mbf{w}_v}\subset P$, we have a fully faithful forgetful functor which sends a triple \eqref{eq:QuadMod} to the pair $((P\to \hat\Sigma, \{x_v\in P\rvert_{v}\}_{v\in V_-\cup  V_+})$. Thus, the moduli space of triples \eqref{eq:QuadMod} imbeds as a subspace of the moduli space $M_{\hat\Sigma,V_-\cup V_+}$ of pairs \eqref{eq:FramedBndls}.
  Now the holonomy of a given connection along the domain wall $\mbf{w}_v$, as measured with respect to the framing $x_{v_-},x_{ v_+}$, is an element of $C_v$. This is the only restriction placed on the holonomies by the structure of the triples \eqref{eq:QuadMod}, and it identifies the moduli space of triples with the imbedded submanifold $\hat M\subseteq M_{\hat\Sigma,V_-\cup V_+}$ (cf. \eqref{eq:hatM}).

The forgetful functor sending a triple \eqref{eq:QuadMod} to the corresponding pair \eqref{eq:PairMod} coincides with the quotient map $\hat M\to \hat M/\hat C\cong M_{\Sigma,V}/C$. Meanwhile,  $\hat\mu^{-1}(\mbf{1})\subset\hat M$ consists of those triples \eqref{eq:QuadMod} for which parallel transport along $a_v$ takes $x_{v_+}$ to $x_{\sigma(v)_-}$.


Let $q:\hat\Sigma\to \hat\Sigma^c$ denote the quotient map which collapses each respective segment $a_v$ to a point, as in Figure~\ref{fig:tilSig}. We denote the image of each domain wall $\mbf{w}$ in $\hat\Sigma^c$ by the same symbol. 
The flat connection along each $a_v$ canonically trivializes $P\rvert_{a_v}$, allowing us to identify $P\to \hat\Sigma$ with the pullback of a flat bundle $P^c\to \hat\Sigma^c$, i.e. $P=q^* P^c$. Moreover, if the parallel transport along $a_v$ takes $x_{v_+}$ to $x_{\sigma(v)_-}$, then $q(x_{v_+})=q(x_{\sigma(v)_-})\in (Q_{\mbf{w}_v}\cap Q_{\mbf{w}_{\sigma(v)}})\rvert_{q(v)}$ is a common framing. Thus the moduli space $\hat\mu^{-1}(\mbf{1})$ of triples \eqref{eq:QuadMod} is canonically isomorphic to the moduli space of triples over $\hat \Sigma^c$:
\begin{equation}\label{eq:TripSigc}(P^c\to \hat\Sigma^c, \{Q_{\mbf{w}}\to \mbf{w}\}_{\mbf{w}\in \mbf{Wall}}, \{x_v\in (Q_{\mbf{w}_v}\cap Q_{\mbf{w}_{\sigma(v)}})\rvert_{q(v)}\}_{v\in V})\end{equation}
where each $x_v$ is a common framing of $Q_{\mbf{w}_v}$ and $Q_{\mbf{w}_{\sigma(v)}}$ over $q(v)$.

Now $\hat C$ acts on the moduli space $\hat M$ of triples \eqref{eq:TripSigc} by forgetting the framing. The essential image of this forgetful functor consists of those pairs
$$(P^c\to\hat\Sigma^c,\{Q_{\mbf{w}}\to \mbf{w}\}_{\mbf{w}\in \mbf{Wall}})$$
of flat $G$-bundles over $\hat\Sigma^c$ equipped with reductions of the structure group from $G$ to $C_\mbf{w}$ over each domain wall $\mbf{w}$ such that $\emptyset\neq(Q_{\mbf{w}_v}\cap Q_{\mbf{w}_{\sigma(v)}})\rvert_{q(v)}$. Equivalently, for any two domain walls $\mbf{w},\mbf{w}'\in\mbf{Wall}$, the fibres of $Q_{\mbf{w}}$ and $Q_{\mbf{w}'}$ intersect non-trivially over every point in the intersection of their bases $\mbf{w}\cap\mbf{w}'$. That is, the reductions of structure along $\mbf{w}$ and $\mbf{w}'$ are compatible (cf. Definition~\ref{def:CompRed}).
This proves Theorem~\ref{thm:CompCornData}
\end{proof}


We now lift the restrictions that $s\in S^2(\g)$ be non-degenerate, and each $C_v\subseteq G$ be Lagrangian; instead we suppose only that $C^\perp_v\subset C_v$ is closed, for each $v\in V$.
To describe the Poisson submanifolds $\hat\mu^{-1}(h)/\on{Stab}(h)\subseteq M_{\Sigma,V}/C$ as moduli spaces, we proceed as above, by first identifying $\hat M$ with a moduli space.

Consider the moduli space of triples 
\begin{equation}\label{eq:QuadModGen}(P\to \hat\Sigma, \{Q_{\mbf{w}}\to \mbf{w}\}_{\mbf{w}\in \mbf{Wall}}, \{x_v\in Q_{\mbf{w}_v}\rvert_{v}\}_{v\in V_-\cup V_+}),\end{equation}
i.e a pair of the form \eqref{eq:PairMod}, but equipped with an additional framing of the $C_v$-bundle $Q_{\mbf{w}_v}$ at either end point $v_-,v_+\subset \mbf{w}_v$ such that parallel transport in the quotient bundle $Q_{\mbf{w}_v}/C_{\mbf{w}_v}^\perp$ identifies their images $\tilde x_{v_-}, \tilde x_{v_+}\in Q_{\mbf{w}_v}/C_{\mbf{w}_v}^\perp$. As in the proof of Theorem~\ref{thm:CompCornData}, the forgetful functor which sends a triple \eqref{eq:QuadModGen} to the pair $$((P\to \hat\Sigma, \{x_v\in P\rvert_{v}\}_{v\in V_-\cup V_+})$$ identifies the moduli space of triples with the imbedded submanifold $\hat M\subseteq M_{\hat\Sigma,V_-\cup V_+}$ (cf. \eqref{eq:hatM}).


The forgetful functor sending a triple \eqref{eq:QuadModGen} to the corresponding pair \eqref{eq:PairMod} coincides with the quotient map $\hat M\to \hat M/\hat C\cong M_{\Sigma,V}/C$; and the moment map $\hat\mu:\hat M\to G^V$ coincides with the map which sends an isomorphism class of triples \eqref{eq:QuadModGen} to the collection of holonomies $(\hol_{a_v})_{v\in V}\in G^V$ (these holonomies are measured with respect to the chosen framings).

Let $q:\hat\Sigma\to \hat\Sigma^c$ denote the quotient map which collapses each respective segment $a_v$ to a point, as in Figure~\ref{fig:tilSig}. 
 As before, $P\to \hat\Sigma$ is pullback of a flat bundle $P^c\to \hat\Sigma^c$. Thus the moduli space $\hat M$ of triples \eqref{eq:QuadModGen} is canonically isomorphic to the moduli space of triples over $\hat \Sigma^c$:
$$(P^c\to \hat\Sigma^c, \{Q_{\mbf{w}}\to \mbf{w}\}_{\mbf{w}\in \mbf{Wall}}, \{x_v\in Q_{\mbf{w}_v}\rvert_{q(v)}\}_{v\in V_-\cup  V_+}),$$
and the moment map $\hat\mu$ sends a given triple to the element $(g_v)_{v\in V}\in G^V$ such that $x_{\sigma(v)_-}=g_v\cdot x_{v_+}$ for each $v\in V$. In particular, the moment map factors through the isomorphism classes of triples over the boundary:
\begin{equation}\label{eq:BndQuadMod}(P'\to \partial\hat\Sigma^c, \{Q_{\mbf{w}}'\to \mbf{w}\}_{\mbf{w}\in \mbf{Wall}}, \{x_v\in Q_{\mbf{w}_v}'\rvert_{q(v)}\}_{v\in V_-\cup V_+}),\end{equation}
where $P'\to \partial\hat\Sigma^c$ is a $G$-bundle, $Q_{\mbf{w}}'\subseteq P'\rvert_{\mbf{w}}$ is a $C_\mbf{w}$-subbundle equipped with
\begin{enumerate}
\item a flat connection on the quotient bundle $Q_{\mbf{w}}/C_{\mbf{w}}^\perp$, and 
\item a framing at either endpoint of $\mbf{w}$ such that parallel transport in the quotient bundle, $Q_{\mbf{w}}/C_{\mbf{w}}^\perp$, identifies the images of the framings.
\end{enumerate}
 An element $(g_{v_-},g_{v_+})_{v\in V}\in \hat C\subseteq G^{V_-\cup  V_+}$ acts on a triple \eqref{eq:BndQuadMod} by $x_v\to g_v\cdot x_v$ for each $v\in V_-\cup V_+$.

Notice that neither $P'$ nor $Q_{\mbf{w}}'$ carry a distinguished flat structure, so the $\hat C$-equivariant map which sends an isomorphism class of such quadruples \eqref{eq:BndQuadMod} to the unique element $(g_v)_{v\in V}\in G^V$ such that $x_{\sigma(v)_-}=g_v\cdot x_{v_+}$ for each $v\in V$ is injective. In particular, a given isomorphism class of quadruples \eqref{eq:BndQuadMod} is naturally identified with an element of $G^V$, and under this identification, the moment map $\hat\mu:\hat M\to G^V$ coincides with the functor which sends a quadruple of the form \eqref{eq:QuadModGen} to a quadruple of the form \eqref{eq:BndQuadMod} by restricting the ambient principal bundle $P^c$ to $\partial\hat\Sigma^c$, and forgetting the superfluous flat connections.

\begin{defn}
\emph{Boundary data} for the surface $\hat\Sigma^c$ is a pair
$$(P\to \partial\hat\Sigma^c, \{Q_{\mbf{w}}\to \mbf{w}\}_{\mbf{w}\in \mbf{Wall}}),$$
where $P\to \partial\hat\Sigma^c$ is a $G$-bundle, and each $Q_{\mbf{w}}\subseteq P\rvert_{\mbf{w}}$ is a $C_\mbf{w}$-subbundle equipped with a choice of flat connection on the quotient bundle $Q_{\mbf{w}}/C_\mbf{w}^\perp$.
\end{defn}

Thus the isomorphism classes of boundary data for $\hat\Sigma^c$ correspond precisely to the $\hat C$-orbits of isomorphism classes of quadruples \eqref{eq:BndQuadMod}. In particular, a given isomorphism class of boundary data is naturally identified with a $\hat C$-orbit $\mc{O}\subseteq G^V$. Of course, given a flat $G$-bundle over $\hat\Sigma^c$ equipped with a flat reduction of the structure group from $G$ to $C_\mbf{w}$ over each domain wall $\{\mbf{w}\}_{\mbf{w}\in \mbf{Wall}}$:
$$(P\to\hat\Sigma^c,\{Q_{\mbf{w}}\to \mbf{w}\}_{\mbf{w}\in \mbf{Wall}}),$$
restricting $P$ to $\partial \hat\Sigma^c$, and forgetting the superfluous flat connections defines a boundary data,
$$(P\rvert_{\partial\hat\Sigma^c},\{Q_{\mbf{w}}\}_{\mbf{w}\in \mbf{Wall}}),$$
 for $\hat\Sigma^c$.
Summarizing, we have:

\begin{thm}
Let $\mc{O}\subseteq G^V$ be a $\hat C$-orbit corresponding to a given isomorphism class of boundary data for $\hat\Sigma^c$.
Suppose that $\mc{O}$  is in a clean position relative to $\hat\mu$, and the $\hat C$-orbits of $\hat\mu^{-1}(\mc{O})$ form a regular foliation. Then the moduli space of pairs 
$$(P\to\hat\Sigma^c,\{Q_{\mbf{w}}\to \mbf{w}\}_{\mbf{w}\in \mbf{Wall}}),$$
a flat $G$-bundles over $\hat\Sigma^c$ equipped with a flat reduction of the structure group from $G$ to $C_\mbf{w}$ over each domain wall $\{\mbf{w}\}_{\mbf{w}\in \mbf{Wall}}$ and such that $(P\rvert_{\partial\hat\Sigma^c},\{Q_{\mbf{w}}\}_{\mbf{w}\in \mbf{Wall}}),$ lies in the given isomorphism class of boundary data, is a Poisson manifold which is canonically isomorphic to $\hat\mu^{-1}(\mc{O})/\hat C$ as a Poisson manifold.
\end{thm}
\begin{rem}
When $s\in S^2(\g)$ is non-degenerate, and each $C_v\subseteq G$ is Lagrangian (i.e. $C^\perp_v= C_v$ for each $v\in V$),
each moduli space described in Theorem~\ref{thm:CompCornData} is a symplectic manifold (cf. \cite{LiBland:2013ue}).
\end{rem}

\section{Some Technical Lemmas}
We expect the following lemmas are common knowledge, but we were unable to find any proofs for them, so we have included these proofs for completeness.
\begin{lem}\label{lem:CleanInter}
Suppose that embedded submanifolds $U,V\subseteq W$ intersect cleanly, $q:W\to Q$ is a surjective submersion, and $V$ is a union of $q$-fibres, i.e $V=q^{-1}\big(q(V)\big)$.
If $q(U),q(V)\subseteq Q$ are embedded submanifolds, and $q\rvert_{U}:U\to q(U)$ is a submersion, then $q(U)$ and $q(V)$ intersect cleanly.
\end{lem}
\begin{proof}We claim that $q(U)\cap q(V)\subseteq Q$ is an embedded manifold.
Assuming this, we need only show that $$T\big(q(U)\cap q(V)\big)=Tq(U)\cap Tq(V).$$ It is clear that the right hand side contains the left hand side. Suppose that $$\xi \in Tq(U)\cap Tq(V)=q_*(TU)\cap q_*(TV).$$ By assumption, $TV$ contains the kernel of $q_*$, so there exists a lift $\tilde \xi\in TU\cap TV$ such that $\xi=q_*(\tilde \xi)$. Since $U$ and $V$ intersect cleanly, $\tilde\xi\in T(U\cap V)$, and hence $\xi=q_*(\tilde\xi)\in q_*\big( T(U\cap V)\big)\subseteq T\big(q(U)\cap q(V)\big)$.

We now prove the claim that $q(U)\cap q(V)\subseteq Q$ is an embedded submanifold. Consider the groupoid $U\times_{q(U)}U$ over $U$ corresponding to the surjective submersion $q:U\to q(U)$. Then $U\times_{q(U)}U$ is Morita equivalent to $q(U)$, thus it is a proper groupoid with trivial isotropy groups \cite[\S~5.6]{Moerdijk:2003vw}. Since $U\cap V\subseteq U$ is an embedded submanifold, it follows that 
 $$U\times_{q(U)}(U\cap V)\subseteq U\times_{q(U)}U$$ is an embedded submanifold. 
 Now $V$ is a union of $q$-fibres, therefore $U\cap V\subseteq U$ is a union of $q\rvert_U$-fibres, and $(U\cap V)\times_{q(U)}(U\cap V)=U\times_{q(U)}(U\cap V)$.
 
In summary, we have shown that $$H:=(U\cap V)\times_{q(U)}(U\cap V) \subseteq U\times_{q(U)}U$$ is an embedded subgroupoid over base $U\cap V$. Since $H$ is a proper Lie groupoid with trivial isotropy groups (it inherits these properties from the inclusion $H\subseteq U\times_{q(U)}U$), it follows that $H$ is Morita equivalent to a manifold \cite[\S~5.6]{Moerdijk:2003vw}. That is to say, the orbit space of $H$, $q(U\cap V)$, is a manifold.
Since $V$ is a union of $q$-fibres, $q(U)\cap q(V)=q(U\cap V)$ is also a manifold. 

By construction, the inclusion $q(U)\cap q(V)\to Q$ is an injective immersion. Now $U\cap V\subset M$ has the subset topology, while $q(U\cap V)$ and $Q$ both have the quotient topology. Since $q$ is a surjective submersion, $q(U)\cap q(V)=q(U\cap V)\subseteq Q$ has the subset topology. Thus $q(U)\cap q(V)\subseteq Q$ is an embedded submanifold.

\end{proof}

\begin{lem}\label{lem:MoritaEqReg}
Suppose that $G$ and $H$ are Morita equivalent Lie groupoids, and the orbit-equivalence relation on $H$ is regular, then the orbit-equivalence relation on $G$ is regular too.
\end{lem}
\begin{proof}By Godement's criterion, the image $R_H=\mbf{t}\times\mbf{s}(H_1)\subseteq H_0\times H_0$ is a closed embedded submanifold.
Let $G\leftarrow K\rightarrow H$ be a Morita equivalence.\footnote{i.e. $K$ is a Lie groupoid, and the two maps are Morita morphisms (cf. \cite{Moerdijk:2003vw,Behrend06}).} Then $K_0\to H_0$ is a surjective submersion, and the pullback $R_K=R_H\times_{H_0\times H_0}(K_0\times K_0)$ is a closed embedded submanifold (since $K_0\to H_0$ is a surjective submersion). But $R_K=\mbf{t}\times\mbf{s}(K_1)$, so the projection from from $R_K$ onto either factor of $K_0$ is a surjective submersion. It follows from Godement's criterion that $R_K$ is a regular equivalence relation.

Now form $Q=K_0/R_K$. The quotient map $q:K_0\to Q$ is a surjective submersion, which factors through the surjective submersion $f:K_0\to G_0$, i.e. there exists a map $p:G_0\to Q$ such that $q=p\circ f$. Since $f$ and $q$ are smooth, surjective submersions, it follows that $q$ is a smooth, surjective submersion. Since $Q=G_0/G_1$ is the orbit space for $G$, this proves the lemma.
\end{proof}


\Addresses
\end{document}